# Generalized Conditional Gradient for Sparse Estimation

**Yaoliang Yu**                                    YAOLIANG@CS.UALBERTA.CA
*Department of Computing Science*
*University of Alberta*
*Edmonton, Alberta T6G 2E8, Canada*

**Xinhua Zhang**                                  XINHUA.ZHANG@NICTA.COM.AU
*Machine Learning Research Group*
*National ICT Australia (NICTA)*

*Research School of Computer Science*
*The Australian National University*
*Canberra, ACT 2601, Australia*

**Dale Schuurmans**                                   DALE@CS.UALBERTA.CA
*Department of Computing Science*
*University of Alberta*
*Edmonton, Alberta T6G 2E8, Canada*

**Editor:** Unknown

## Abstract

Structured sparsity is an important modeling tool that expands the applicability of convex formulations for data analysis, however it also creates significant challenges for efficient algorithm design. In this paper we investigate the generalized conditional gradient (GCG) algorithm for solving structured sparse optimization problems—demonstrating that, with some enhancements, it can provide a more efficient alternative to current state of the art approaches. After providing a comprehensive overview of the convergence properties of GCG, we develop efficient methods for evaluating polar operators, a subroutine that is required in each GCG iteration. In particular, we show how the polar operator can be efficiently evaluated in two important scenarios: dictionary learning and structured sparse estimation. A further improvement is achieved by interleaving GCG with fixed-rank local subspace optimization. A series of experiments on matrix completion, multi-class classification, multi-view dictionary learning and overlapping group lasso shows that the proposed method can significantly reduce the training cost of current alternatives.

## 1. Introduction

Sparsity is an important concept in high-dimensional statistics ([Bühlmann and van de Geer, 2011](#)) and signal processing ([Eldar and Kutyniok, 2012](#)), which has led to important application successes by reducing model complexity and improving interpretability of the results. Although it is common to promote sparsity by adding appropriate regularizers, such as the $l_1$ norm, sparsity can also present itself in more sophisticated forms, such as group sparsity ([Yuan and Lin, 2006](#); [Huang et al., 2011](#); [Zhao et al., 2009](#)), graph sparsity ([Kim](#) and Xing, 2009; Jacob et al., 2009), matrix rank sparsity (Candès and Recht, 2009), and other combinatorial sparsity patterns (Obozinski and Bach, 2012). These recent notions of *structured* sparsity (Bach et al., 2012b; Micchelli et al., 2013) have greatly enhanced the ability to model complex structural relationships in data; however, they also create significant computational challenges.

A currently popular optimization scheme for training structured sparse models is *accelerated proximal gradient* (APG) (Beck and Teboulle, 2009; Nesterov, 2013), which enjoys an optimal rate of convergence among black-box first-order procedures (Nesterov, 2013). An advantage of APG is that each iteration consists solely of computing a *proximal update* (PU), which for simple regularizers can have a very low complexity. For example, under $l_1$ norm regularization each iterate of APG reduces to a soft-shrinkage operator that can be computed in linear time, which partly explains its popularity for such problems. Unfortunately, for *structured* sparse regularizers the PU is highly nontrivial and can become a computational bottleneck. For example, the trace norm is often used to promote low rank solutions in matrix variable problems such as matrix completion (Candès and Recht, 2009); but here the associated PU requires a *full* singular value decomposition (SVD) of the gradient matrix on each iteration, which prevents APG from being applied to large problems. Not only does the PU require nontrivial computation for structured regularizers (Kim and Xing, 2010)—for example, requiring tailored network flow algorithms in some cases (Jenatton et al., 2011; Mairal et al., 2011; Mairal and Yu, 2013)—it yields *dense* intermediate iterates. Recently, Obozinski and Bach (2012) have demonstrated a class of regularizers where the corresponding PUs can be computed by a sequence of submodular function minimizations, but such an approach remains expensive.

To address the primary shortcomings with APG, we investigate the *generalized conditional gradient* (GCG) strategy for sparse optimization, which has been motivated by the promise of recent sparse approximation methods (Hazan, 2008; Clarkson, 2010; Zhang, 2003). A special case of GCG, known as the conditional gradient (CG), was originally proposed by Frank and Wolfe (1956) and has received significant renewed interest (Bach, 2013a; Clarkson, 2010; Freund and Grigas, 2013; Hazan, 2008; Jaggi and Sulovsky, 2010; Jaggi, 2013; Shalev-Shwartz et al., 2010; Tewari et al., 2011; Yuan and Yan, 2013). The key advantage of GCG for structured sparse estimation is that it need only compute the *polar* of the regularizer in each iteration, which sometimes admits a much more efficient update than the PU in APG. For example, under trace norm regularization, GCG only requires the *spectral norm* of the gradient to be computed in each iteration, which is an order of magnitude cheaper than evaluating the full SVD as required by APG. Furthermore, the greedy nature of GCG affords explicit control over the sparsity of intermediate iterates, which is not available in APG. Although existing work on GCG has generally been restricted to *constrained* optimization (i.e. enforcing an upper bound on the sparsity-inducing regularizer), in this paper we consider the more general *regularized* version of the problem. Despite their theoretical equivalence, the regularized form allows more efficient local optimization to be interleaved with the primary update, which provides a significant acceleration in practice.

This paper is divided into two major parts: first we present a general treatment of GCG and establish its convergence properties in Sections 2–3; then we apply the method to important case studies to demonstrate its effectiveness on large-scale, structured sparse estimation problems in Sections 4–6.

In particular, in the first part, after establishing notation and briefly discussing relevant optimization schemes (Section 2), we present our first main contribution (Section 3): a unified treatment of GCG and its convergence properties. More specifically, we begin with a general setup Section 3.1 that considers even nonconvex loss functions, then progressively establish tighter convergence analyses by strengthening the assumptions. The treatment concludes in Section 3.4 with the case of gauge function regularizers (e.g. norms), where we prove that, under standard assumptions, GCG converges to the global optimum at the rate of $O(1/t)$—which is a significant improvement over previous bounds (Dudik et al., 2012). Although the convergence rate for GCG remains inferior to that of APG, these two algorithms demonstrate alternative trade-offs between the number of iterations required and the per-step complexity. Indeed, as shown in our experiments, GCG can be significantly faster than APG in terms of *overall* computation time for large problems. While some of these results follow from existing work, many of the results we establish for GCG (such as Theorem 7, Theorem 8 and Theorem 12) are new.

Equipped with these technical results, the second part of the paper then applies the GCG algorithm to two important application areas: low rank learning (Section 4) and structured sparse estimation (Section 5).

For low rank learning (Section 4), since imposing a hard bound on matrix rank generally leads to an intractable problem, in Section 4.2 we first present a generic relaxation based on gauge functions (Chandrasekaran et al., 2012; Tewari et al., 2011) that yields a convex program (Bach et al., 2008; Bradley and Bagnell, 2009; Lee et al., 2009; Zhang et al., 2012). Conveniently, the resulting problem can be easily optimized using the GCG algorithm from Section 3.4. To further reduce computation time, we introduce in Section 4.3 an auxiliary fixed rank subspace optimization within each iteration of GCG. Although similar hybrid approaches have been previously suggested (Burer and Monteiro, 2005; Journee et al., 2010; Laue, 2012; Mishra et al., 2013), we propose an efficient new alternative that, instead of locally optimizing a *constrained* fixed rank problem, optimizes an *unconstrained* surrogate objective based on a variational representation of matrix norms. This alternative strategy allows a far more efficient local optimization without compromising GCG's convergence properties. In Section 4.4 we show that the approach can be applied to dictionary learning and matrix completion problems, as well as extended to a non-trivial multi-view form of dictionary learning (White et al., 2012; Lee et al., 2009).

Next, for structured sparse estimation (Section 5), we consider a family of structured regularizers that are induced by cost functions defined on variable subsets (Obozinski and Bach, 2012). Here our main contribution is to characterize a rich class of such regularizers that allow efficient computation of their polar operators, which enables the GCG algorithm from Section 3 to be efficiently applied. In Section 5.2 we introduce a "lifting" construction that expresses these structured sparse regularizers as "marginalized" linear functions, whose polar, after some reformulation, can be efficiently evaluated by a simple linear program (LP). An important example is the overlapping group lasso (Jenatton et al., 2011), treated in detail in Section 5.3. By exploiting problem structure, we are able to further reduce the LP to piecewise linear minimization over a weighted simplex, enabling further speed-up via smoothing (Nesterov, 2005).

Finally, in Section 6 we provide an extensive experimental evaluation that compares the performance of GCG (augmented with interleaved local search) to state-of-the-art op-timization strategies, across a range of problems, including matrix completion, multi-class classification, multi-view dictionary learning and overlapping group lasso.

## 1.1 Extensions over Previously Published Work

Some of the results in this article have appeared in a preliminary form in two conference papers (Zhang et al., 2013, 2012); however, several specific extensions have been added in this paper, including:

- A more extensive and unified treatment of GCG in general Banach spaces.

- A more refined analysis of the convergence properties of GCG.

- A new application of GCG to structured sparse training for multi-view dictionary learning.

- New experimental evaluations on larger data sets, including new results for matrix completion and image multi-class classification.

We have also released open source Matlab code, which is available for download at http://users.cecs.anu.edu.au/~xzhang/GCG.

## 2. Preliminaries

We first establish some of the definitions and notation that will be used throughout the paper, before providing a brief background on the relevant optimization methods that establish a context for discussing GCG.

## 2.1 Notation and Definitions

Throughout this paper we will assume that the underlying space, unless otherwise stated, is a real Banach space $\mathcal{B}$ with norm $\|\cdot\|$; for example, $\mathcal{B}$ could be the Euclidean space $\mathbb{R}^d$ with any norm. The dual space (the set of all real-valued continuous linear functions on $\mathcal{B}$) is denoted $\mathcal{B}^*$ and is equipped with the dual norm $\|\mathbf{g}\|^\circ = \sup\{\mathbf{g}(\mathbf{w}) : \|\mathbf{w}\| \leq 1\}$. The generality we enjoy here does not complicate any of our later analyses hence is preferred.

We use bold lowercase letters to denote vectors, where the $i$-th component of a vector $\mathbf{w}$ is denoted $w_i$, while $\mathbf{w}_t$ denotes some other vector. We use the shorthand $\langle \mathbf{w}, \mathbf{g} \rangle := \mathbf{g}(\mathbf{w})$ for any $\mathbf{g} \in \mathcal{B}^*$ and $\mathbf{w} \in \mathcal{B}$. The $\mathsf{l_p}$ norm (with $\mathsf{p} \geq 1$) on $\mathbb{R}^d$ is defined as $\|\mathbf{w}\|_\mathsf{p} = (\sum_{i=1}^d |w_i|^\mathsf{p})^{1/\mathsf{p}}$, with its dual norm given by $\mathsf{l_q}$ with $1/\mathsf{p} + 1/\mathsf{q} = 1$. The vectors of all 1s and all 0s are denoted $\mathbf{1}$ and $\mathbf{0}$ respectively, with size clear from context. We also let $(t)_+ := \max\{t, 0\}$, and adopt the abbreviation $\overline{\mathbb{R}} := \mathbb{R} \cup \{\infty\}$.

The function $f : \mathcal{B} \to \overline{\mathbb{R}}$ is called $\sigma$-strongly convex (w.r.t. the norm $\|\cdot\|$) if there exists some $\sigma \geq 0$ such that for all $0 \leq \lambda \leq 1$ and $\mathbf{w}, \mathbf{z} \in \mathcal{B}$,

$$f(\lambda\mathbf{w} + (1-\lambda)\mathbf{z}) + \frac{\sigma\lambda(1-\lambda)}{2}\|\mathbf{w} - \mathbf{z}\|^2 \leq \lambda f(\mathbf{w}) + (1-\lambda)f(\mathbf{z}).$$

In the case where $\sigma = 0$ we simply say $f$ is convex. For a convex function $f : \mathcal{B} \to \overline{\mathbb{R}}$, we denote its subdifferential at $\mathbf{w}$ by $\partial f(\mathbf{w}) = \{\mathbf{g} \in \mathcal{B}^* : f(\mathbf{z}) \geq f(\mathbf{w}) + \langle \mathbf{z} - \mathbf{w}, \mathbf{g} \rangle, \forall \mathbf{z} \in \mathcal{B}\}$.

For a differentiable function $\ell : \mathcal{B} \to \overline{\mathbb{R}}$, we use $\nabla \ell(\mathbf{w})$ to denote its (Fréchet) derivative at $\mathbf{w}$. The Fenchel conjugate of the function $f : \mathcal{B} \to \overline{\mathbb{R}}$, denoted $f^* : \mathcal{B}^* \to \overline{\mathbb{R}}$, is defined by $f^*(\mathbf{g}) = \sup_{\mathbf{w}} \langle \mathbf{w}, \mathbf{g} \rangle - f(\mathbf{w})$ for $\mathbf{g} \in \mathcal{B}^*$. Note that $f^*$ is always convex. For any set $C \subseteq \mathcal{B}$, we define the indicator function

$$\iota_C(\mathbf{w}) = \begin{cases} 0, & \text{if } \mathbf{w} \in C \\ \infty, & \text{otherwise} \end{cases}.$$

Conveniently, $C$ is a convex set iff $\iota_C$ is a convex function. We use $\text{int}(C)$ to denote the interior of $C$, i.e. the largest open subset of $C$. A convex function $f : \mathcal{B} \to \overline{\mathbb{R}}$ is *closed* iff its sublevel sets $\{ \mathbf{w} \in \mathcal{B} : f(\mathbf{w}) \le \alpha \}$ are closed sets for all $\alpha \in \mathbb{R}$, and *proper* iff its domain $\text{dom} \, f := \{ \mathbf{w} \in \mathcal{B} : f(\mathbf{w}) < \infty \}$ is nonempty. Collectively we denote $\Gamma_0$ as the set of all closed proper and convex functions on $\mathcal{B}$. A map from a Banach space $(\mathcal{B}, \|\cdot\|)$ to another Banach space $(\mathcal{B}', |\!|\!|\cdot|\!|\!|)$ is said to be $L$-Lipschitz continuous (with Lipschitz constant $L < \infty$) if for all $\mathbf{w}, \mathbf{z} \in \mathcal{B}$, $|\!|\!| f(\mathbf{w}) - f(\mathbf{z}) |\!|\!| \le L \|\mathbf{w} - \mathbf{z}\|$. The identity map $\mathsf{Id} : \mathcal{B} \to \mathcal{B}$ maps any $\mathbf{w} \in \mathcal{B}$ to itself.

Finally, we will use uppercase letters to denote matrices. For a matrix $W \in \mathbb{R}^{m \times n}$, $\|W\|_{\text{tr}}, \|W\|_{\text{sp}}$ and $\|W\|_{\text{F}}$ denote its trace norm (sum of singular values), spectral norm (largest singular value) and Frobenius norm (the $l_2$ norm of the singular values), respectively. The trace norm is dual to the spectral norm, while the Frobenius norm is dual to itself. We use $W_{i:}$ and $W_{:j}$ to denote the $i$-th row and $j$-th column of $W$, respectively. The transpose of $W$ is denoted as $W^\top$. For (real) symmetric matrices $W$ and $Z$, the notation $W \succeq Z$ asserts that $W - Z$ is positive semidefinite. We use $I$ to denote the identity matrix (with size clear from the context). Finally, for a set $A$, we let $|A|$ denote its cardinality (i.e., the number of elements in $A$).

Below we will also need to make use of the following results.

**Proposition 1.** *For a function $\ell : \mathcal{B} \to \mathbb{R}$, if the gradient $\nabla \ell : \mathcal{B} \to \mathcal{B}^*$ is $L$-Lipschitz continuous, then for all $\mathbf{w}, \mathbf{z} \in \mathcal{B}$, we have $\ell(\mathbf{w}) \le \tilde{\ell}_L(\mathbf{w}; \mathbf{z})$, where*

$$\tilde{\ell}_L(\mathbf{w}; \mathbf{z}) := \ell(\mathbf{z}) + \langle \mathbf{w} - \mathbf{z}, \nabla \ell(\mathbf{z}) \rangle + \frac{L}{2} \|\mathbf{w} - \mathbf{z}\|^2. \tag{1}$$

That is, $\tilde{\ell}_L(\,\cdot\,; \mathbf{z})$ is a quadratic upper approximation of $\ell(\cdot)$, centered at $\mathbf{z}$. This well-known result can be proved by a simple application of the mean value theorem. We will make repeated use of this upper approximation in our analysis below. When $\ell$ is convex, the inequality in Proposition 1 furthermore implies that $\nabla \ell$ is $L$-Lipschitz continuous. A further result we will make use of is the following.

**Proposition 2.** *For a convex function $\ell : \mathcal{B} \to \mathbb{R}$, the gradient $\nabla \ell : \mathcal{B} \to \mathcal{B}^*$ is $L$-Lipschitz continuous if and only if its Fenchel conjugate $\ell^* : \mathcal{B}^* \to \overline{\mathbb{R}}$ is $1/L$-strongly convex.*

A proof of Proposition 2 can be found in (Zălinescu, 2002, Corollary 3.5.11).

## 2.2 Composite Minimization Problem

Many machine learning problems, particularly those formulated in terms of regularized risk minimization, can be formulated as:

$$\inf_{\mathbf{w} \in \mathcal{B}} F(\mathbf{w}), \text{ such that } F(\mathbf{w}) := \ell(\mathbf{w}) + f(\mathbf{w}), \tag{2}$$

where $f$ is closed, proper and convex, and $\ell$ is continuously differentiable. Typically, $\ell$ is a loss function and $f$ is a regularizer although their roles can be reversed in some common examples, such as support vector machines (Cortes and Vapnik, 1995) and least absolute shrinkage (Tibshirani, 1996).

Due to the importance of this problem in machine learning, significant effort has been devoted to designing efficient algorithms for solving (2). Standard methods are generally based on computing an update sequence $\{\mathbf{w}_t\}$ that converges to a global minimizer of $F$ when $F$ is convex, or converges to a stationary point otherwise. For example, a popular example is the mirror descent (MD) algorithm (Beck and Teboulle, 2003), which takes an arbitrary subgradient $\mathbf{g}_t \in \partial F(\mathbf{w}_t)$ (provided $F$ is convex) and successively performs the update

$$\mathbf{w}_{t+1} = \arg\min_{\mathbf{w}} \left\{ \langle \mathbf{w}, \mathbf{g}_t \rangle + \frac{1}{2\eta_t} \mathsf{D}(\mathbf{w}, \mathbf{w}_t) \right\} \tag{3}$$

with a suitable step size $\eta_t \geq 0$; here $\mathsf{D}(\mathbf{w}, \mathbf{z}) := \mathsf{d}(\mathbf{w}) - \mathsf{d}(\mathbf{z}) - \langle \mathbf{w} - \mathbf{z}, \nabla \mathsf{d}(\mathbf{z}) \rangle$ denotes a Bregman divergence induced by some differentiable 1-strongly convex function $\mathsf{d} : \mathcal{B} \to \overline{\mathbb{R}}$. Under weak assumptions on the subdifferential, it can be shown that MD converges at a rate of $O(1/\sqrt{t})$ (Beck and Teboulle, 2003). However, this algorithm ignores the composite structure in (2) and can be impractically slow for many applications.

Another widely used algorithm is the proximal gradient (PG) of Fukushima and Mine (1981), also known as the forward-backward splitting procedure (Bauschke and Combettes, 2011), which iteratively performs the proximal update (PU)

$$\begin{aligned}
\mathbf{w}_{t+1} &= \arg\min_{\mathbf{w}} \left\{ \tilde{\ell}_0(\mathbf{w}; \mathbf{w}_t) + f(\mathbf{w}) + \frac{1}{2\eta_t} \mathsf{D}(\mathbf{w}, \mathbf{w}_t) \right\} \\
&= \arg\min_{\mathbf{w}} \left\{ \langle \mathbf{w}, \nabla \ell(\mathbf{w}_t) \rangle + f(\mathbf{w}) + \frac{1}{2\eta_t} \mathsf{D}(\mathbf{w}, \mathbf{w}_t) \right\}.
\end{aligned} \tag{4}$$

The difference between the MD and PU updates is that the latter does not linearize the regularizer $f$. The downside is that (4) is usually harder to solve than (3) for each iterate, with an upside that a faster $O(1/t)$ rate of convergence is achieved (Tseng, 2010).[1] In particular, when $\mathcal{B} = \mathbb{R}^d$, $\mathsf{D}(\mathbf{w}, \mathbf{z}) = \frac{1}{2} \|\mathbf{w} - \mathbf{z}\|_2^2$ and $f(\mathbf{w}) = \|\mathbf{w}\|_1$, the soft-shrinkage operator is recovered, which plays an important role in sparse estimation:

$$\mathbf{z}_t = \mathbf{w}_t - \eta_t \nabla \ell(\mathbf{w}_t) \tag{5}$$

$$\mathbf{w}_{t+1} = (1 - \eta_t / |\mathbf{z}_t|)_+ \, \mathbf{z}_t, \tag{6}$$

where the algebraic operations in (6) are performed componentwise. Similarly, when $\mathcal{B} = \mathbb{R}^{m \times n}$, $\mathsf{D}(W, Z) = \frac{1}{2} \|W - Z\|_{\mathrm{F}}^2$ and $f(W) = \|W\|_{\mathrm{tr}}$, one recovers the matrix analogue of soft-shrinkage used for matrix completion (Candès and Recht, 2009): here the same gradient update (5) is applied, but the shrinkage operator (6) is applied to the singular values of $Z_t$ while keeping the singular vectors intact. Unfortunately, such a PU can become prohibitively expensive for large matrices, since it requires a *full* singular value decomposition (SVD) at each iteration, at a cubic time cost. We will see in the next section that an alternative algorithm only requires computing the dual of the trace norm (namely the spectral norm), which requires only quadratic time per iteration.

---

1. Moreover, an additional extrapolation variant, accelerated proximal gradient (APG), can further improve the convergence rate to $O(1/t^2)$ (Beck and Teboulle, 2009; Nesterov, 2013).

## 3. Generalized Conditional Gradient

The first main contribution of this paper is to expand the understanding of the convergence properties of the generalized conditional gradient (GCG) algorithm, which has recently been receiving renewed attention (Bach, 2013a; Bredies et al., 2009; Clarkson, 2010; Freund and Grigas, 2013; Hazan, 2008; Jaggi and Sulovsky, 2010; Jaggi, 2013; Shalev-Shwartz et al., 2010; Tewari et al., 2011; Yuan and Yan, 2013; Zhang et al., 2012). After some enhancements, the algorithm and convergence guarantees we establish in this section will then be applied to dictionary learning and sparse estimation in Sections 4 and 5 respectively.

### 3.1 General Setting

Our focus in this paper is on developing efficient solution methods for (2). We begin by considering a general setting where $\ell$ is not assumed to be convex; in particular, we initially assume only that $\ell$ is continuously differentiable and $f$ is closed, proper and convex.

First observe that if $\ell$ is continuously differentiable and $f$ is convex, then $F$ is locally Lipschitz continuous, hence any local minimum $\mathbf{w}$ in (2) must satisfy the necessary condition

$$0 \in \partial F(\mathbf{w}) = \nabla \ell(\mathbf{w}) + \partial f(\mathbf{w}), \tag{7}$$

where $\partial F$ denotes Clarke's generalized gradient (Clarke, 1990, Proposition 2.3.2).[2] Recalling that $f^*$ denotes the Fenchel conjugate of $f$, and using the fact that $\mathbf{g} \in \partial f(\mathbf{w})$ iff $\mathbf{w} \in \partial f^*(\mathbf{g})$ when $f$ is closed, proper and convex (Zălinescu, 2002, Theorem 2.4.2 (iii)), one can observe that the necessary condition for a minimum can also be expressed as

$$\begin{aligned}
(7) &\iff -\nabla \ell(\mathbf{w}) \in \partial f(\mathbf{w}) \\
&\iff \mathbf{w} \in \partial f^*(-\nabla \ell(\mathbf{w})) \\
&\iff \mathbf{w} - (1-\eta)\mathbf{w} \in \eta \partial f^*(-\nabla \ell(\mathbf{w})), \tag{8} \\
&\iff \mathbf{w} \in (1-\eta)\mathbf{w} + \eta \partial f^*(-\nabla \ell(\mathbf{w})), \tag{9}
\end{aligned}$$

when $\eta \in (0, 1)$. Then, by the definition of $f^*$, we have

$$\begin{aligned}
\mathbf{d} \in \partial f^*(-\nabla \ell(\mathbf{w})) &\iff \mathbf{d} \in \arg\max_{\mathbf{d}} \left\{ \langle \mathbf{d}, -\nabla \ell(\mathbf{w}) \rangle - f(\mathbf{d}) \right\} \\
&\iff \mathbf{d} \in \arg\min_{\mathbf{d}} \left\{ \langle \mathbf{d}, \nabla \ell(\mathbf{w}) \rangle + f(\mathbf{d}) \right\}, \tag{10}
\end{aligned}$$

thus the necessary condition for a minimum can be characterized by the fixed-point inclusion

$$\mathbf{w} \in (1-\eta)\mathbf{w} + \eta \mathbf{d} \ \text{ for some } \ \mathbf{d} \in \arg\min_{\mathbf{d}} \left\{ \langle \mathbf{d}, \nabla \ell(\mathbf{w}) \rangle + f(\mathbf{d}) \right\}. \tag{11}$$

---

2. The notation $\partial$ can be sensibly used to denote both the subdifferential and generalized gradient in (7), since the correspondence will be maintained when convexity is present.

---

**Algorithm 1** Generalized Conditional Gradient (GCG).

---
1: Initialize: $\mathbf{w}_0 \in \operatorname{dom} f$.
2: **for** $t = 0, 1, \ldots$ **do**
3:     $\mathbf{g}_t = \nabla \ell(\mathbf{w}_t)$
4:     $\mathbf{d}_t \in \arg\min_{\mathbf{d}} \{\langle \mathbf{d}, \mathbf{g}_t \rangle + f(\mathbf{d})\}$
5:     Choose step size $\eta_t \in [0, 1]$
6:     $\tilde{\mathbf{w}}_{t+1} = (1 - \eta_t)\mathbf{w}_t + \eta_t \mathbf{d}_t$
7:     $\mathbf{w}_{t+1} = \texttt{Improve}(\tilde{\mathbf{w}}_{t+1}, \ell, f)$         ▷ Subroutine, see Definition 1
8: **end for**

---

This particular fixed-point condition immediately suggests an update that provides the foundation for the generalized conditional gradient algorithm (GCG) outlined in Algorithm 1. Each iteration of this procedure involves linearizing the smooth loss $\ell$, solving the subproblem (10), selecting a step size, taking a convex combination, then conducting some form of local improvement.[3]

Note that the subproblem (10) shares some similarity with the proximal update (4): both choose to leave the potentially nonsmooth function $f$ intact, but here the smooth loss $\ell$ is replaced by its plain linearization rather than its linearization plus a (strictly convex) proximal regularizer (e.g. the quadratic upper bound $\tilde{\ell}_L(\cdot)$ in (1)). Consequently, the subproblem (10) might have multiple solutions, in which case we simply adopt any one, or no solution (e.g. divergence), in which case we will impose extra assumptions to ensure boundedness (details below).

An important component of Algorithm 1 is the final step, $\texttt{Improve}$, in Line 7, which is particularly important in practice: it allows for heuristic local improvement to be conducted on the iterates without sacrificing any theoretical convergence properties of the overall procedure. Since $\texttt{Improve}$ has access both $\ell$ and $f$ it can be very powerful in principle—we will consider the following variants that have respective consequences for practice and subsequent analysis.

**Definition 1.** *The subroutine* $\texttt{Improve}$ *is called* $\texttt{Null}$ *if for all $t$, $\mathbf{w}_{t+1} = \tilde{\mathbf{w}}_{t+1}$; $\texttt{Descent}$ if $F(\mathbf{w}_{t+1}) \leq F(\tilde{\mathbf{w}}_{t+1})$; and* $\texttt{Monotone}$ *if $F(\mathbf{w}_{t+1}) \leq F(\mathbf{w}_t)$. Later (when $\nabla \ell$ will be assumed to be $L$-Lipschitz continuous) the subroutine* $\texttt{Improve}$ *will also be called* $\texttt{Relaxed}$ *if*

$$F(\mathbf{w}_{t+1}) \leq \tilde{\ell}_L(\tilde{\mathbf{w}}_{t+1}; \mathbf{w}_t) + (1 - \eta_t)f(\mathbf{w}_t) + \eta_t f(\mathbf{d}_t).[4] \tag{12}$$

---

Our main goal remains to establish that Algorithm 1 indeed converges to a stationary point under general conditions. To do so, we will make the following three assumptions.

**Assumption 1.** *$\ell$ is continuously differentiable in an open set that contains* dom $f$; *$f$ is closed, proper and convex; and* $-\nabla\ell(\text{dom } f) \subseteq \text{Range}(\partial f)$. *Note that we are using the notation* $-\nabla\ell(\text{dom } f) := \bigcup_{\mathbf{w}\in\text{dom } f}\{-\nabla\ell(\mathbf{w})\}$ *and* $\text{Range}(\partial f) := \bigcup_{\mathbf{w}\in\text{dom } f}\partial f(\mathbf{w})$ *in this definition.*

The last condition simply ensures that the subproblem (10) has at least one solution.

**Assumption 2.** *$\nabla\ell$ is uniformly continuous on bounded sets.*

Recall that a set $C \subseteq \mathcal{B}$ is bounded if $\sup_{\mathbf{w},\mathbf{z}\in C}\|\mathbf{w}-\mathbf{z}\| < \infty$. Also recall that a function $f : \mathcal{B} \to \overline{\mathbb{R}}$ is uniformly continuous if for all $\epsilon > 0$ there exists $\delta > 0$ such that for all $\mathbf{w}, \mathbf{z} \in \text{dom } f$, $\|\mathbf{w}-\mathbf{z}\| < \delta \implies |f(\mathbf{w}) - f(\mathbf{z})| < \epsilon$. Uniform continuity is stronger than continuity since here $\delta$ is chosen independently of $\mathbf{w}$ and $\mathbf{z}$ (hence the name "uniform"), but it is weaker than Lipschitz continuity. When $\mathcal{B}$ is of finite dimension, Assumption 2 is automatically satisfied whenever $\nabla\ell$ is continuous (since the closure of a bounded set in a finite dimensional space is compact and continuous functions on compact sets are uniformly continuous). The uniform continuity of $\nabla\ell$ also follows if $\nabla\ell$ is Lipschitz continuous on bounded sets.

**Assumption 3.** *The sequences $\{\mathbf{d}_t\}$ and $\{\mathbf{w}_t\}$ generated by Algorithm 1 are bounded for any choice of $\eta_t \in [0,1]$.*

Although this assumption is more stringent, it can fortunately be achieved under various conditions, which we summarize as follows.

**Proposition 3.** *Let Assumption 1 and Assumption 2 hold, then Assumption 3 is satisfied if any of the following hold.*

(a) *The subroutine* `Improve` *is* `Monotone`; *the sublevel set $\{\mathbf{w} \in \text{dom } f : F(\mathbf{w}) \leq F(\mathbf{w}_0)\}$ is compact; and $-\nabla\ell(\text{dom } f) \subseteq \text{int}(\text{dom } f^*)$. The latter condition holds, in particular, when $f$ is cofinite, i.e., $f^*$ has full domain.*

(b) *The subroutine* `Improve` *is* `Monotone`; *the sublevel set $\{\mathbf{w} \in \text{dom } f : F(\mathbf{w}) \leq F(\mathbf{w}_0)\}$ is bounded; and $f$ is super-coercive, i.e., $\lim_{\|\mathbf{w}\|\to\infty} f(\mathbf{w})/\|\mathbf{w}\| \to \infty$.*

(c) *dom $f$ is bounded.*

The proof is given in Appendix A.1.[5]

---

[5]. Slightly more restricted forms of these three conditions have been previously considered by (Mine and Fukushima, 1981; Bredies et al., 2009), and (Frank and Wolfe, 1956; Dem'yanov and Rubinov, 1967; Levitin and Polyak, 1966), respectively. Note that condition (a) implies $F$ is bounded from below, while condition (b) implies the same conclusion if $F$ is convex and $\mathcal{B}$ is a Hilbert space (or more generally a reflexive space, i.e. such that $(\mathcal{B}^*)^* = \mathcal{B}$). It is interesting to compare condition (a) and (b): there appears to be a trade-off between restricting the sublevel sets of $F$ versus restricting the behavior of $f^*$. In particular, $f$ being super-coercive implies it is cofinite, but the converse is only true in finite dimensional spaces (Borwein and Vanderwerff, 2010).

Finally, we define a quantity, referred to as the *duality gap*, that will be useful in understanding the convergence properties of GCG:

$$\mathsf{G}(\mathbf{w}) := F(\mathbf{w}) - \inf_{\mathbf{d}} \{\ell(\mathbf{w}) + \langle \mathbf{d} - \mathbf{w}, \nabla\ell(\mathbf{w})\rangle + f(\mathbf{d})\}$$
$$= \langle \mathbf{w}, \nabla\ell(\mathbf{w})\rangle + f(\mathbf{w}) - \inf_{\mathbf{d}} \{\langle \mathbf{d}, \nabla\ell(\mathbf{w})\rangle + f(\mathbf{d})\}$$
$$= \langle \mathbf{w}, \nabla\ell(\mathbf{w})\rangle + f(\mathbf{w}) + f^*(-\nabla\ell(\mathbf{w})). \tag{13}$$

As long as $f$ is convex (which we are assuming throughout this paper) the duality gap provides an upper bound on the suboptimality of any search point, as established in the following proposition.

**Proposition 4.** *For $F = \ell + f$, if $\ell$ is differentiable and $f$ is convex, then $\mathsf{G}(\mathbf{w}) \geq 0$ with equality iff $\mathbf{w}$ satisfies the necessary condition (7). If $\ell$ is also convex, then $\mathsf{G}(\mathbf{w}) \geq F(\mathbf{w}) - \inf_{\mathbf{d}} F(\mathbf{d}) \geq 0$ for all $\mathbf{w} \in \mathcal{B}$.*

A proof is given in Appendix A.2 and an illustration is given in Figure 1. Since $\mathsf{G}(\mathbf{w}_t)$ provides an upper bound on the suboptimality of $\mathbf{w}_t$, the duality gap gives a natural stopping criterion for (2) that can be easily computed in each iteration of Algorithm 1 (provided $f^*$ can be conveniently evaluated).

We are now ready for the first convergence result regarding Algorithm 1.

**Theorem 5.** *Let Assumption 1 hold, equip Algorithm 1 with a `Descent` subroutine for the local improvement phase, and choose step sizes $\eta_t$ such that*

$$F(\bar{\mathbf{w}}_{t+1}) \leq \min_{0 \leq \eta \leq 1} \{\ell((1-\eta)\mathbf{w}_t + \eta\mathbf{d}_t) + (1-\eta)f(\mathbf{w}_t) + \eta f(\mathbf{d}_t)\}. \tag{14}$$

*Then for all $t$, either $F(\mathbf{w}_{t+1}) < F(\mathbf{w}_t)$ or $\mathsf{G}(\mathbf{w}_t) = 0$. Additionally, if Assumption 2 and Assumption 3 also hold, then either $F(\mathbf{w}_t) \downarrow -\infty$, or $\mathsf{G}(\mathbf{w}_t) = 0$ indefinitely, or $\mathsf{G}(\mathbf{w}_t) \to 0$.*

The proof, which can be found in Appendix A.3, slightly generalizes the argument of (Bredies et al., 2009). Recall that the $\ell$ in Theorem 5 need not be convex, which partly explains why one can only establish convergence to a stationary point in this case.

Some remarks about Theorem 5 are in order.

**Remark 1.** *Due to the possible non-uniqueness of the solution to (10), $\mathsf{G}(\mathbf{w}_t) = 0$ does not imply $\mathsf{G}(\mathbf{w}_s) = 0$ for all $s > t$. However, when $F$ is convex, $\mathsf{G}(\mathbf{w}_t) = 0$ implies $\mathbf{w}_t$ is globally optimal (cf. Proposition 4), in which case it follows that $\mathsf{G}(\mathbf{w}_s) = 0$ for all $s > t$, since the monotonicity $F(\mathbf{w}_s) \leq F(\mathbf{w}_t)$ implies the global optimality of $\mathbf{w}_s$. On the other hand, if $f$ is strictly convex, then achieving $\mathsf{G}(\mathbf{w}_t) = 0$ for some $t$ implies $\mathbf{d}_t = \mathbf{w}_t$, hence $\mathbf{w}_s = \mathbf{w}_t$ for all $s \geq t$, provided the `Null` subroutine is employed in Algorithm 1.*

**Remark 2.** *It is easily seen that the duality gap $\mathsf{G}$ in (13) is lower semicontinuous, hence if $\mathsf{G}(\mathbf{w}_t) \to 0$ and any subsequence of $\mathbf{w}_t$ converges to some point, say $\mathbf{w}$, it must follow that $0 \leq \mathsf{G}(\mathbf{w}) \leq \liminf \mathsf{G}(\mathbf{w}_{t_k}) = 0$; that is, the limit point $\mathbf{w}$ is a stationary point. Note that if $F$ has compact level sets, $\mathbf{w}_t$ is guaranteed to have a convergent subsequence.*

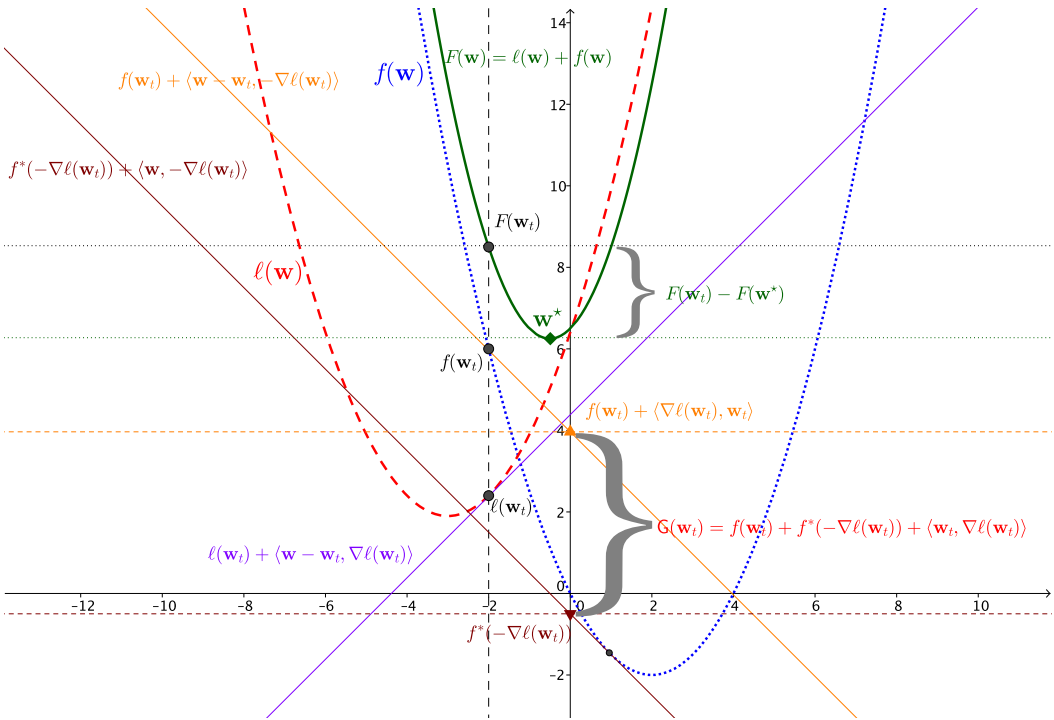

Figure 1: Here both $\ell$ (the red dashed parabolic) and $f$ (the blue dotted parabolic) are convex (quadratic) functions, and $\mathbf{w}^\star$ denotes the minimizer of $F = \ell + f$ (the green solid parabolic). The dashed vertical line represents our current iterate $\mathbf{w}_t$. As predicted, the duality gap $\mathsf{G}(\mathbf{w}_t) \geq F(\mathbf{w}_t) - F(\mathbf{w}^\star)$.

**Remark 3.** *One of the small improvements contributed by Theorem 5 over previous work (Mine and Fukushima, 1981; Bredies et al., 2009) is the step size rule (14). Earlier work insisted on choosing*

$$\eta_t \in \arg\min_{0 \leq \eta \leq 1} \left\{ \ell\big((1-\eta)\mathbf{w}_t + \eta\mathbf{d}_t\big) + f\big((1-\eta)\mathbf{w}_t + \eta\mathbf{d}_t\big) \right\}, \tag{15}$$

*which is a special case of (14), owing to the convexity of $f$. However, the sufficiency of (14) has a practical consequence, since minimizing the right-hand side of (14) can be significantly easier than dealing with (15). Indeed, in an application to a linear inverse problem, Bredies et al. (2009) needed to develop specialized subroutines (under additional assumptions on a parameter $p$) to solve (15), while the right-hand side of (14) would be trivial to apply in the same setting (without any assumption on the parameter $p$). Nevertheless, from the perspective of greedily decreasing the objective, (15) considers a tighter approximation.*

**Remark 4.** *Perhaps surprisingly, the proof of Theorem 5 does not use the convexity of $f$ explicitly. Convexity is implicitly needed only in two places: the tractability of the subproblem (10) and the satisfiability of the step size rule (14).*

### 3.2 Refined Analysis when $\nabla\ell$ is Lipschitz Continuous

The analysis of Theorem 5 can be tightened if we make a slightly more restrictive assumption.

**Assumption 4.** *There exists some positive constant $L < \infty$ such that for the sequences $\{\mathbf{w}_t\}$ and $\{\mathbf{d}_t\}$ generated by Algorithm 1, and for all $\eta \in [0, 1]$, we have*

$$\ell(\mathbf{w}_t + \eta(\mathbf{d}_t - \mathbf{w}_t)) \leq \ell(\mathbf{w}_t) + \eta\langle\mathbf{d}_t - \mathbf{w}_t, \nabla\ell(\mathbf{w}_t)\rangle + \frac{L\eta^2}{2}\|\mathbf{d}_t - \mathbf{w}_t\|^2. \quad (16)$$

The inequality (16) essentially asserts that the gradient $\nabla\ell : (\mathcal{B}, \|\cdot\|) \to (\mathcal{B}^*, \|\cdot\|^\circ)$ is $L$-Lipschitz continuous, imposing a bound on $\ell$'s curvature. Note that under Assumption 1 and Assumption 4, a `Descent` subroutine is automatically `Relaxed`. This additional assumption allows a slightly improved convergence result over Theorem 5.

**Theorem 6.** *Let Assumptions 1–4 hold, equip Algorithm 1 with the subroutine `Null`, and choose step sizes $\eta_t$ such that*

$$\eta_t = \min\left\{\frac{\mathsf{G}(\mathbf{w}_t)}{L\|\mathbf{w}_t - \mathbf{d}_t\|^2}, 1\right\}. \quad (18)$$

*Then as $t \to \infty$, either $F(\mathbf{w}_t) \downarrow -\infty$ or $G(\mathbf{w}_t) \to 0$.*

The adaptive step size (18) is in some sense "optimal" in that it minimizes the quadratic upper bound (16) from Assumption 4. Theorem 6 extends (Levitin and Polyak, 1966, Theorem 6.1 (1)) to accommodate a general (closed, proper and convex) $f$, instead of restricting $f$ to an indicator; see Appendix A.4 for a proof. The main advantage of Theorem 6 over Theorem 5 is that the step size can now be determined by the simple formula (18) rather than requiring a nonconvex optimization to be approximately solved as in (14). In particular, when $\ell$ is not convex, achieving (14) can be highly non-trivial.

Although the results, so far, do not offer any *rates*, a non-adaptive step size can be shown to yield a convergence guarantee along with a simple bound on the asymptotic rate.

**Theorem 7.** *Let Assumptions 1–4 hold, equip Algorithm 1 with the subroutine `Relaxed`, and assume that at time step $t$ the subproblem (10) is solved up to some additive error $\varepsilon_t$. Let $H_t = \sum_{s=0}^{t} \eta_s$ be the partial sum. Then for all $t \geq 0$*

$$\sum_{s=0}^{t}\left[\eta_s\mathsf{G}(\mathbf{w}_s) - \frac{L\eta_s^2}{2}\|\mathbf{w}_s - \mathbf{d}_s\|^2 - \eta_s\varepsilon_s\right] \leq F(\mathbf{w}_0) - F(\mathbf{w}_{t+1}). \quad (19)$$

*Moreover, if $F$ is bounded from below, $\sum_t \eta_t = \infty$, $\sum_t \eta_t^2 < \infty$, and $\varepsilon_t = O(1/H_t^{1+\delta})$ for some $\delta > 0$, then $\liminf_{t\to\infty} \mathsf{G}(\mathbf{w}_t)H_t = 0$.*

---

6. In fact, a weaker inequality is sufficient (which does not even require a topology on the space $\mathcal{B}$):

$$\ell(\mathbf{w}_t + \eta(\mathbf{d}_t - \mathbf{w}_t)) \leq \ell(\mathbf{w}_t) + \eta\langle\mathbf{d}_t - \mathbf{w}_t, \nabla\ell(\mathbf{w}_t)\rangle + \frac{L_F\eta^2}{2} \text{ for some constant } L_F < \infty. \quad (17)$$

Indeed, let $\rho$ be the smallest number so that the ball with radius $\rho$ contains the sequences $\{\mathbf{w}_t\}$ and $\{\mathbf{d}_t\}$. Then $L_F \leq L\rho^2$ provided that $\ell$ satisfies Assumption 4 with constant $L$. According to Proposition 1, (16) holds as long as the gradient $\nabla\ell : (\mathcal{B}, \|\cdot\|) \to (\mathcal{B}^*, \|\cdot\|^\circ)$ is $L$-Lipschitz continuous, in which case Assumption 2 also becomes trivial. The reason to use the stronger condition (16) is that it does not require $\|\mathbf{d}_t - \mathbf{w}_t\|$ to be bounded *a priori* (which, although can be done under Proposition 3, is usually loose). A second reason is that almost all examples we are aware of deduce (17) from (16).

The proof is given in Appendix A.5. The significance of this result is that it gives an asymptotic $O(1/H_t)$ rate for convergence to a stationary point. Note that there appears to be a trade-off between the asymptotic rate of $\mathsf{G}(\mathbf{w}_t)$ approaching 0 and the error tolerance in each subproblem (10). Of course there are many admissible choices of the step size, for instance $\eta_t = O(1/t^\beta)$ for any $1/2 < \beta \le 1$ would suffice. A minor advantage of the non-adaptive step size rule over the "optimal" adaptive step size (18) is that it does not require the constant $L$ explicitly; it is merely sufficient that an appropriate $L < \infty$ exists.

Finally, we show how one additional assumption allows a stronger convergence property to be established: namely, convergence in terms of the iterates $\{\mathbf{w}_t\}$, which obviously also implies convergence of the corresponding objective values $\{F(\mathbf{w}_t)\}$ (since $F$ is continuous).

**Assumption 5.** *The underlying space $\mathcal{B}$ is a Hilbert space (with its norm $\|\cdot\|_{\mathcal{H}}$), $\nabla\ell$ is $L$-Lipschitz continuous, and $f$ is $L$-strongly convex, i.e., $f - \frac{L}{2}\|\cdot\|_{\mathcal{H}}^2$ is convex.*

It is easy to observe that, as long as Assumption 3 holds (verified below), Assumption 5 implies both Assumption 2 and Assumption 4. Recall that a sequence $\{\mathbf{w}_t\}$ in $\mathcal{B}$ converges *weakly* to $\mathbf{w}$ if for all $\mathbf{g} \in \mathcal{B}^*$, $\langle\mathbf{w}_t, \mathbf{g}\rangle \to \langle\mathbf{w}, \mathbf{g}\rangle$. As the name suggests, weak convergence is less stringent than the conventional notion (namely $\|\mathbf{w}_t - \mathbf{w}\| \to 0$), but the two coincide when $\mathcal{B}$ has finite dimension.

**Theorem 8.** *Let Assumption 1 and Assumption 5 hold. Assume that $\mathrm{dom}\, f$ is closed, that the subroutine is `Null`, that $F$ has at least one stationary point (i.e. some $\mathbf{w}^\star$ such that $\mathbf{0} \in \partial F(\mathbf{w}^\star)$), and that the step size $\eta_t \in [0,1]$ satisfies $\sum_t \eta_t(1 - \eta_t) = \infty$. Then the iterates $\{\mathbf{w}_t\}$ generated by Algorithm 1 will converge weakly to some stationary point $\mathbf{w}^\star$.*

The proof is given in Appendix A.6. To verify that Assumption 3 is met, note that $\{\mathbf{w}_t\}$, as a weakly convergent sequence, must be bounded, and the strong convexity of $f$ implies super-coerciveness, hence a similar argument as in condition (b) of Proposition 3 establishes Assumption 3. Clearly, any step size rule such that $\eta_t = \Omega(1/t)$ and eventually bounded away from 1 (i.e. for all large $t$, $\eta_t \le c < 1$ for some constant $c$) will suffice here.

### 3.3 Refined Analysis when $\ell$ is Convex

Thus far, we have remained satisfied with determining convergence without emphasizing *rates*, since we have been focusing on generality; that is, we have been investigating the weakest assumptions on $\ell$ and $f$ that ensure convergence is still achieved. To pursue a sharper analysis where specific rates of convergence can be established, we now add the assumption that $\ell$ is *convex*.

**Theorem 9.** *Let Assumption 1 and Assumption 4 hold. Also assume $\ell$ is convex, the subroutine is `Relaxed`, and the subproblem (10) is solved up to some additive error $\varepsilon_t \ge 0$. Then, for any $\mathbf{w} \in \mathrm{dom}\, F$ and $t \ge 0$, Algorithm 1 yields*

$$F(\mathbf{w}_{t+1}) \le F(\mathbf{w}) + \pi_t(1 - \eta_0)(F(\mathbf{w}_0) - F(\mathbf{w})) + \sum_{s=0}^{t} \frac{\pi_t}{\pi_s}\eta_s^2(\varepsilon_s/\eta_s + \frac{L}{2}\|\mathbf{d}_s - \mathbf{w}_s\|^2), \quad (20)$$

where $\pi_t := \prod_{s=1}^{t}(1 - \eta_s)$ with $\pi_0 = 1$. Furthermore, for all $t \geq k \geq 0$, the minimal duality gap $\tilde{\mathsf{G}}_k^t := \min_{k \leq s \leq t} \mathsf{G}(\mathbf{w}_s)$ satisfies

$$\tilde{\mathsf{G}}_k^t \leq \frac{1}{\sum_{s=k}^{t} \eta_s} \left[ F(\mathbf{w}_k) - F(\mathbf{w}_{t+1}) + \sum_{s=k}^{t} \eta_s^2 (\varepsilon_s / \eta_s + \tfrac{L}{2} \| \mathbf{d}_s - \mathbf{w}_s \|^2) \right]. \tag{21}$$

From Theorem 9 we derive the following useful result.

**Corollary 10.** *Under the same assumptions as Theorem 9, also let $\eta_t = 2/(t+2)$, $\varepsilon_t \leq \delta \eta_t / 2$ for some $\delta \geq 0$, and $L_F := \sup_t L \| \mathbf{d}_t - \mathbf{w}_t \|^2$. Then, for all $t \geq 1$, Algorithm 1 yields*

$$F(\mathbf{w}_t) \leq F(\mathbf{w}) + \frac{2(\delta + L_F)}{t+3}, \quad \text{and} \tag{22}$$

$$\tilde{\mathsf{G}}_1^t \leq \frac{3(\delta + L_F)}{t \ln 2} \leq \frac{4.5(\delta + L_F)}{t}. \tag{23}$$

The proofs of Theorem 9 and Corollary 10 can be found in Appendix A.8 and Appendix A.9 respectively. Note that the simple step size rule $\eta_t = 2/(t+2)$ already leads to an $O(1/t)$ bound on the objective value attained.[7] Of course, it is possible to use other step size rules. For instance, both $\eta_s = 1/(s+1)$ and the constant rule $\eta_s = 1 - (t+1)^{1/t}$ lead to an $O(\frac{1+\log t}{t+1})$ rate; see (Freund and Grigas, 2013) for detailed calculations. Similar polynomial-decay rules in (Shamir and Zhang, 2013) can also be used.

**Remark 5.** *One catch in Corollary 10 is that $L_F$ might be infinite. Fortunately, $L_F < \infty$ can be easily ensured by Proposition 3. It is important to emphasize that GCG with the simple step size rule $\eta_t = 2/(t+2)$ requires neither the Lipschitz constant $L$ nor the norm $\|\cdot\|$ to be known (the analysis can therefore adopt the best choices for the problem at hand). Note that, even though the rate does not depend on the initial point $\mathbf{w}_0$ provided $\eta_0 = 1$, by letting $\eta_0 \neq 1$ the bound can be slightly improved (Freund and Grigas, 2013).*

**Remark 6.** *Theorem 9 and Corollary 10 demonstrate the usefulness of the Relaxed subroutine and provide a much simpler approach to establishing the convergence rates achieved by GCG. In particular, let "Algorithm 007" denote Algorithm 1 configured with the "optimal" step size rule (18) and the Null subroutine. Although Theorem 6 establishes that Algorithm 007 indeed converges asymptotically, it has been a nontrivial exercise to establish that it also achieves an $O(1/t)$ rate under the additional assumption that $\ell$ is convex (Frank and Wolfe, 1956; Levitin and Polyak, 1966) (see also (Bach, 2013a) for a slightly sharper constant). However, here we can make the simple observation that the behavior of Algorithm 007 is upper bounded by the behavior of the configuration, "Algorithm 008", given by Algorithm 1 with step size $\eta_t = 2/(t+2)$ and any Relaxed subroutine. Since Corollary 10 establishes an $O(1/t)$ rate for Algorithm 008, we immediately have the same result for Algorithm 007.*

---

7. *Bibliographic remark*: The observation that an $O(1/t)$ rate can be obtained by the simple $\eta_t = 2/(t+2)$ step size rule appears to have been first made by Clarkson (2010) in the setting where $f = \iota_C$ for some specific compact set $C$. A similar rate on the minimal duality gap was also given in (Clarkson, 2010), and extended by (Jaggi, 2013). The particular case when $f$ is strongly convex appeared in (Bach, 2013a).

Interestingly, it turns out that GCG also solves the dual problem in the convex setting

$$\inf_{\mathbf{g}} \left\{ \ell^*(\mathbf{g}) + f^*(-\mathbf{g}) \right\},\tag{24}$$

since when $\ell$ and $f$ are both convex (subject to some mild regularity conditions) we have by the Fenchel-Rockafellar duality (Zălinescu, 2002, Corollary 2.8.5):

$$\inf_{\mathbf{w}} \left\{ \ell(\mathbf{w}) + f(\mathbf{w}) \right\} = -\inf_{\mathbf{g}} \left\{ \ell^*(\mathbf{g}) + f^*(-\mathbf{g}) \right\}.$$

In particular, Bach (2013a) has shown that the averaged gradient $\bar{\mathbf{g}}_T$ automatically solves the dual problem (24) at the rate of $O(1/t)$, provided that the sequences $\{\mathbf{d}_t\}$ and $\{\mathbf{w}_t\}$ generated in Algorithm 1 can be bounded.

**Theorem 11.** *Adopt Assumption 1, and furthermore assume $\ell$ is convex with $L$-Lipschitz continuous gradient $\nabla \ell$, the subroutine is `Null`, and the step size $\eta_t = \frac{2}{t+2}$. Also let $\bar{\mathbf{g}}_{t+1} := \frac{2}{(t+1)(t+2)} \sum_{s=0}^{t} (s+1) \mathbf{g}_s$. Then for all $\mathbf{g}$ and $t \geq 0$, Algorithm 1 yields*

$$\ell^*(\bar{\mathbf{g}}_{t+1}) + f^*(-\bar{\mathbf{g}}_{t+1}) \leq \ell^*(\mathbf{g}) + f^*(-\mathbf{g}) + \frac{2L}{(t+1)(t+2)} \sum_{s=0}^{t} \frac{s+1}{s+2} \|\mathbf{w}_s - \mathbf{d}_s\|^2.\tag{25}$$

Theorem 11 was first observed in (Bach, 2013a) by identifying the iterates with those of a modified mirror descent. We provide a direct proof in Appendix A.7.

### 3.4 Improved Algorithm and Refined Analysis when $f$ is a Gauge

We now present one of our main contributions, which we will exploit in the second part of this paper: an efficient specialization of GCG that is applicable to positively homogeneous convex regularizers; i.e., *gauge* regularizers. Such a specialization is motivated by the fact that regularizers in machine learning are often (semi) norms, hence *bona fide* gauges. In this setting, we develop a specialized form of GCG that exploits properties of gauges to significantly reduce computational overhead. The resulting GCG formulation will be used to achieve efficient algorithms in important case studies, such as matrix completion under trace norm regularization (Section 4.2.2 below).

It is important to first observe that Algorithm 1 *cannot* always be applied in this setting, since the direction subproblem (10) might diverge to $-\infty$ leaving the update in Line 4 of Algorithm 1 undefined. Two main approaches have been used to bypass this unboundedness in the literature; unfortunately, both remain unsatisfactory in different ways. First, one obvious way to restore boundedness to (10) is to reformulate the regularized problem (2) in its equivalent constrained form

$$\left\{ \inf_{\mathbf{w}} \ \ell(\mathbf{w}) \ \text{ s.t. } \ f(\mathbf{w}) \leq \zeta \right\} \ = \ \inf_{\mathbf{w}} \left\{ \ell(\mathbf{w}) + \iota_{\{\mathbf{w} : f(\mathbf{w}) \leq \zeta\}}(\mathbf{w}) \right\},\tag{26}$$

where it is well-known that a correspondence between (2) and (26) can be achieved if the constant $\zeta$ is chosen appropriately. In this case, Algorithm 1 can be directly applied as long as (26) has a bounded domain. Much recent work in machine learning has investigated this variant (Clarkson, 2010; Hazan, 2008; Jaggi and Sulovsky, 2010; Jaggi, 2013;

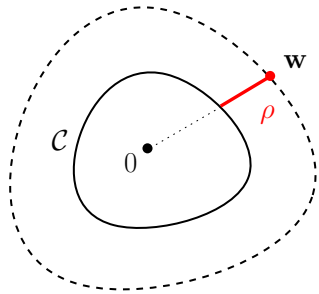

Figure 2: Defining a gauge via the Minkowski functional.

Shalev-Shwartz et al., 2010; Tewari et al., 2011; Yuan and Yan, 2013). However, finding the appropriate value of $\zeta$ is costly in general, and often one cannot avoid considering the penalized formulation (e.g. when it is nested in another problem). Furthermore, the constraints in (26) preclude the application of many efficient local `Improve` techniques designed for unconstrained problems.

A second, alternative, approach is to square $f$, which ensures that it becomes supercoercive (Bradley and Bagnell, 2009). Although squaring indeed works in a finite dimensional space when $f$ is a norm, it can fail to preserve solution equivalence in infinite dimensions (since there not all norms are equivalent). Moreover, the insertion of a local `Improve` step in such an approach requires the regularizer $f$ to be evaluated at all iterates, which is expensive for applications like matrix completion. Furthermore, squaring is a somewhat arbitrary modification, in the sense that supercoerciveness can be achieved (in finite dimensional spaces) by raising the norm to the $p$th power for any $p \geq 2$.

### 3.4.1 Generalized Gauge Regularization and Minkowski Functional

To address the issue of ensuring well defined iterates, while also reducing their cost, we develop an alternative modification of GCG. The algorithm we develop considers a mild generalization of gauge regularization by introducing a composition with an increasing convex function. In particular, let $\kappa$ denote a closed and convex gauge and let $h : \mathbb{R}_+ \to \overline{\mathbb{R}}$ denote an increasing convex function. Then regularizing by $f(\cdot) = h(\kappa(\cdot))$ yields the particular form of (2) we will address with this modified approach:

$$\inf_{\mathbf{w}} \left\{ \ell(\mathbf{w}) + h(\kappa(\mathbf{w})) \right\}. \tag{27}$$

To achieve an efficient GCG algorithm for (27), we exploit the fact that a gauge can be reconstructed from its "unit ball" via the Minkowski functional; that is, letting $\mathcal{C} := \{\mathbf{w} \in \mathcal{B} : \kappa(\mathbf{w}) \leq 1\}$ (which must be a closed and convex set), the gauge $\kappa$ must satisfy

$$\kappa(\mathbf{w}) = \inf\{\rho : \mathbf{w} \in \rho \, \mathcal{C}\}. \tag{28}$$

Intuitively, (28) states that $\kappa(\mathbf{w})$ is the smallest rescaling of the set $\mathcal{C}$ that contains $\mathbf{w}$; see Figure 2. A related function that underpins our algorithm is its *polar*, defined as

$$\kappa^{\circ}(\mathbf{g}) := \sup_{\mathbf{w} \in \mathcal{C}} \langle \mathbf{w}, \mathbf{g} \rangle. \tag{29}$$

(which, for example, would be the dual norm if $\kappa$ is a norm). For the structured sparse regularizers typically considered in machine learning, the "unit ball" $\mathcal{C}$ has additional structure that admits efficient computation (Chandrasekaran et al., 2012). In particular, $\mathcal{C}$ is often constructed by taking the (closed) convex hull of some compact set of "atoms" $\mathcal{A}$ specific to the problem; i.e., $\mathcal{C} = \mathrm{conv}\,\mathcal{A}$. Such a structured formulation allows one to re-express the gauge and polar, (28) and (29), respectively as

$$\kappa(\mathbf{w}) = \inf\left\{\rho : \mathbf{w} = \rho\sum_i \sigma_i \mathbf{a}_i, \sigma_i \geq 0, \sum_i \sigma_i = 1, \mathbf{a}_i \in \mathcal{A}\right\},^8 \tag{30}$$

$$\kappa^\circ(\mathbf{g}) = \sup_{\mathbf{a}\in\mathcal{A}} \langle \mathbf{a}, \mathbf{g}\rangle. \tag{31}$$

Intuitively, $\kappa^\circ$ has a much simpler form than $\kappa$ if the structure of $\mathcal{A}$ is "simple". In such cases, it is advantageous to design an optimization algorithm that circumvents evaluation of $\kappa$, by instead using $\kappa^\circ$ which is generally far less expensive to evaluate.

### 3.4.2 MODIFIED GCG ALGORITHM FOR GENERALIZED GAUGE REGULARIZATION

The key to the modified algorithm is to avoid the possibly unbounded solution to (10) while also avoiding explicit evaluation of the gauge $\kappa$. First, to avoid unboundedness, we adopt the technique of moving the regularizer to the constraint, as in (26), but rather than applying the modification directly to (27) we only use it in the direction problem (10) via:

$$\mathbf{d}_t \in \arg\min_{\mathbf{d}:h(\kappa(\mathbf{d}))\leq\zeta} \langle \mathbf{d}, \nabla\ell(\mathbf{w}_t)\rangle. \tag{32}$$

Then, to bypass the complication of dealing with $h$ and the unknown bound $\zeta$, we decompose $\mathbf{d}_t$ into its normalized direction $\mathbf{a}_t$ and scale $\theta_t$ (i.e., such that $\mathbf{d}_t = \mathbf{a}_t\theta_t/\eta_t$), determining each separately. In particular, from this decomposition, the update from Line 6 of Algorithm 1 can be equivalently expressed as

$$\mathbf{w}_{t+1} = (1-\eta_t)\mathbf{w}_t + \eta_t\mathbf{d}_t \tag{33}$$

$$= (1-\eta_t)\mathbf{w}_t + \theta_t\mathbf{a}_t, \tag{34}$$

hence the modified algorithm need only work with $\mathbf{a}_t$ and $\theta_t$ directly.

**Determining the direction** First, the normalized direction, $\mathbf{a}_t$, can be recovered via

$$\mathbf{a}_t \in \arg\min_{\mathbf{a}:\kappa(\mathbf{a})\leq 1} \langle \mathbf{a}, \nabla\ell(\mathbf{w}_t)\rangle \tag{35}$$

$$\Leftrightarrow \quad \mathbf{a}_t \in \arg\min_{\mathbf{a}\in\mathcal{A}} \langle \mathbf{a}, \nabla\ell(\mathbf{w}_t)\rangle \qquad (\text{by } (30)) \tag{36}$$

$$\Leftrightarrow \quad \mathbf{a}_t \in \arg\max_{\mathbf{a}\in\mathcal{A}} \langle \mathbf{a}, -\nabla\ell(\mathbf{w}_t)\rangle. \tag{37}$$

Observe that (37) effectively involves computation of the polar $\kappa^\circ(-\nabla\ell(\mathbf{w}))$ only, by (31). Since this optimization might still be challenging, we further allow it to be solved only

---

8. Note that although $\mathcal{A}$ can have an infinite or even uncountable cardinality in this construction, the summation in (30) is always well-defined since $\{\sigma_i\}$ is absolutely summable; that is, for any $\mathbf{w} \in \mathcal{B}$, (30) can be expressed by taking a sequence $\{\sigma_i\}$ with at most countably many nonzero entries.

---

**Algorithm 2** GCG for positively homogeneous regularizers.

---

**Require:** The set $\mathcal{A}$ whose (closed) convex hull $\mathcal{C}$ defines the gauge $\kappa$.
1: Initialize $\mathbf{w}_0$ and $\rho_0 \geq \kappa(\mathbf{w}_0)$.
2: **for** $t = 0, 1, \ldots$ **do**
3:    Choose normalized direction $\mathbf{a}_t$ that satisfies (38)
4:    Choose step size $\eta_t \in [0, 1]$ and set the scaling $\theta_t \geq 0$ by (41)
5:    $\tilde{\mathbf{w}}_{t+1} = (1 - \eta_t)\mathbf{w}_t + \theta_t \mathbf{a}_t$
6:    $\tilde{\rho}_{t+1} = (1 - \eta_t)\rho_t + \theta_t$
7:    $(\mathbf{w}_{t+1}, \rho_{t+1}) = \texttt{Improve}(\tilde{\mathbf{w}}_{t+1}, \tilde{\rho}_{t+1}, \ell, f)$    ▷ Subroutine, see Theorem 12
8: **end for**

---

*approximately* to within an additive error $\varepsilon_t \geq 0$ and a multiplicative factor $\alpha_t \in (0, 1]$; that is, by relaxing (36) the modified algorithm only requires an $\mathbf{a}_t \in \mathcal{A}$ to be found that satisfies

$$\langle \mathbf{a}_t, \nabla\ell(\mathbf{w}_t) \rangle \leq \alpha_t \left( \varepsilon_t + \min_{\mathbf{a} \in \mathcal{A}} \langle \mathbf{a}, \nabla\ell(\mathbf{w}_t) \rangle \right) = \alpha_t \left( \varepsilon_t - \kappa^\circ(-\nabla\ell(\mathbf{w}_t)) \right). \tag{38}$$

Intuitively, this formulation computes the normalized direction that has (approximately) the greatest negative "correlation" with the gradient $\nabla\ell$, yielding the steepest local decrease in $\ell$. Importantly, this formulation does not require evaluation of the gauge function $\kappa$.

**Determining the scale**   Next, to chose the scale $\theta_t$, it is natural to consider minimizing a simple upper bound on the objective that is obtained by replacing $\ell(\cdot)$ with its quadratic upper approximation $\tilde{\ell}_L(\cdot; \mathbf{w}_t)$ given in (1):

$$\min_{\theta \geq 0} \left\{ \tilde{\ell}_L((1 - \eta_t)\mathbf{w}_t + \theta\mathbf{a}_t; \mathbf{w}_t) + h(\kappa((1 - \eta_t)\mathbf{w}_t + \theta\mathbf{a}_t)) \right\}. \tag{39}$$

Unfortunately, $\kappa$ still participates in (39), so a further approximation is required. Here is where the decomposition of $\mathbf{d}_t$ into $\mathbf{a}_t$ and $\theta_t$ is particularly advantageous. Observe that

$$\begin{aligned} h(\kappa((1 - \eta_t)\mathbf{w}_t + \theta\mathbf{a}_t)) &\leq h((1 - \eta_t)\kappa(\mathbf{w}_t) + \eta_t\kappa(\theta\mathbf{a}_t/\eta_t)) \\ &\leq (1 - \eta_t)h(\kappa(\mathbf{w}_t)) + \eta_t h(\kappa(\theta\mathbf{a}_t/\eta_t)) \\ &\leq (1 - \eta_t)h(\kappa(\mathbf{w}_t)) + \eta_t h(\theta/\eta_t)), \end{aligned} \tag{40}$$

where the first inequality follows from the convexity of $\kappa$ and the fact that $h$ is increasing, the second inequality follows from the convexity of $h$, and (40) follows from the fact that $\kappa(\mathbf{a}_t) \leq 1$ by construction ($\mathbf{a}_t \in \mathcal{A}$ in (30)). Although $\kappa$ still participates in (40) we can now bypass evaluation of $\kappa$ by maintaining a simple upper bound $\rho_t \geq \kappa(\mathbf{w}_t)$, yielding the relaxed update

$$\theta_t = \arg\min_{\theta \geq 0} \left\{ \tilde{\ell}_L((1 - \eta_t)\mathbf{w}_t + \theta\mathbf{a}_t; \mathbf{w}_t) + (1 - \eta_t)h(\rho_t) + \eta_t h(\theta/\eta_t) \right\}. \tag{41}$$

Crucially, the upper bound $\rho_t \geq \kappa(\mathbf{w}_t)$ can be maintained by the simple update $\rho_{t+1} = (1 - \eta_t)\rho_t + \theta_t$, which ensures

$$\rho_{t+1} = (1 - \eta_t)\rho_t + \theta_t \geq (1 - \eta_t)\kappa(\mathbf{w}_t) + \theta_t\kappa(\mathbf{a}_t) \geq \kappa((1 - \eta_t)\mathbf{w}_t + \theta_t\mathbf{a}_t) \geq \kappa(\mathbf{w}_{t+1}). \tag{42}$$

From these components we obtain the final modified procedure, Algorithm 2: a variant of GCG (Algorithm 1) that avoids the potential unboundedness of the subproblem (10) by replacing it with (38), while also completely avoiding explicit evaluation of the gauge $\kappa$. For problems of the form (27), Algorithm 2 can provide significantly more efficient updates than Algorithm 1. Despite the relaxations introduced in Algorithm 2 to achieve faster iterates, the following theorem establishes that the procedure is still sound, preserving essentially the same convergence guarantees as the more expensive Algorithm 1.

**Theorem 12.** *Let Assumption 4 hold and assume that $\ell$ is convex. Let $h : \mathbb{R}_+ \to \overline{\mathbb{R}}$ be an increasing convex function and $\kappa$ be a closed gauge induced by the atomic set $\mathcal{A}$. Let $\alpha_t > 0$, $0 \leq \eta_t \leq 1$, and the subroutine `Improve` be `Relaxed` in the sense akin to (41):*

$$\begin{cases} \ell(\mathbf{w}_{t+1}) + h(\rho_{t+1}) \leq \tilde{\ell}_L(\tilde{\mathbf{w}}_{t+1}; \mathbf{w}_t) + (1 - \eta_t)h(\rho_t) + \eta_t h(\theta_t/\eta_t) \\ \quad\quad and \quad \rho_{t+1} \geq \kappa(\mathbf{w}_{t+1}) \end{cases}. \tag{43}$$

*Then for any $\mathbf{w} \in \operatorname{dom} F$ and $t \geq 0$, Algorithm 2 yields*

$$F(\mathbf{w}_{t+1}) \leq F(\mathbf{w}) + \pi_t(1 - \eta_0)\big(F(\mathbf{w}_0) - F(\mathbf{w})\big)$$
$$+ \sum_{s=0}^{t} \frac{\pi_t}{\pi_s} \eta_s^2 \Big( \big(\rho\varepsilon_s + h(\rho/\alpha_s) - h(\rho)\big)/\eta_s + \frac{L}{2}\left\| \frac{\rho}{\alpha_s}\mathbf{a}_s - \mathbf{w}_s \right\|^2 \Big), \tag{44}$$

*where $\rho := \kappa(\mathbf{w})$ and $\pi_t := \prod_{s=1}^{t}(1 - \eta_s)$ with $\pi_0 = 1$.*

The proof is given in Appendix A.10. From this theorem, we obtain the following rate.

**Corollary 13.** *Under the same setting as in Theorem 12, let $\eta_t = 2/(t+2)$, $\varepsilon_t \leq \delta\eta_t/2$ for some absolute constant $\delta > 0$, $\alpha_t = \alpha > 0$, $\rho = \kappa(\mathbf{w})$, and $L_F := L \sup_t \left\| \frac{\rho}{\alpha}\mathbf{a}_t - \mathbf{w}_t \right\|^2$. Then for all $t \geq 1$, Algorithm 2 yields*

$$F(\mathbf{w}_t) \leq F(\mathbf{w}) + \frac{2(\rho\delta + L_F)}{t+3} + h(\rho/\alpha) - h(\rho). \tag{45}$$

*Moreover, if $\ell \geq 0$ and $h = \mathsf{Id}$, then*

$$F(\mathbf{w}_t) \leq F(\mathbf{w})/\alpha + \frac{2(\rho\delta + L_F)}{t+3}. \tag{46}$$

Some remarks concerning this algorithm and its convergence properties are in order.

**Remark 7.** *When $\alpha_t = 1$ for all $t$, and $h$ is the indicator $\iota_{[0,\zeta]}$ for some $\zeta > 0$, Corollary 13 implies the same convergence rate for the constrained problem (26), which immediately recovers (some of) the more specialized results discussed in (Clarkson, 2010; Hazan, 2008; Jaggi and Sulovsky, 2010; Jaggi, 2013; Shalev-Shwartz et al., 2010; Tewari et al., 2011; Yuan and Yan, 2013).*

**Remark 8.** *The additional factor $\rho$ in the above bounds is necessary because the approximation in (38) is not invariant to scaling, requiring some compensation. Note also that $\alpha \leq 1$, since the right-hand side of (38) must become negative as $\varepsilon_t \to 0$. The result in (46) states roughly that an $\alpha$-approximate subroutine (for computing the polar of $\kappa$) leads to an $\alpha$-approximate "minimizer",[9] also at the rate of $O(1/t)$.*

---

9. Independent of our work, which was documented earlier in (Yu, 2013b, Section 4.1.4), Bach (2013b) has also considered a similar multiplicative approximate polar oracle.

**Remark 9.** *Again, $L_F < \infty$ can be guaranteed by Proposition 3. Assumption 4 is implicitly required to ensure the* `Relaxed` *condition is satisfiable. For instance, the rule*

$$\mathbf{w}_{t+1} := \tilde{\mathbf{w}}_{t+1} = (1 - \eta_t)\mathbf{w}_t + \theta_t \mathbf{a}_t, \quad \rho_{t+1} := \tilde{\rho}_{t+1} = (1 - \eta_t)\rho_t + \theta_t, \tag{47}$$

*which we call* `Null`, *suffices. The update for $\theta$ in (41) requires knowledge of the Lipschitz constant $L$ and the norm $\|\cdot\|$. If the norm $\|\cdot\|$ is Hilbertian and $h = \mathsf{Id}$, we have the formula*

$$\theta_t = \left( \frac{\langle \mathbf{a}_t, L\eta_t \mathbf{w}_t - \nabla \ell(\mathbf{w}_t) \rangle - 1}{L \|\mathbf{a}_t\|^2} \right)_+.$$

*However, it is also easy to devise other scaling and step size updates that do not require the knowledge of $L$ or the norm $\|\cdot\|$, such as*

$$(\eta_t^*, \theta_t^*) \in \arg\min_{1 \geq \eta \geq 0, \theta \geq 0} \left\{ \ell((1 - \eta)\mathbf{w}_t + \theta\mathbf{a}_t) + (1 - \eta)h(\rho_t) + \eta h(\theta/\eta) \right\}, \tag{48}$$

*followed by taking (47). Note that (48) is jointly convex in $\eta$ and $\theta$, since $\eta h(\theta/\eta)$ is a perspective function. Evidently this update can be treated as a* `Relaxed` *subroutine, since for any $(\eta_t, \theta_t)$ chosen in step 4 of Algorithm 2, it follows that*

$$\begin{aligned}
\ell(\mathbf{w}_{t+1}) + h(\rho_{t+1}) &\leq \ell((1 - \eta_t^*)\mathbf{w}_t + \theta_t^*\mathbf{a}_t) + (1 - \eta_t^*)h(\rho_t) + \eta_t^* h(\theta_t^*/\eta_t^*) \\
&\leq \ell((1 - \eta_t)\mathbf{w}_t + \theta_t\mathbf{a}_t) + (1 - \eta_t)h(\rho_t) + \eta_t h(\theta_t/\eta_t) \quad (by \ (48)) \\
&\leq \tilde{\ell}_L(\tilde{\mathbf{w}}_{t+1}; \mathbf{w}_t) + (1 - \eta_t)h(\rho_t) + \eta_t h(\theta_t/\eta_t).
\end{aligned}$$

*The upper bound $\rho_{t+1} \geq \kappa(\mathbf{w}_{t+1})$ can be verified in the same way as for (42); hence the update (48) enjoys the same convergence guarantee as in Corollary 13.*

**Remark 10.** *An alternative workaround for solving (2) with gauge regularization would be to impose an upper bound $\zeta \geq \kappa(\mathbf{w})$ and consider*

$$\min_{\mathbf{w}:\kappa(\mathbf{w}) \leq \zeta} \left\{ \ell(\mathbf{w}) + \kappa(\mathbf{w}) \right\}, \tag{49}$$

*while dynamically adjusting $\zeta$. A standard GCG approach would then require solving*

$$\mathbf{d}_t \in \arg\min_{\mathbf{d}:\kappa(\mathbf{d}) \leq \zeta} \left\{ \langle \mathbf{d}, \nabla \ell(\mathbf{w}_t) \rangle + \kappa(\mathbf{d}) \right\} \tag{50}$$

*as a subproblem, which can be as easy to solve as (38). If, however, one wished to include a local improvement step* `Improve` *(which is essential in practice), it is not clear how to efficiently maintain the constraint $\kappa(\mathbf{w}) \leq \zeta$. Moreover, solving (50) up to some multiplicative factor might* not *be as easy as (38).*

Surprisingly, there is little work on investigating the gauge regularized problem (27). One exception is Dudík et al. (2012), who propose a totally corrective variant (cf. (51) below). However, the analysis in Dudík et al. (2012) only addresses the diminishment of the gradient, and establishes a suboptimal $O(1/\sqrt{t})$ rate.[10]

---

10. Shortly after the appearance of our preliminary work Zhang et al. (2012), Harchaoui et al. (2014) established a similar $O(1/t)$ rate.

### 3.5 Additional Discussion

To establish convergence rates for the convex case we have assumed the loss $\ell$ has a Lipschitz continuous gradient. Although this holds for a variety of losses such as the square and logistic loss, it does not hold for the hinge loss used in support vector machines. Although one can always approximate a nonsmooth loss by a smooth function (Nesterov, 2005) before applying GCG, this usually results in a slower $O(1/\sqrt{t})$ rate of convergence.

Several of the convergence proofs rely on the step size rule $\eta_t = O(1/t)$, which appears to be "optimal" among non-adaptive rules in the following sense. On the one hand, large steps often lead to faster convergence, while on the other hand, Algorithm 2 needs to remove atoms $\mathbf{a}_t$ that are perhaps "incidentally" added, which requires the discount factor $\prod_{t=1}^{\infty}(1 - \eta_t)$ to be small. It is an easy exercise to prove that the latter condition holds iff $\sum_{t=1}^{\infty} \eta_t = \infty$, thus the step size $O(1/t)$ is (almost) the largest non-adaptive rule that still allows irrelevant atom removal.

As given, Algorithm 2 adds a single atom in each iteration, with a conic combination between the new atom and the previous aggregate. An even more aggressive scheme is to re-optimize the weights of each individual atom in each iteration. Such a procedure was first studied in (Meyer, 1974), which is often known as the "totally (or fully) corrective update" in the boosting literature. For example, given a gauge $\kappa$ and $h = \mathsf{Id}$ each iteration of a totally corrective procedure would involve solving

$$\min_{\boldsymbol{\sigma} \geq 0} \left\{ \ell \left( \sum_{\tau=1}^{t} \sigma_\tau \mathbf{a}_\tau \right) + \sum_{\tau=1}^{t} \sigma_\tau \right\}. \tag{51}$$

Not surprisingly, the totally corrective variant can be seen as a `Relaxed` subroutine in Algorithm 2, hence converges at least as fast as Corollary 13 suggests. Much faster convergence is usually observed in practice, although this advantage must be countered by the extra effort spent in solving (51), which itself need not be trivial. In a finite dimensional setting, provided that the atoms are linearly independent and some restricted strong convexity is present, it is possible to prove that the totally corrective variant (51) converges at a linear rate; see (Shalev-Shwartz et al., 2010; Yuan and Yan, 2013; Zhang et al., 2012).

Finally, we remark that the convergence rate established for GCG is on par with PG, and cannot be improved even in the presence of strong convexity (Canon and Cullum, 1968). In this sense, GCG is slower than an "optimal" algorithm like APG. GCG's advantage, however, is that it only needs to solve a linear subproblem (10) in each iteration (i.e., computing the polar $\kappa^{\circ}$), whereas APG (or PG) needs to compute the proximal update, which is a quadratic problem for least square prox-function. The two approaches are complementary in the sense that the polar can be easier to compute in some cases, whereas in other cases the proximal update can be computed analytically. An additional advantage GCG has over APG, however, is its greedy sparse nature: Each iteration of GCG amounts to a single atom, hence the number of atoms never exceeds the total number of iterations for GCG. By contrast, APG can produce a dense update in each iteration, although in later stages the estimates may become sparser due to the shrinkage effect of the proximal update. Moreover, GCG is "robust" with respect to $\alpha$-approximate atom selection (cf. (46) in Corollary 13), whereas similar results do not appear to be available for PG or APG when the proximal update is computed up to a multiplicative factor.

## 4. Application to Low Rank Learning

As a first application of GCG, we consider the problem of optimizing low rank matrices; a problem that forms the cornerstone of many machine learning tasks such as dictionary learning and matrix completion. Since imposing a bound on matrix rank generally leads to an intractable problem, much recent work has investigated approaches that relax the hard rank constraint with carefully designed regularization penalties that indirectly encourage low rank structure (Bach et al., 2008; Lee et al., 2009; Bradley and Bagnell, 2009). In Section 4.2 we first formulate a generic convex relaxation scheme for low rank problems arising in matrix approximation tasks, including dictionary learning and matrix completion. Conveniently, the modified GCG algorithm developed in Section 3.4, Algorithm 2, provides a general solution method for the resulting problems. Next, in Section 4.3 we show how the local improvement component of the modified GCG algorithm can be provided in this case by a procedure that optimizes an unconstrained surrogate (instead of a constrained problem) to obtain greater efficiency. This modification significantly improves the practical efficiency of GCG without affecting its convergence properties. We illustrate these contributions by highlighting the example of matrix completion under trace norm regularization, where the improved GCG algorithm demonstrates state of the art performance. Finally, we demonstrate how the convex relaxation can be generalized to handle latent multi-view models (White et al., 2012; Lee et al., 2009) in Section 4.4, showing how the improved GCG algorithm can still be applied once the major challenge of efficiently computing the dual norm (polar $\kappa^\circ$) has been solved.

### 4.1 Low Rank Learning Problems

Many machine learning tasks, such as dictionary learning and matrix completion, can be formulated as seeking low rank matrix approximations of a given data matrix. In these problems, one is presented with an $n \times m$ matrix $X$ (perhaps only partially observed) where each column corresponds to a training example and each row corresponds to a feature across examples; the goal is to recover a low rank matrix $W$ that approximates (the observed entries of) $X$. To achieve this, we assume we are given a convex loss function $\ell(W) := L(W, X)$ that measures how well $W$ approximates $X$.[11] If we impose a hard bound, $t$, on the rank of $W$, the core training problem can then be expressed as

$$\inf_{W:\text{rank}(W) \le t} \ell(W) = \inf_{U \in \mathbb{R}^{n \times t}} \inf_{V \in \mathbb{R}^{t \times m}} \ell(UV). \tag{52}$$

Unfortunately, this problem is not convex, and except for special cases (such as SVD) is believed to be intractable.

The two forms of this problem that we will explicitly illustrate in this paper are dictionary learning and matrix completion. First, for dictionary learning (Olshausen and Field, 1996), the $n \times m$ data matrix $X$ is fully observed, the $n \times t$ matrix $U$ is interpreted as a "dictionary" consisting of $t$ basis vectors, and the $t \times m$ matrix $V$ is interpreted as the "coefficients" consisting of $m$ code vectors. For this problem, the dictionary and coefficient matrices need to be inferred *simultaneously*, in contrast to traditional signal approximation schemes where the dictionary (e.g., the Fourier or a wavelet basis) is fixed *a priori*. Unfortunately, learning the dictionary with the coefficients creates significant computational challenge since the formulation is no longer *jointly* convex in the variables $U$ and $V$, even when the loss $\ell$ is convex. Indeed, for a fixed dictionary size $t$, the problem is known to be computationally tractable only for losses induced by unitarily invariant norms (Yu and Schuurmans, 2011). With nonnegative constraints on $U$ and $V$ (Lee and Seung, 1999) the problem is NP-hard even for the squared loss (Vavasis, 2010).

---

11. Note that the matrix $W$ need not reconstruct the entries of $X$ directly: a nonlinear transfer could be interposed while maintaining convexity of $\ell$; for example, by using an appropriate Bregman divergence.

The second learning problem we illustrate is low rank matrix completion. Here, one observes a small number of entries in $X \in \mathbb{R}^{n \times m}$ and the task is to infer the remaining entries. In this case, the optimization objective can be expressed as $\ell(W) := L(\mathcal{P}(W - X))$, where $L$ is a reconstruction loss (such as Frobenius norm squared) and $\mathcal{P} : \mathbb{R}^{n \times m} \to \mathbb{R}^{n \times m}$ is a "mask" operator that simply fills the entries corresponding to unobserved indices in $X$ with zero. Various recommendation problems, such as the Netflix prize,[12] can be cast as matrix completion. Due to the ill-posed nature of the problem it is standard to assume that the matrix $X$ can be well approximated by a low rank reconstruction $W$. Unfortunately, the rank function is not convex, hence hard to minimize in general. Conveniently, both the dictionary learning and matrix completion problems are subsumed by the general formulation we consider in (52).

## 4.2 General Convex Relaxation and Solution via GCG

To develop a convex relaxation of (52), we first consider the factored formulation expressed in terms of $U$ and $V$ on the right hand side, which only involves an unconstrained minimization. Although unconstrained minimization can be convenient, care is needed to handle the scale invariance introduced by the factored form: since $U$ can always be scaled and $V$ counter-scaled to preserve the product $UV$, some form of penalty or constraint is required to restore well-posedness. A standard strategy is to constrain $U$, for which a particularly common choice is to constrain its columns to have unit length; i.e. $\|U_{:i}\|_c \leq 1$ for all $i$, where $c$ simply denotes some norm on the columns. (For example, for dictionary learning, this amounts to constraining the basis vectors in the dictionary $U$ to lie within the unit ball of $\|\cdot\|_c$.) Of course, even given scale-invariance the problem remains non-convex and current work on dictionary learning and matrix completion still often resorts to local minimization.

However, it has recently been observed that a simple relaxation of the rank constraint in (52) allows a convex reformulation to be achieved (Argyriou et al., 2008; Bach et al., 2008; Bradley and Bagnell, 2009; Zhang et al., 2011; Haeffele et al., 2014). First, observe that replacing the hard constraint with a regularization penalty on the magnitude of rows in the coefficient matrix $V$ yields a relaxed form of (52):

$$\inf_{U:\|U_{:i}\|_c \leq 1} \inf_{\tilde{V}} \left\{ \ell(U\tilde{V}) + \lambda \sum_i \|\tilde{V}_{i:}\|_r \right\}, \tag{53}$$

where $U_{:i}$ and $\tilde{V}_{i:}$ denote the $i$th column of $U$ and $i$th row of $\tilde{V}$ respectively, $\|\cdot\|_r$ denotes some norm on the rows, and $\lambda \geq 0$ balances the trade-off between the loss and the regularizer. Here, the number of columns in $U$ (and also the number of rows in $\tilde{V}$) is not restricted;

---

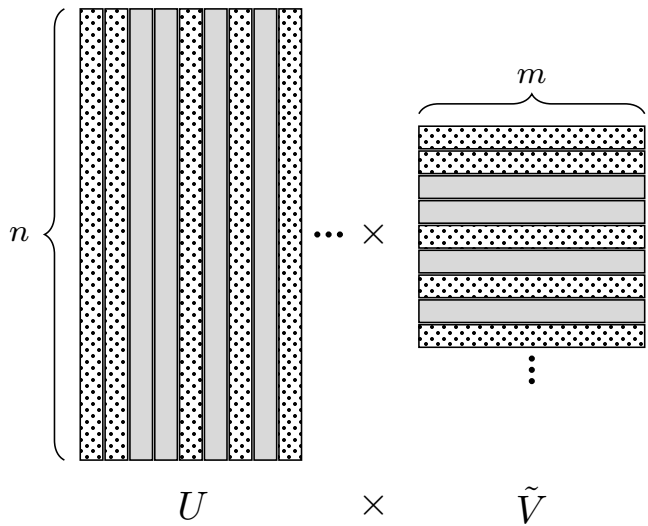

Figure 3: The idea behind the convex relaxation for dictionary learning: Due to the row-wise norm regularizer on $\tilde{V}$, its many rows (grayed) will be zeroed out, therefore the corresponding columns of $U$ will be dropped, resulting in a small dictionary.

instead, by minimizing the regularized problem (53) many of the rows in $\tilde{V}$ will become zero due to the effect of the row-norm regularizer; see Figure 3. Hence, the corresponding columns of $U$ can be dropped, implying that the rank in (53) is chosen adaptively; a scheme that can be motivated by the fact that one usually does not know the exact rank *a priori*. The specific form of the row norm $\|\cdot\|_r$ also provides additional flexibility in promoting different structures; for example, the $l_1$ norm leads to sparse solutions, the $l_2$ norm yields low rank solutions, and block structured norms generate group sparsity. The specific form of the column norm $\|\cdot\|_c$ also has an effect, as discussed in Section 4.4 below.

Importantly, although the problem (53) is still not jointly convex in the factors $U$ and $\tilde{V}$, it is known that it can be exactly reformulated as a convex optimization via a simple change of variables $W := U\tilde{V}$; an observation that has been made, for example, by (Argyriou et al., 2008; Bach et al., 2008; Zhang et al., 2011; White et al., 2012). A key development in these earlier works was to infer the structure of an implicit (tractable) regularizer on $W$ implied by the constraint on $U$ and regularizer on $\tilde{V}$. However, rather than repeat these earlier derivations, we provide an alternative derivation that uses gauge regularization to achieve a simpler and more flexible view (Zhang et al., 2012). Not only does this alternative approach yield more concise proofs, it also establishes a direct connection to the GCG method outlined in Section 3.4, yielding an efficient solution method with no additional effort.

### 4.2.1 Convex Reformulation via Gauge Regularization

The following treatment was originally developed in (Zhang et al., 2012), which we now present with some clarifications. Starting with (53), consider a reparameterization $\tilde{V}_{i:} = \sigma_i V_{i:}$, where $\sigma_i \geq 0$ and $\|V_{i:}\|_r \leq 1$. This change of variable allows (53) to be equivalently expressed, using $W = U\tilde{V}$, as:

$$(53) = \min_W \ \ell(W) + \lambda \ \inf \left\{ \sum_i \left\| \tilde{V}_{i:} \right\|_r : \left\| U_{:i} \right\|_c \leq 1, U\tilde{V} = W \right\} \tag{54}$$

$$= \min_W \ \ell(W) + \lambda \ \inf \left\{ \sum_i \sigma_i : \boldsymbol{\sigma} \geq 0, W = \sum_i \sigma_i U_{:i} V_{i:}, \left\| U_{:i} \right\|_c \leq 1, \left\| V_{i:} \right\|_r \leq 1 \right\} \tag{55}$$

$$= \min_W \ \ell(W) + \lambda \ \kappa(W), \tag{56}$$

where the gauge $\kappa$ is defined by the set

$$\mathcal{A} := \{ \mathbf{u}\mathbf{v}^\top : \mathbf{u} \in \mathbb{R}^n, \mathbf{v} \in \mathbb{R}^m, \|\mathbf{u}\|_c \leq 1, \|\mathbf{v}\|_r \leq 1 \} \tag{57}$$

via the Minkowski functional (30). Note that the set $\mathcal{A}$ contains uncountably many elements. Clearly (56) is a convex problem, provided that the loss $\ell$ is convex.

The reformulation (56) illuminates a few important aspects. First, it reveals that the relaxed learning problem (53) is equivalent to considering a rank-one decomposition of $W$, where penalization by the sum of "singular values" $\{\sigma_i\}$ provides a convex proxy for rank. From this perspective, it is clear the formulation (55) generalizes the singular value decomposition. Second, this formulation now reveals a clearer understanding of the computational properties of the regularizer in (56) via its polar:

$$\kappa^\circ(G) = \sup_{A \in \mathcal{A}} \ \langle A, G \rangle = \sup_{\|\mathbf{u}\|_c \leq 1, \|\mathbf{v}\|_r \leq 1} \mathbf{u}^\top G \mathbf{v} \tag{58}$$

$$= \sup_{\|\mathbf{u}\|_c \leq 1} \ \left\| G^\top \mathbf{u} \right\|_r^\circ \tag{59}$$

$$= \sup_{\|\mathbf{v}\|_r \leq 1} \ \|G\mathbf{v}\|_c^\circ . \tag{60}$$

Recall that for gauge regularized problems in the form (56), Algorithm 2 only requires approximate computation of the polar (38) in Line 3; hence an efficient method for (approximately) solving (58) immediately yields an efficient implementation. However, this simplicity must be countered by the fact that the polar is not always tractable, since it involves *maximizing* a norm subject to a different norm constraint.[13] Nevertheless, there exist important cases where the polar can be efficiently evaluated, which we illustrate below.

The polar also allows one to gain insight into the structure of the original gauge regularizer, since by duality $\kappa = (\kappa^\circ)^\circ$ yields an explicit formula for $\kappa$. For example, consider the important special case where $\|\cdot\|_r = \|\cdot\|_c = \|\cdot\|_2$. In this case, using (58) one can infer that the polar $\kappa^\circ(\cdot) = \|\cdot\|_{\mathrm{sp}}$ is the spectral norm, and therefore the corresponding gauge regularizer $\kappa(\cdot) = \|\cdot\|_{\mathrm{tr}}$ is the trace norm. In this way, the trace norm regularizer, which is known to provide a convex relaxation of the matrix rank (Fazel et al., 2001; Chandrasekaran et al., 2012), can be viewed from the perspective of dictionary learning as an *induced* norm that arises from a simple 2-norm constraint on the magnitude of the dictionary entries $U_{:i}$ and a simple 2-norm penalty on the magnitudes of the code vectors $V_{i:}$.

---

13. In this regard, the $\alpha$-approximate polar oracle introduced in Algorithm 2 (38) and Corollary 13 might be useful, although we shall not develop this idea further here.

This convex relaxation strategy based on gauge functions is quite flexible and has been recently studied for some structured sparse problems (Chandrasekaran et al., 2012; Tewari et al., 2011; Zhang et al., 2012). Importantly, the generality of the column and row norms in (53) proffer considerable flexibility to finely characterize the structures sought by low rank learning methods. We will exploit this flexibility in Section 4.4 below to achieve a novel extension of this framework to multi-view dictionary learning (White et al., 2012).

### 4.2.2 Application to Matrix Completion

We end this section with a brief illustration of how the framework can be applied to matrix completion. Recall that in this case one is given a partially observed data matrix $X$, where the goal is to find a low rank approximator $W$ that minimizes the reconstruction error on observed entries. In particular, we consider the following standard specialization of (53) to this problem:

$$\min_{W \in \mathbb{R}^{n \times m}} \ \tfrac{1}{2} \left\| \mathcal{P}(X - W) \right\|_{\mathrm{F}}^2 + \lambda \left\| W \right\|_{\mathrm{tr}}, \tag{61}$$

where, as noted above, $\mathcal{P} : \mathbb{R}^{n \times m} \to \mathbb{R}^{n \times m}$ is the mask operator that fills entries where $X$ is unobserved with zero. Surprisingly, Candès and Recht (2009) have proved that, under the low rank assumption, the solution of the convex relaxation (61) will recover the true matrix $X$ with high probability, even when only a small *random* portion of $X$ is observed.

Note that, since (61) is a specialization of (56), the modified GCG algorithm, Algorithm 2, can be immediately applied. Interestingly, such an approach contrasts with the most popular algorithm favored in the current literature: PG (or its accelerated variant APG; cf. Section 2). The two strategies, GCG versus PG, exhibit an interesting trade-off. Each iteration of PG (or APG) involves solving the PU (4), which in this case requires a *full* singular value decomposition (SVD) of the current approximant $W$. The modified GCG method, Algorithm 2, by contrast, only requires the polar of the gradient matrix $\nabla \ell(W)$ to be (approximately) computed in each iteration (i.e. dual of the trace norm; namely, the spectral norm), which is an order of magnitude cheaper to evaluate than a full SVD (cubic versus squared worst case time).[14] On the other hand, GCG has a slower theoretical rate than APG, $O(1/t)$ versus $O(1/t^2)$, in terms of the number of iterations required to achieve a given approximation error. Once the enhancements in the next section have been adopted, our experiments in Section 6.1 will demonstrate that GCG is far more efficient than APG for solving the matrix completion problem (61) to tolerances used in practice.

## 4.3 Fixed rank Local Optimization

Although the modified GCG algorithm, Algorithm 2, can be readily applied to the reformulated low rank learning problem (56), the sublinear rate of convergence established in

---

14. The per-iteration cost of computing the spectral norm can be further reduced for the sparse matrices occuring naturally in matrix completion, since the gradient matrix in (61) contains at most $k$ nonzero elements for $k$ equal to the number of observed entries in $X$. In such cases, an $\epsilon$ accurate estimate of the squared spectral norm can be obtained in $O\left(\frac{k \log(n+m)}{\sqrt{\epsilon}}\right)$ time (Jaggi, 2013; Kuczyński and Woźniakowski, 1992), offering a significant speed up over the general case. Although some speed ups might also be possible when computing the SVD for sparse matrices, the gap between the computational cost for evaluating the spectral versus trace norm for sparse matrices remains significant.

Corollary 13 is still too slow in practice. Nevertheless, an effective local improvement scheme can be developed that satisfies the properties of a `Relaxed` subroutine in Algorithm 2, which allows empirical convergence to be greatly accelerated. Recall that in the `Null` version of Algorithm 2 (cf. (47)), $\mathbf{w}_{t+1}$ is determined by a linear combination of the previous iterate $\mathbf{w}_t$ and the newly added atom $\mathbf{a}_t$. We first demonstrate that we can further improve $\mathbf{w}_{t+1}$ by solving a reformulated local problem, then show how this can still preserve the convergence guarantees of Algorithm 2.

The key reformulation relies on the following fact.

**Proposition 14.** *The gauge $\kappa$ induced by the set $\mathcal{A}$ in (57) can be re-expressed as*

$$\kappa(W) = \min\left\{\sum_{i=1}^{t} \|U_{:i}\|_c \ \|V_{i:}\|_r : UV = W\right\} = \min\left\{\frac{1}{2}\sum_{i=1}^{t}\left(\|U_{:i}\|_c^2 + \|V_{i:}\|_r^2\right) : UV = W\right\},$$

*for $U \in \mathbb{R}^{n \times t}$ and $V \in \mathbb{R}^{t \times m}$, as long as $t \geq mn$.*

The first equivalence in Proposition 14 is well-known since Grothendieck's work on tensor products of Banach spaces (where it is usually called the projective tensor norm). The second equality follows directly from the arithmetic-geometric mean inequality. In Appendix B.1 we provide a direct proof that allows us to bound the number of terms involved in the summation, which is made possible through the finite dimensional nature of the matrices. We note that the above lower bound $mn$ on $t$ may not be sharp. For example, if $\|\cdot\|_r = \|\cdot\|_c = \|\cdot\|_2$ then $\kappa$ is the trace norm, where $\sum_{i=1}^{t}(\|U_{:i}\|_c^2 + \|V_{i:}\|_r^2)$ is simply $\|U\|_F^2 + \|V\|_F^2$, the sum of the squared Frobenius norms; that is, Proposition 14 yields the well-known variational form of the trace norm in this case (Srebro et al., 2005). Observe that $t$ in this case can be as low as the rank of $W$. Based on Proposition 14 the objective in (56) can then be approximated at iteration $t$ as

$$\mathcal{F}_t(U, V) := \ell(UV) + \frac{\lambda}{2}\sum_{i=1}^{t}\left(\|U_{:i}\|_c^2 + \|V_{i:}\|_r^2\right), \tag{62}$$

which provides a particularly convenient objective for local improvement since the problem is now both smooth and unconstrained, allowing efficient unconstrained minimization routines to be applied. Note however that (62) is not *jointly* convex in $U$ and $V$ since the size $t$ (i.e. the iteration counter) is fixed. When $t$ is sufficiently large, any local minimizer[15] of (62) *globally* solves the original problem (56) (Burer and Monteiro, 2005), but we locally minimize (62) in each iteration, even for small $t$.

The main issue in adding a local optimization step in Algorithm 2 is that the `Relaxed` condition (43) in Theorem 12 must be maintained. Fortunately, the locally improved iterate $W_{t+1}$ and corresponding bound $\rho_{t+1}$ needed to achieve the desired condition can easily be recovered from any local minimizer $(U^*, V^*)$ of (62) via:

$$W_{t+1} = U^*V^* \quad \text{and} \quad \rho_{t+1} = \frac{1}{2}\sum_{i=1}^{t+1}\left(\|U_{:i}^*\|_c^2 + \|V_{i:}^*\|_r^2\right). \tag{63}$$

---

15. The catch here is that seeking a local minimizer can still be highly non-trivial, sometimes even hard to certify the local optimality.

---

**Algorithm 3** GCG variant for low rank learning.

---

**Require:** The atomic set $\mathcal{A}$.

1: Initialize $W_0 = \mathbf{0}, \quad \rho_0 = 0, \quad U_0 = V_0 = \Lambda_0 = \emptyset.$

2: **for** $t = 0, 1, \ldots$ **do**

3: $\quad (\mathbf{u}_t, \mathbf{v}_t) \leftarrow \arg \min\limits_{\mathbf{u}\mathbf{v}^\top \in \mathcal{A}} \left\langle \nabla \ell(W_t), \mathbf{u}\mathbf{v}^\top \right\rangle$

4: $\quad (\eta_t, \theta_t) \leftarrow \arg \min\limits_{0 \leq \eta \leq 1, \theta \geq 0} \ell((1-\eta)W_t + \theta\,\mathbf{u}_t\mathbf{v}_t^\top) + \lambda((1-\eta)\rho_t + \theta)$

5: $\quad U_{\mathrm{init}} \leftarrow (\sqrt{1-\eta_t}U_t, \sqrt{\theta_t}\mathbf{u}_t), \; V_{\mathrm{init}} \leftarrow (\sqrt{1-\eta_t}V_t^\top, \sqrt{\theta_t}\mathbf{v}_t)^\top$

6: $\quad (U_{t+1}, V_{t+1}) = \texttt{Improve}(\mathcal{F}_{t+1}, U_{\mathrm{init}}, V_{\mathrm{init}}),$
$\quad$ i.e. find a local minimum of $\mathcal{F}_{t+1}$ with $U$ and $V$ initialized to $U_{\mathrm{init}}$ and $V_{\mathrm{init}}$ resp.

7: $\quad W_{t+1} \leftarrow U_{t+1}V_{t+1}$

8: $\quad \rho_{t+1} \leftarrow \frac{1}{2}\sum_{i=1}^{t+1}(\|(U_{t+1})_{:i}\|_c^2 + \|(V_{t+1})_{i:}\|_r^2)$

9: **end for**

---

This improves the quality of the current iterate without increasing the size of the dictionary.

Once this modification is added to Algorithm 2 we obtain the final variant, Algorithm 3. In this algorithm, Lines 3 and 4 remain equivalent to before. Line 5 splits the iterate into two parts $U_{\mathrm{init}}$ and $V_{\mathrm{init}}$, while line 6 finds a local minimizer of the reformulated local objective (62) with the designated initialization. The last two steps restore the GCG iterate. To show that Algorithm 3 still has the same $O(1/t)$ rate of convergence established in Corollary 13, it suffices to prove that the introduced $\texttt{Improve}$ subroutine is $\texttt{Relaxed}$ in the sense required by (43). Indeed, by construction

$$
\begin{aligned}
\ell(W_{t+1}) + \lambda\rho_{t+1} &= \ell(U_{t+1}V_{t+1}) + \frac{\lambda}{2}\sum_{i=1}^{t+1}\left(\|(U_{t+1})_{:i}\|_c^2 + \|(V_{t+1})_{i:}\|_r^2\right)\\
&= \mathcal{F}_{t+1}(U_{t+1}, V_{t+1})\\
&\leq \mathcal{F}_{t+1}(U_{\mathrm{init}}, V_{\mathrm{init}})\\
&\leq \ell((1-\eta_t)U_tV_t + \theta_t\mathbf{u}_t\mathbf{v}_t^\top) + \lambda\theta_t + \lambda\frac{(1-\eta_t)}{2}\sum_{i=1}^{t}\left(\|(U_t)_{:i}\|_c^2 + \|(V_t)_{i:}\|_r^2\right)\\
&= \ell((1-\eta_t)W_t + \theta_t\mathbf{u}_t\mathbf{v}_t^\top) + \lambda(1-\eta_t)\rho_t + \lambda\theta_t,
\end{aligned}
\tag{64}
$$

where in addition $\kappa(W_{t+1}) \leq \rho_{t+1}$ follows from (63) and Proposition 14. Although the $\texttt{Improve}$ subroutine is not guaranteed to achieve *strict* improvement in every iteration, we observe in our experimental analysis in Section 6 below that Algorithm 3 is almost always significantly faster than the $\texttt{Null}$ version of Algorithm 2. In other words, local acceleration appears to make a big difference in practice.

**Remark 11.** *Interleaving local improvement with a globally convergent procedure has been considered by (Mishra et al., 2013) and (Laue, 2012). In particular, Laue (2012) considered the constrained problem* (26), *which unfortunately leads to considering a constrained minimization for local improvement that can be less efficient than the unconstrained approach considered in* (62). *In addition, (Laue, 2012) used the original objective $F$ directly in local optimization. This leads to significant challenge when the regularizer $f$ is complicated such as trace norm, and we avoid this issue by using a smooth surrogate objective.*

*Focusing on trace norm regularization only, Mishra et al. (2013) propose a trust-region procedure that locally minimizes the original objective on the Stiefel manifold and the positive semidefinite cone, which unfortunately also requires constrained minimization and dynamic maintenance of the singular value decomposition of a small matrix. The rate of convergence of this procedure has not been analyzed.*

**Remark 12.** *Finally, it is useful to note that local minimization of* (62) *constitutes the primary bottleneck of Algorithm 3, consuming over 50% of the computation. Since* (63) *is the same objective considered by the standard fixed rank approach to collaborative filtering, the overall performance of Algorithm 3 can be substantially improved by leveraging the techniques that have been developed for this specific objective. These improvements include system level optimizations such as stochastic solvers on distributed, parallel and asynchronous infrastructures (Zhuang et al., 2013; Yun et al., 2014). In the experiments conducted in Section 6 below, we merely solve* (63) *using a generic L-BFGS solver to highlight the efficiency of the GCG framework itself. Nevertheless, such a study leaves significant room for further straightforward acceleration.*

### 4.4 Multi-view Dictionary Learning

Finally, in this section we show how low rank learning can be generalized to a multi-view setting where, once again, the modified GCG algorithm can efficiently provide optimal solutions. The key challenge in this case is to develop column and row norms that can capture the more complex problem structure.

Consider a two-view learning task where one is given $m$ paired observations $\left\{ \begin{bmatrix} \mathbf{x}_j \\ \mathbf{y}_j \end{bmatrix} \right\}$ consisting of an $x$-view and a $y$-view with lengths $n_1$ and $n_2$ respectively. The goal is to infer a set of latent representations, $\mathbf{h}_j$ (of dimension $t < \min(n_1, n_2)$), and reconstruction models parameterized by the matrices $A$ and $B$, such that $A\mathbf{h}_j \approx \mathbf{x}_j$ and $B\mathbf{h}_j \approx \mathbf{y}_j$ for all $j$. In particular, let $X$ denote the $n_1 \times m$ matrix of $x$-view observations, $Y$ the $n_2 \times m$ matrix of $y$-view observations, and $Z = \begin{bmatrix} X \\ Y \end{bmatrix}$ the concatenated $(n_1 + n_2) \times m$ matrix. The problem can then be expressed as recovering an $(n_1 + n_2) \times t$ matrix of model parameters, $C = \begin{bmatrix} A \\ B \end{bmatrix}$, and a $t \times m$ matrix of latent representations, $H$, such that $Z \approx CH$.

The key assumption behind multi-view learning is that each of the two views, $\mathbf{x}_j$ and $\mathbf{y}_j$, is conditionally independent given the shared latent representation, $\mathbf{h}_j$. Although multi-view data can always be concatenated hence treated as a single view, explicitly representing multiple views enables more accurate identification of the latent representation (as we will see), when the conditional independence assumption holds. To better motivate this model, it is enlightening to first reinterpret the classical formulation of multi-view dictionary learning, which is given by *canonical correlation analysis* (CCA). Typically, it is expressed as the problem of projecting two views so that the correlation between them is maximized (Hardoon et al., 2004). Assuming the data is centered (i.e. $X\mathbf{1} = \mathbf{0}$ and $Y\mathbf{1} = \mathbf{0}$), the sample covariance of $X$ and $Y$ is given by $XX^\top/m$ and $YY^\top/m$ respectively. CCA can then be expressed as an optimization over matrix variables

$$\max_{U,V} \ \mathrm{tr}(U^\top XY^\top V) \qquad \text{s.t.} \qquad U^\top XX^\top U = V^\top YY^\top V = I \qquad (65)$$

for $U \in \mathbb{R}^{n_1 \times t}$, $V \in \mathbb{R}^{n_2 \times t}$ (De Bie et al., 2005). Although this classical formulation (65) does not make the shared latent representation explicit, CCA can be expressed by a generative model: given a latent representation, $\mathbf{h}_j$, the observations $\mathbf{x}_j = A\mathbf{h}_j + \boldsymbol{\epsilon}_j$ and $\mathbf{y}_j = B\mathbf{h}_j + \boldsymbol{\nu}_j$ are generated by a linear mapping plus independent zero mean Gaussian noise, $\boldsymbol{\epsilon} \sim \mathcal{N}(\mathbf{0}, \Sigma_x)$, $\boldsymbol{\nu} \sim \mathcal{N}(\mathbf{0}, \Sigma_y)$ (Bach and Jordan, 2006). Interestingly, one can show that the classical CCA problem (65) is equivalent to the following multi-view dictionary learning problem.

**Proposition 15** ((White et al., 2012, Proposition 1)). *Fix $t$, let*

$$
\begin{aligned}
(A, B, H) &= \arg \min_{A: \|A_{:,i}\|_2 \leq 1} \min_{B: \|B_{:,i}\|_2 \leq 1} \min_H \left\| \tilde{Z} - \begin{bmatrix} A \\ B \end{bmatrix} H \right\|_{\mathrm{F}}^2, \\
where \quad \tilde{Z} &= \begin{bmatrix} (XX^\top)^{-1/2} X \\ (YY^\top)^{-1/2} Y \end{bmatrix}.
\end{aligned}
\tag{66}
$$

*Then $U = (XX^\top)^{-\frac{1}{2}} A$ and $V = (YY^\top)^{-\frac{1}{2}} B$ provide an optimal solution to (65).*

Proposition 15 demonstrates how formulation (66) respects the conditional independence of the separate views: given a latent representation $\mathbf{h}_j$, the reconstruction losses on the two views, $\mathbf{x}_j$ and $\mathbf{y}_j$, cannot influence each other, since the reconstruction models $A$ and $B$ are *individually* constrained. By contrast, in single-view dictionary learning (i.e. *principal components analysis*) $A$ and $B$ are concatenated in the larger variable $C = \begin{bmatrix} A \\ B \end{bmatrix}$, where $C$ as a whole is constrained but $A$ and $B$ are not. $A$ and $B$ must then compete against each other to acquire magnitude to explain their respective "views" given $\mathbf{h}_j$ (i.e. conditional independence is not enforced). Such sharing can be detrimental if the two views really are conditionally independent given $\mathbf{h}_j$.

This matrix factorization viewpoint has recently allowed dictionary learning to be extended from the single-view setting to the multi-view setting (Quadrianto and Lampert, 2011; Jia et al., 2010). However, despite its elegance, the computational tractability of CCA hinges on its restriction to squared loss under a particular normalization, precluding other important losses. To combat these limitations, we apply the same convex relaxation principle as in Section 4.2. In particular, we reformulate (66) by first incorporating an arbitrary loss function $\ell$ that is convex in its first argument, and then relaxing the rank constraint by replacing it with a rank-reducing regularizer on $H$. This amounts to the following training problem that extends (Jia et al., 2010):

$$
\begin{aligned}
\min_{A, B, H} \ell \left( \begin{bmatrix} A \\ B \end{bmatrix} H; Z \right) + \lambda \sum_i \|H_{i:}\|_2, \quad \text{s.t.} \quad \begin{bmatrix} A_{:i} \\ B_{:i} \end{bmatrix} \in \mathcal{C} \ \text{ for all } i, \\
where \quad \mathcal{C} := \left\{ \begin{bmatrix} \mathbf{a} \\ \mathbf{b} \end{bmatrix} : \|\mathbf{a}\|_2 \leq \beta, \ \|\mathbf{b}\|_2 \leq \gamma \right\}.
\end{aligned}
\tag{67}
$$

The $l_2$ norm on the *rows* of $H$ is a regularizer that encourages the rows to become sparse, hence reducing the size of the learned representation (Argyriou et al., 2008). Since the size of the latent representation (the number of rows of $H$) is not fixed, our results in Section 4.2 immediately allows the problem to be solved globally and efficiently for $A$, $B$ and $H$. This considerably improves the previous local solutions in (Jia et al., 2010; Quadrianto and Lampert, 2011) which fixed the latent dimension and approached the nonconvex problem by alternating block coordinate descent between $A$, $B$ and $H$.

Indeed, (67) fits exactly in the framework of (53), with $\|\cdot\|_r$ being the $l_2$ norm and

$$\|\mathbf{c}\|_c = \max\left(\frac{1}{\beta}\|\mathbf{a}\|_2, \frac{1}{\gamma}\|\mathbf{b}\|_2\right) \quad \text{for} \quad \mathbf{c} = \begin{bmatrix} \mathbf{a} \\ \mathbf{b} \end{bmatrix}. \tag{68}$$

Next we will show how to efficiently compute the polar operator (59) that underpins Algorithm 3 in this multi-view model.

### 4.4.1 EFFICIENT POLAR OPERATOR

The polar operator (59) for problem (67) turns out to admit an efficient but nontrivial solution. For simplicity, we assume $\beta = \gamma = 1$ while the more general case can be dealt with in exactly the same way. The proof of the following key proposition is in Appendix B.2.

**Proposition 16.** *Let $I_1 := \text{diag}\left(\begin{bmatrix}\mathbf{1}\\\mathbf{0}\end{bmatrix}\right)$ and $I_2 := \text{diag}\left(\begin{bmatrix}\mathbf{0}\\\mathbf{1}\end{bmatrix}\right)$ be the identity matrices on the subspaces spanned by x-view and y-view, respectively. For any $\mu > 0$, denote $D_\mu = \sqrt{1+\mu}I_1 + \sqrt{1+1/\mu}I_2$ as a diagonal scaling of the two identity matrices. Then, for the row norm $\|\cdot\|_r = \|\cdot\|_2$ and column norm $\|\cdot\|_c$ defined as in (68), the corresponding polar (59) is given by*

$$\kappa^\circ(G) = \min_{\mu \geq 0} \|D_\mu G\|_{\text{sp}}. \tag{69}$$

Despite its variational form, the polar (69) can still be efficiently computed by considering its squared form and expanding:

$$\|D_\mu G\|_{\text{sp}}^2 = \left\|G^\top D_\mu D_\mu G\right\|_{\text{sp}} = \left\|G^\top \left[(1+\mu)I_1 + (1+1/\mu)I_2\right] G\right\|_{\text{sp}}.$$

Clearly, the term inside the spectral norm $\|\cdot\|_{\text{sp}}$ is convex in $\mu$ on $\mathbb{R}_+$ while the spectral norm is convex and increasing on the cone of positive semidefinite matrices, hence by the usual rules of convex composition the left hand side above is a *convex* function of $\mu$. Consequently we can use a subgradient algorithm, or even golden section search, to find the optimal $\mu$ in (69). Given an optimal $\mu$, an optimal pair of "atoms" $\mathbf{a}$ and $\mathbf{b}$ can be recovered via the KKT conditions; see Appendix B.3 for details. Thus, by using (69) Algorithm 3 is able to solve the resulting optimization problem far more efficiently than the initial procedure offered in (White et al., 2012).

### 4.4.2 A MODIFIED POWER ITERATION

In practice, one can compute the polar (69) much more efficiently by deploying a heuristic variant of the power iteration method, whose motivation arises from the KKT conditions in the proof of Proposition 16 (see (106) in Appendix B.3). In particular, two normalized vectors $\mathbf{a}$ and $\mathbf{b}$ are optimal for the polar (59) if and only if there exist nonnegative scalars $\mu_1$ and $\mu_2$ such that

$$\begin{bmatrix} \mu_1 \mathbf{a} \\ \mu_2 \mathbf{b} \end{bmatrix} = GG^\top \begin{bmatrix} \mathbf{a} \\ \mathbf{b} \end{bmatrix}, \tag{70}$$

$$\mu_1 I_1 + \mu_2 I_2 \succeq GG^\top. \tag{71}$$

Therefore, after randomly initializing $\mathbf{a}$ and $\mathbf{b}$ with unit norm, it is natural to extend the power iteration method by repeating the update:

$$\begin{bmatrix} \mathbf{s} \\ \mathbf{t} \end{bmatrix} \leftarrow GG^\top \begin{bmatrix} \mathbf{a} \\ \mathbf{b} \end{bmatrix}, \quad \text{then} \quad \begin{bmatrix} \mathbf{a}^{\text{new}} \\ \mathbf{b}^{\text{new}} \end{bmatrix} \leftarrow \begin{bmatrix} \frac{1}{\|\mathbf{s}\|_2}\mathbf{s} \\ \frac{1}{\|\mathbf{t}\|_2}\mathbf{t} \end{bmatrix}.$$

This procedure usually converges to a solution where the Lagrange multipliers $\mu_1$ and $\mu_2$, recovered from (70), satisfy (71) confirming their optimality. When the occasional failure occurs we revert to (69) and follow the recovery rule detailed in Appendix B.3.

## 5. Application to Structured Sparse Estimation

In addition to dictionary learning, the GCG approach can also be very effective in optimizing other sparse models. In this section we illustrate a second key application: *structured sparse estimation* (Bach et al., 2012a). Once again, the major challenge is developing an efficient algorithm for the polar operator required by GCG. Motivated by the work of (Obozinski and Bach, 2012) we first consider in Section 5.1 a family of structured sparse regularizers induced by a cost function on subsets of variables. Section 5.2 then introduces a novel "lifting" construction that enables us to express these regularizers as linear functions, which, after some reformulation, allows efficient evaluation via a simple linear program (LP). We then illustrate the particularly important example of overlapping group lasso (Jenatton et al., 2011) in Section 5.3. By exploiting the structure of the overlapping group lasso we show that the corresponding LP can be reduced to a piecewise linear objective over a simple domain, which allows a further significant reduction in computation time via smoothing (Nesterov, 2005).

### 5.1 Convex Regularizers for Structured Sparse Estimation

Recently, Obozinski and Bach (2012) have established a principled method for deriving regularizers from a subset cost function $J : 2^{[n]} \to \overline{\mathbb{R}}_+$, where $2^{[n]}$ denotes the power set of $[n] := \{1, \ldots, n\}$ and $\overline{\mathbb{R}}_+ := \mathbb{R}_+ \cup \{\infty\}$. In particular, given such a cost function they define the gauge regularizer

$$\kappa_J(\mathbf{w}) = \inf\{\rho \geq 0 : \mathbf{w} \in \rho \, \mathrm{conv}(S_J)\}, \tag{72}$$

$$\text{where} \quad S_J = \left\{ \mathbf{w}_A : \|\mathbf{w}_A\|_{\mathsf{p}}^{\mathsf{q}} = \frac{1}{J(A)}, \ \emptyset \neq A \subseteq [n] \right\}.$$

Note that each set $A \subseteq [n]$ denotes a subset of the coordinates for weight vectors $\mathbf{w}$. Here $\rho$ is a scalar, $\mathsf{p}$ and $\mathsf{q}$ satisfy $\mathsf{p}, \mathsf{q} \geq 1$ with $\frac{1}{\mathsf{p}} + \frac{1}{\mathsf{q}} = 1$, and $\mathbf{w}_A$ denotes a duplicate of $\mathbf{w}$ with all coordinates not in $A$ set to 0. Furthermore, we have tacitly assumed that $J(A) = 0$ iff $A = \emptyset$ in (72). The gauge $\kappa_J$ defined in (72) is also the *atomic norm* with the set of atoms $S_J$ (Chandrasekaran et al., 2012).

The *polar* of $\kappa_J$ is simply the support function of $S_J$:

$$\kappa_J^\circ(\mathbf{g}) := \sup_{\mathbf{w}:\kappa_J(\mathbf{w}) \leq 1} \langle \mathbf{g}, \mathbf{w} \rangle = \max_{\mathbf{w} \in S_J} \langle \mathbf{g}, \mathbf{w} \rangle = \max_{\emptyset \neq A \subseteq [n]} \frac{1}{J(A)^{1/\mathsf{q}}} \|\mathbf{g}_A\|_{\mathsf{q}}. \tag{73}$$

Here, the second equality uses the definition of support function, and the third follows from (72). By varying $\mathsf{p}$ and $J$, one can generate a class of sparsity inducing regularizers that includes most current proposals (Obozinski and Bach, 2012). For example, if $J(A) = |A|$ (the cardinality of $A$), then $\kappa_J^\circ$ is the $l_\infty$ norm and $\kappa_J$ is the usual $l_1$ norm. More importantly, one can encode structural information through the cost function $J$, which selects and establishes preferences over the set of atoms $S_J$.

Using a standard duality argument we can deduce the gauge (72) from its polar (73):

$$\kappa(\mathbf{w}) = (\kappa^\circ)^\circ(\mathbf{w}) = \inf \left\{ \sum_{\emptyset \neq A \subseteq [n]} J(A)^{1/\mathsf{q}} \left\| \mathbf{u}^A \right\|_{\mathsf{p}} : \mathbf{w} = \sum_{\emptyset \neq A \subseteq [n]} \mathbf{u}^A, \ (\mathbf{u}^A)_i = 0, \ \forall i \notin A \right\}. \quad (74)$$

(A more general duality result along this line can be found in (Zhang et al., 2013, Proposition 6)). From (74) and (73) it is clear that a naive way to compute the gauge or its polar requires enumerating all subsets of the ground set $[n]$, too expensive to afford even when $n$ is moderate. This is where useful structural information on the cost function $J$ comes into play from a computational perspective, as we will see shortly.

As discussed in Section 2, one of the standard approaches for solving the composite minimization problem

$$\min_{\mathbf{w} \in \mathbb{R}^n} \ell(\mathbf{w}) + \lambda \, \kappa_J(\mathbf{w}) \quad (75)$$

is to use the *accelerated proximal gradient* (APG) algorithm, which is guaranteed to find an $\epsilon$ accurate solution in $O(1/\sqrt{\epsilon})$ iterations (Beck and Teboulle, 2009; Nesterov, 2013). Unfortunately, the proximal update (PU) (4) in each iteration can be quite difficult to compute when the gauge $\kappa_J$ encodes combinatorial structure.

Here we instead consider solving (75) by the *generalized conditional gradient* (GCG) algorithm developed above in Section 3. Unlike APG, GCG only requires the *polar operator* of $\kappa_J$ to be computed in each iteration, given by the argument of (73):

$$\mathbb{P}_J(\mathbf{g}) = \arg \max_{\mathbf{w} \in S_J} \langle \mathbf{g}, \mathbf{w} \rangle = J(C)^{\frac{-1}{\mathsf{q}}} \arg \max_{\mathbf{w}:\|\mathbf{w}\|_{\mathsf{p}}=1} \langle \mathbf{g}_C, \mathbf{w} \rangle \quad (76)$$

$$\text{where} \quad C = C(\mathbf{g}) = \arg \max_{\emptyset \neq A \subseteq [n]} \frac{1}{J(A)} \left\| \mathbf{g}_A \right\|_{\mathsf{q}}^{\mathsf{q}}. \quad (77)$$

Algorithm 4 outlines a GCG procedure that has been further specialized to solving (75). In detail, Line 4 evaluates the polar operator, which provides a descent direction $\mathbf{v}_t$; Line 5 finds the optimal step sizes for combining the current iterate $\mathbf{w}_t$ with the direction $\mathbf{v}_t$; and Line 6 locally improves the objective (75) with the gauge $\kappa$ being replaced by the equivalent equation (74). Here, to circumvent the combinatorial complexity in evaluating the gauge $\kappa$, we simply maintain the same support patterns $\{A_i\}_{i=0}^t$ in previous rounds but re-optimize the parameters $\{\mathbf{u}^i\}_{i=0}^t$. This is essentially the totally corrective variant we mentioned in Section 3.5, cf. (51). By Corollary 13, GCG can find an $\epsilon$ accurate solution to (75) in $O(1/\epsilon)$ steps, provided only that the polar (76)-(77) is computed to $\epsilon$ accuracy.

Although GCG has a slower theoretical convergence rate than APG, the local optimization improvement (Line 5) and the totally corrective step (Line 6) often yield faster

---

**Algorithm 4** Generalized conditional gradient (GCG) for optimizing (75).

---

1: Initialize $\mathbf{w}_0 \leftarrow \mathbf{0}, \quad \rho_0 \leftarrow 0$.
2: **for** $t = 0, 1, \dots$ **do**
3: $\quad \mathbf{g}_t \leftarrow -\nabla \ell(\mathbf{w}_t)$.
4: $\quad$ Polar operator: $\mathbf{v}_t \leftarrow \mathbb{P}_J(\mathbf{g}_t), A_t \leftarrow C(\mathbf{g}_t)$, and $C(\cdot)$ is defined in (77).
5: $\quad$ 2-D Conic search: $(\alpha, \beta) := \arg\min_{\alpha \geq 0, \beta \geq 0} \ell(\alpha \mathbf{w}_t + \beta \mathbf{v}_t) + \lambda(\alpha \rho_t + \beta)$.
6: $\quad$ Local re-optimization: $\{\mathbf{u}^i\}_0^t := \arg\min_{\{\mathbf{u}^i = \mathbf{u}^i_{A_i}\}} \ell(\sum_i \mathbf{u}^i) + \lambda \sum_i J(A_i)^{\frac{1}{q}} \|\mathbf{u}^i\|_{\mathsf{p}}$
$\quad\quad\quad\quad$ where the $\{\mathbf{u}^i\}$ are initialized by $\mathbf{u}^i = \alpha \mathbf{q}_i$ for $i < t$ and $\mathbf{u}^i = \beta \mathbf{v}_i$ for $i = t$.
7: $\quad \mathbf{w}_{t+1} \leftarrow \sum_i \mathbf{u}^i, \quad \mathbf{q}_i \leftarrow \mathbf{u}^i$ for $i \leq t, \quad \rho_{t+1} \leftarrow \sum_i J(A_i)^{\frac{1}{q}} \|\mathbf{u}^i\|_{\mathsf{p}}$.
8: **end for**

---

convergence in practice. Our main goal in this section therefore is to extend GCG to structured sparse models by developing efficient algorithms that compute the polar operator in (76)-(77). As pointed out by Obozinski and Bach, 2012), the polar (73) can be evaluated by a secant method which is tractable when $J$ is submodular (Bach, 2010, §8.4). However, we show that it is possible to do significantly better by bypassing submodular function minimization entirely.

## 5.2 Polar Operators for Atomic Norms

Our first main contribution in this section is to develop a general class of atomic norm regularizers whose polar operator (76)-(77) can be computed efficiently. To begin, consider the case of a (partially) *linear* function $J$ where there exists a $\mathbf{c} \in \mathbb{R}^n$ such that $J(A) = \langle \mathbf{c}, \mathbf{1}_A \rangle$ for all $A \in \text{dom } J$ (note that the domain need not be a lattice). A few useful regularizers can be generated by such functions: for example, the $\mathsf{I}_1$ norm can be derived from $J(A) = \langle \mathbf{1}, \mathbf{1}_A \rangle$, which is linear. Unfortunately, linearity is too restrictive to capture most structured regularizers of interest, therefore we will need to expand the space of functions $J$ we consider. To do so, we introduce the more general class of *marginalized linear functions*: we say that $J$ is marginalized linear if there exists a nonnegative linear function $M$ on an extended domain $2^{[n+l]}$ such that its marginalization to $2^{[n]}$ is exactly $J$:

$$ J(A) = \min_{B: A \subseteq B \subseteq [n+l]} M(B), \quad \forall A \subseteq [n]. \tag{78} $$

Essentially, such a function $J$ is "lifted" to a larger domain where it becomes linear. The key question is whether the polar $\kappa_J^\circ$ can be efficiently evaluated for such functions.

To develop an efficient procedure for computing the polar $\kappa_J^\circ$ for such functions, first consider the simpler case of computing the polar $\kappa_M^\circ$ for a nonnegative linear function $M$. Note that, by linearity, the function $M$ can be expressed as $M(B) = \langle \mathbf{b}, \mathbf{1}_B \rangle$ for some $\mathbf{b} \in \mathbb{R}_+^{n+l}$ such that $B \in \text{dom } M \subseteq 2^{[n+l]}$. Since the effective domain of $M$ need not be the entire space $2^{[n+l]}$ in general, we make use of the specialized polytope:

$$ P := \text{conv}\{\mathbf{1}_B : B \in \text{dom } M\} \subseteq [0, 1]^{n+l}. \tag{79} $$

Note that $P$ may have exponentially many faces. From the definition (73) one can then re-express the polar $\kappa_M^\circ$ as:

$$\kappa_M^\circ(\mathbf{g}) = \max_{\emptyset \neq B \in \text{dom } M} \frac{\|\mathbf{g}_B\|_{\mathsf{q}}}{M(B)^{1/\mathsf{q}}} = \left( \max_{\mathbf{0} \neq \mathbf{w} \in P} \frac{\langle \tilde{\mathbf{g}}, \mathbf{w} \rangle}{\langle \mathbf{b}, \mathbf{w} \rangle} \right)^{1/\mathsf{q}} \quad \text{where } \tilde{g}_i = |g_i|^{\mathsf{q}} \ \forall i, \qquad (80)$$

where we have used the fact that the linear-fractional objective must attain its maximum at vertices of $P$; that is, at $\mathbf{1}_B$ for some $B \in \text{dom } M$. Although the linear-fractional program (80) can be reduced to a sequence of LPs (Dinkelbach, 1967), a single LP suffices for our purposes, as we now demonstrate. Indeed, let us first remove the constraint $\mathbf{w} \neq \mathbf{0}$ by considering the alternative polytope:

$$Q := P \cap \{\mathbf{w} \in \mathbb{R}^{n+l} : \langle \mathbf{1}, \mathbf{w} \rangle \geq 1\}. \qquad (81)$$

As shown in Appendix C.1, all vertices of $Q$ are scalar multiples of the nonzero vertices of $P$. Since the objective in (80) is scale invariant, one can restrict the constraints to $\mathbf{w} \in Q$. Then, by applying transformations $\tilde{\mathbf{w}} = \mathbf{w}/\langle \mathbf{b}, \mathbf{w} \rangle$, $\sigma = 1/\langle \mathbf{b}, \mathbf{w} \rangle$, problem (80) can be equivalently re-expressed by the LP:

$$\max_{\tilde{\mathbf{w}}, \sigma} \ \langle \tilde{\mathbf{g}}, \tilde{\mathbf{w}} \rangle, \ \text{subject to} \ \tilde{\mathbf{w}} \in \sigma Q, \ \langle \mathbf{b}, \tilde{\mathbf{w}} \rangle = 1, \ \sigma \geq 0. \qquad (82)$$

Of course, whether this LP can be solved efficiently depends on the structure of $Q$, which we consider below.

Finally, we note that the same formulation allows the polar to be efficiently computed for a *marginalized* linear function $J$ via a simple reduction: Consider any $\mathbf{g} \in \mathbb{R}^n$ and let $[\mathbf{g}; \mathbf{0}] \in \mathbb{R}^{n+l}$ denote $\mathbf{g}$ padded by $l$ zeros. Then $\kappa_J^\circ(\mathbf{g}) = \kappa_M^\circ([\mathbf{g}; \mathbf{0}])$ for all $\mathbf{g} \in \mathbb{R}^n$, since

$$(\kappa_J^\circ(\mathbf{g}))^{\mathsf{q}} = \max_{\emptyset \neq A \subseteq [n]} \frac{\|\mathbf{g}_A\|_{\mathsf{q}}^{\mathsf{q}}}{J(A)} = \max_{\emptyset \neq A \subseteq [n]} \frac{\|\mathbf{g}_A\|_{\mathsf{q}}^{\mathsf{q}}}{\min_{B:A \subseteq B \subseteq [n+l]} M(B)}$$

$$= \max_{\emptyset \neq A \subseteq B} \frac{\|\mathbf{g}_A\|_{\mathsf{q}}^{\mathsf{q}}}{M(B)} = \max_{B: \emptyset \neq B \subseteq [n+l]} \frac{\|[\mathbf{g}; \mathbf{0}]_B\|_{\mathsf{q}}^{\mathsf{q}}}{M(B)} = (\kappa_M^\circ([\mathbf{g}; \mathbf{0}]))^{\mathsf{q}}. \qquad (83)$$

To understand the middle equality in (83), observe that for a fixed $B$ the optimal $A$ is attained at $A = B \cap [n]$. If $B \cap [n]$ is empty, then $\|[\mathbf{g}; \mathbf{0}]_B\| = 0$ and the corresponding $B$ cannot be the maximizer of the last term, unless $\kappa_J^\circ(\mathbf{g}) = 0$, in which case it is easy to see $\kappa_M^\circ([\mathbf{g}; \mathbf{0}]) = 0$.

Although we have kept our development general so far, the idea is clear: Once an appropriate "lifting" has been found so that the polytope $Q$ in (81) can be compactly represented, the polar (76) can be reformulated as the LP (82), for which efficient implementations can be sought. We now demonstrate this new methodology for an important structured regularizer: overlapping group sparsity.

## 5.3 Group Sparsity

For a general formulation of group sparsity, let $\mathcal{G} \subseteq 2^{[n]}$ be a set of variable groups (subsets) that possibly overlap (Huang et al., 2011; Obozinski and Bach, 2012; Zhao et al., 2009).

Here we use $i \in [n]$ to index variables and $G \in \mathcal{G}$ to index groups. Consider the cost function over variable groups $J_{\mathsf{g}} : 2^{[n]} \to \mathbb{R}_+$ defined by:

$$J_{\mathsf{g}}(A) = \sum_{G \in \mathcal{G}} c_G \, \mathbb{1}(A \cap G \neq \emptyset), \tag{84}$$

where $c_G$ is a nonnegative cost and $\mathbb{1}$ is an indicator such that $\mathbb{1}(\cdot) = 1$ if its argument is true, and 0 otherwise. The value $J_{\mathsf{g}}(A)$ provides a weighted count of how many groups overlap with $A$. Unfortunately, $J_{\mathsf{g}}$ is not linear, so we need to re-express it to recover an efficient polar operator. To do so, augment the domain by adding $l = |\mathcal{G}|$ variables such that each new variable $G$ corresponds to a group $G$. Then define a weight vector $\mathbf{b} \in \mathbb{R}_+^{n+l}$ such that $b_i = 0$ for $i \leq n$ and $b_G = c_G$ for $n < G \leq n + l$. Finally, consider the linear cost function $M_{\mathsf{g}} : 2^{[n+l]} \to \overline{\mathbb{R}}_+$ defined by:

$$M_{\mathsf{g}}(B) = \begin{cases} \langle \mathbf{b}, \mathbf{1}_B \rangle & \text{if } i \in B \Rightarrow G \in B, \ \forall \, i \in G \in \mathcal{G}; \\ \infty & \text{otherwise.} \end{cases} \tag{85}$$

The constraint ensures that if a variable $i \leq n$ appears in the set $B$, then every variable $G$ corresponding to a group that contains $i$ must also appear in $B$. By construction, $M_{\mathsf{g}}$ is a nonnegative linear function. It is also easy to verify that $J_{\mathsf{g}}$ satisfies (78) with respect to $M_{\mathsf{g}}$.

To compute the corresponding polar, observe that the effective domain of $M_{\mathsf{g}}$ is a lattice, hence (73) can be solved by combinatorial methods. However, we can do better by exploiting problem structure in the LP. For example, observe that the polytope (79) can now be compactly represented as:

$$P_{\mathsf{g}} = \{ \mathbf{w} \in \mathbb{R}^{n+l} : \mathbf{0} \leq \mathbf{w} \leq \mathbf{1}, w_i \leq w_G, \forall \, i \in G \in \mathcal{G} \}. \tag{86}$$

Indeed, it is easy to verify that the integral vectors in $P_{\mathsf{g}}$ are precisely $\{ \mathbf{1}_B : B \in \mathrm{dom}\, M_{\mathsf{g}} \}$. Moreover, the linear constraint in (86) is totally unimodular (TUM) since it is the incidence matrix of a bipartite graph (variables and groups), hence $P_{\mathsf{g}}$ is the convex hull of its integral vectors (Schrijver, 1986). Using the fact that the scalar $\sigma$ in (82) admits a closed form solution $\sigma = \langle \mathbf{1}, \tilde{\mathbf{w}} \rangle$ in this case, the LP (82) can be reduced to:

$$\max_{\tilde{\mathbf{w}}} \ \sum_{i \in [n]} \tilde{g}_i \min_{G : i \in G \in \mathcal{G}} \tilde{w}_G, \quad \text{subject to} \ \ \tilde{\mathbf{w}} \geq 0, \ \sum_{G \in \mathcal{G}} b_G \tilde{w}_G = 1. \tag{87}$$

This is now just a piecewise linear objective over a (reweighted) simplex. Since projecting to a simplex can be performed in linear time, the smoothing method of (Nesterov, 2005) can be used to obtain a very efficient implementation. We illustrate a particular case where each variable $i \in [n]$ belongs to at most $r > 1$ groups.

**Proposition 17.** *Let $h(\tilde{\mathbf{w}})$ denote the* negated *objective of* (87)*. Then for any $\epsilon > 0$,*

$$h_\epsilon(\tilde{\mathbf{w}}) := \frac{\epsilon}{n \log r} \sum_{i \in [n]} \log \sum_{G : i \in G} r^{-n \tilde{g}_i \tilde{w}_G / \epsilon}$$

*enjoys the following properties*

1. *the gradient of $h_\epsilon$ is $\left(\frac{n}{\epsilon} \|\tilde{\mathbf{g}}\|_\infty^2 \log r\right)$-Lipschitz continuous,*

2. *$h(\tilde{\mathbf{w}}) - h_\epsilon(\tilde{\mathbf{w}}) \in (-\epsilon, 0]$ for all $\tilde{\mathbf{w}}$, and*

3. *the gradient of $h_\epsilon$ can be computed in $O(nr)$ time.*

The proof is given in Appendix C.2. With this construction, APG can be run on $h_\epsilon$ to achieve a $2\epsilon$ accurate solution to (87) within $O(\frac{1}{\epsilon}\sqrt{n \log r})$ steps (Nesterov, 2005), using a total time cost of $O(\frac{nr}{\epsilon}\sqrt{n \log r})$. For readability we defer the subtle issue of recovering an integral solution for (80) to Appendix C.3.

As a comparison to existing work, consider first (Mairal et al., 2011). The Algorithm 2 in (Mairal et al., 2011) proceeds in loops, with each iteration involving a max-flow problem on the canonical graph. The loop can take at most $n$ iterations, while each max-flow problem can be solved with $O(|V| |E|)$ cost where $|V|$ and $|E|$ are the number of nodes and edges in the canonical graph, respectively. By construction, $|V| = n + l$, and $|E| \le nr$ since each pair of $(G, i)$ with the node $i$ belonging to the group $G$ contributes an edge. Therefore the total cost is upper bounded by $O(n^2(n + l)r)$. Note that in the worst case $l = \Theta(nr)$. By contrast, the approach we developed above costs $O(\frac{nr}{\epsilon}\sqrt{n \log r})$, significantly cheaper in the regime where $n$ is big and $\epsilon$ is moderate. More importantly, we gain explicit control of the trade-off between accuracy and computational cost.

Another comparison is (Obozinski and Bach, 2012), which developed a line search scheme to compute the polar. The major computational step in that work is to solve

$$\tilde{\mathbf{w}}_\sigma = \arg\max_{\tilde{\mathbf{w}} \in Q} \langle \tilde{\mathbf{g}}, \tilde{\mathbf{w}} \rangle - \sigma \langle \mathbf{b}, \tilde{\mathbf{w}} \rangle, \tag{88}$$

recursively, each time with an updated $\sigma > 0$. In the case of bounded degree groups, this is again a max-flow problem, which costs $O(n(n + l)r)$ hence the overall cost is $O(n(n + l)r \log \frac{1}{\epsilon})$. This improves on the approach of (Mairal et al., 2011) but remains worse than the approach proposed here (when $n$ and $r$ are large and $\epsilon$ is not too small).

# 6. Experimental Evaluation

We evaluate the extended GCG approach developed above in four sparse learning problems: matrix completion, multi-class classification, multi-view dictionary learning, and structured sparse estimation. All algorithms were implemented in Matlab unless otherwise noted.[16] All experiments were run on a single core of NICTA's cluster with AMD 6-Core Opteron 4184 (2.8GHz, 3M L2/6M L3 Cache, 95W TDP, 64 GB memory).

## 6.1 Matrix Completion with Trace Norm Regularization

Our first set of experiment is on matrix completion (61), using five data sets whose statistics are given in Table 1. The data sets `Yahoo252M` (Dror et al., 2012) and `Yahoo60M`[17] provide music ratings, and the data sets `Netflix` (Dror et al., 2012), `MovieLens1M`, and `MovieLens10M`[18] all provide movie ratings. For `MovieLens1M` and `Yahoo60M`, we randomly sampled 75% ratings as training data, using the rest for testing. The other three data sets came with training and testing partition.

---

16. All code is available for download at http://users.cecs.anu.edu.au/~zhang/GCG.
17. `Yahoo60M` is the R1 data set at http://webscope.sandbox.yahoo.com/catalog.php?datatype=r.
18. Both MovieLens1M and MovieLens10M were downloaded from http://www.grouplens.org/datasets/movielens.

Table 1: Summary of the data sets used in matrix completion experiments. Here $n$ is the total number of users, $m$ is the number of products (movies or music), #train and #test are the number of ratings used for training and testing respectively.

| data set | $n$ (#user) | $m$ (#product) | #train | #test | $\lambda$ |
|---|---|---|---|---|---|
| MovieLens1M | 6,040 | 3,952 | 750,070 | 250,024 | 50 |
| MovieLens10M | 71,567 | 65,133 | 9,301,174 | 698,599 | 100 |
| Netflix | 2,549,429 | 17,770 | 99,072,112 | 1,408,395 | 1 |
| Yahoo60M | 98,130 | 1,738,442 | 60,021,577 | 20,872,219 | 300 |
| Yahoo252M | 1,000,990 | 624,961 | 252,800,275 | 4,003,960 | 100 |

We compared the extended GCG approach (Algorithm 3) with three state-of-the-art solvers for trace norm regularized objectives: MMBS[19] (Mishra et al., 2013), DHM (Dudik et al., 2012), and JS (Jaggi and Sulovsky, 2010). The local optimization in GCG was performed by L-BFGS with the maximum number of iteration set to 20.[20] The implementation of DHM is streamlined by buffering the quantities $\langle \mathcal{P}(\mathbf{u}_t \mathbf{v}_t'), \mathcal{P}(\mathbf{u}_s \mathbf{v}_s') \rangle$ and $\langle \mathcal{P}(\mathbf{u}_t \mathbf{v}_t'), \mathcal{P}(X) \rangle$, where $(\mathbf{u}_t, \mathbf{v}_t)$ is the result of polar operator at iteration $t$ (Line 3 of Algorithm 3). This allows totally corrective updates (51) to be performed efficiently, with no need to compute matrix inner products when evaluating the objective. JS was proposed for solving the constrained problem: $\min_W \ell(W)$ s.t. $\|W\|_{\mathrm{tr}} \leq \zeta$, which makes it hard to directly compare with solvers for the penalized problem (61). As a workaround, given $\lambda$ we first found the optimal solution $W^*$ for the penalized problem, and then we set $\zeta = \|W^*\|_{\mathrm{tr}}$ and finally solved the constrained problem by JS. Traditional solvers such as proximal methods (Toh and Yun, 2010; Pong et al., 2010) were not included in this comparison because they are much slower.

Figures 4 to 8 show how quickly the various algorithms drive the training objective down, while also providing the root mean square error (RMSE) on test data (i.e. comparing the reconstructed matrix each iteration to the ground truth). The regularization constant $\lambda$, given in Table 1, was chosen from $\{1, 10, 50, 100, 200, 300, 400, 500, 1000\}$ to minimize the test RMSE. On the MovieLens1M and MovieLens10M data sets, GCG and DHM exhibit similar efficiency in reducing the objective value with respect to CPU time (see Figure 4(a) and 5(a)). However, GCG achieves comparable test RMSE an order of magnitude faster than DHM, as shown in Figure 4(b) and 5(b). Interestingly, this discrepancy can be explained by investigating the rank of the solution, i.e. the number of outer iterations taken. In Figure 4(c) and 5(c), DHM clearly converges much more slowly in terms of the number of iterations, which is expected because it does not re-optimize the basis at each iteration, unlike GCG and MMBS. Therefore, although the overall computational cost for DHM appears on par with that for GCG, its solution has a much higher rank, resulting in a much higher test RMSE in these cases.

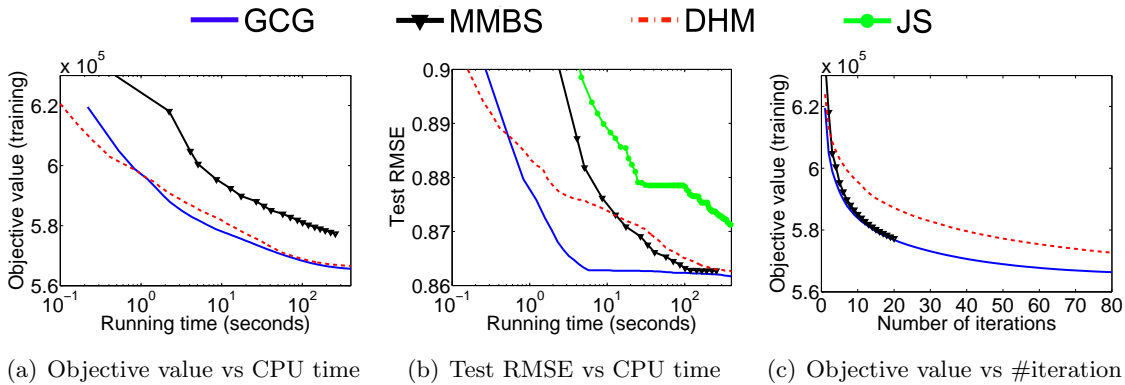

(a) Objective value vs CPU time     (b) Test RMSE vs CPU time     (c) Objective value vs #iteration

Figure 4: Matrix completion on `MovieLens1M`

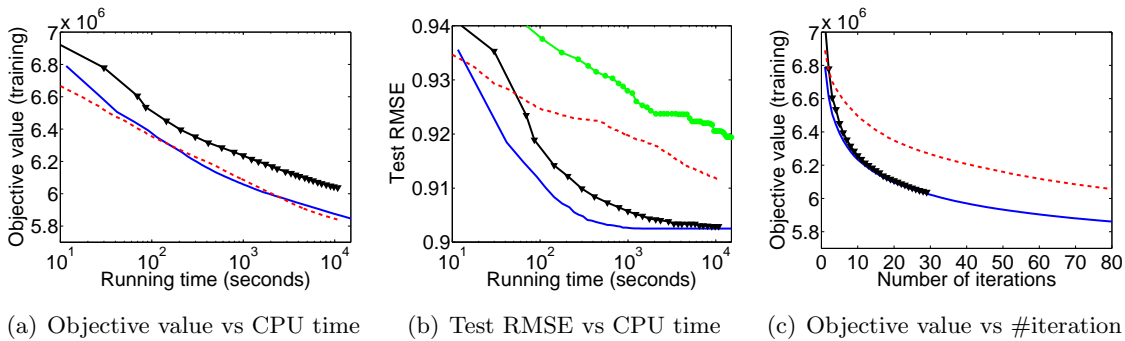

(a) Objective value vs CPU time     (b) Test RMSE vs CPU time     (c) Objective value vs #iteration

Figure 5: Matrix completion on `MovieLens10M`

MMBS takes nearly the same number of iterations as GCG to achieve the same objective value (see Figure 4(c) and 5(c)). However, its local search at each iteration is conducted on a *constrained* manifold, making it much more expensive than locally optimizing the unconstrained smooth surrogate (62). As a result, the overall computational complexity of MMBS is an order of magnitude greater than that of GCG, with respect to optimizing both the objective value and the test RMSE.

The observed discrepancies between the relative efficiency of optimizing the objective value and test RMSE becomes much smaller on the other three data sets, which are much larger in size. In Figure 6 to 8, it is clear that GCG takes significantly less CPU time to optimize both criteria. The declining trend of the objective value has a similar shape to that of the test RMSE. Here we observe that DHM is substantially slower than GCG and MMBS, confirming the importance of optimizing the bases at the same time as optimizing their weights. JS, which does not exploit the penalized form of the objective, is considerably slower with respect to reducing the test RMSE.

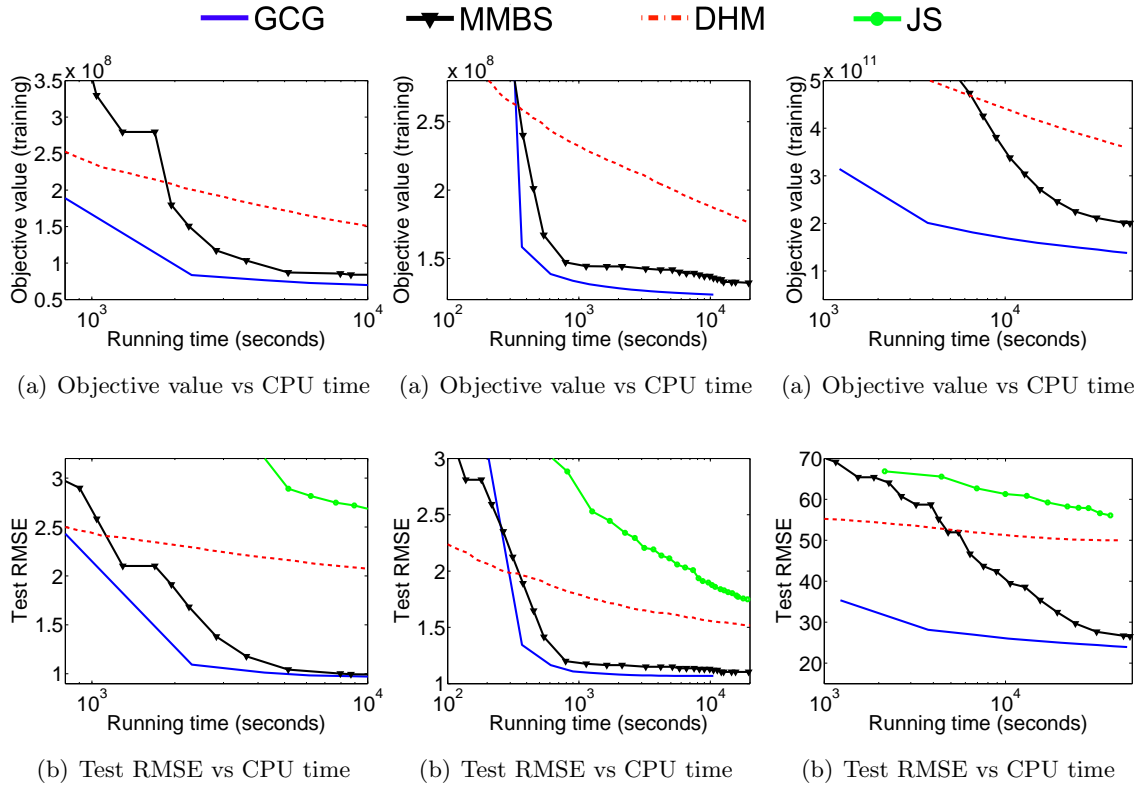

(a) Objective value vs CPU time  (a) Objective value vs CPU time  (a) Objective value vs CPU time

(b) Test RMSE vs CPU time  (b) Test RMSE vs CPU time  (b) Test RMSE vs CPU time

Figure 6: Matrix completion     Figure 7: Matrix completion     Figure 8: Matrix completion
on `Netflix`     on `Yahoo100m`     on `Yahoo252m`

## 6.2 Multi-class Classification with Trace Norm Regularization

Next we compared the four algorithms on multi-class classification problems also in the context of trace norm regularization. Here the task is to predict the class membership of input instances, $\mathbf{x}_i \in \mathbb{R}^n$, where each instance is associated with a unique label $y_i \in \{1, \ldots, C\}$. In particular, we consider a standard linear classification model where each class $c$ is associated with a weight vector in $\mathbb{R}^n$, stacked in a model matrix $W \in \mathbb{R}^{n \times C}$. Then for each training example $(\mathbf{x}_i, y_i)$, the individual logistic loss $L(W; \mathbf{x}_i, y_i)$ is given by

$$L(W; \mathbf{x}_i, y_i) = -\log \frac{\exp(\mathbf{x}_i' W_{:,y_i})}{\sum_{c=1}^{C} \exp(\mathbf{x}_i' W_{:c})},$$

and the complete loss is given by averaging over the training set, $\ell(W) := \frac{1}{m} \sum_{i=1}^{m} L(W; \mathbf{x}_i, y_i)$. It is reasonable to assume that $W$ has a low rank; i.e. the weight vector of all classes lie in a low dimensional subspace (Harchaoui et al., 2012), which motivates the introduction of a trace norm regularizer on $W$, yielding a training problem in the form of (56).

Among the comparison methods (MMBS, DHM and JS), since the original MMBS package is specialized to matrix completion, we implemented its extension to general loss functions, which required an optimized computational scheme for the Hessian-vector multiplication (the most time-consuming step), as outlined in Appendix D.

Table 2: Summary of the data sets used in multi-class classification experiments. Here $n$ is the number of features, $C$ is the number of classes, #train and #test are the number of training and test images respectively.

| data set | $n$ | $C$ | #train | #test | $\lambda$ |
|----------|-----|-----|--------|-------|-----------|
| Fungus10 | 4,096 | 10 | 10 per class | 10 per class | $10^{-2}$ |
| Fungus134 | 4,096 | 134 | 50 per class | 50 per class | $10^{-3}$ |
| k1024 | 131,072 | 50 | 50 per class | 50 per class | $10^{-3}$ |

We conducted experiments on three data sets extracted from the ImageNet repository (Deng et al., 2009), with characteristics shown in Table 2 (Akata et al., 2014).[21] k1024 is from the ILSVRC2011 benchmark, and Fungus10 and Fungus134 are from the Fungus group. All images are represented by the Fisher descriptor.

As in the matrix completion case above, we again compared how fast the algorithms reduced the objective value and test accuracy. Here the value of $\lambda$ was chosen from $\{10^{-4}, 10^{-3}, 10^{-2}, 10^{-1}\}$ to maximize test accuracy. Figures 9 to 11 show that the models produced by GCG achieve comparable (or better) training objectives and test accuracies in far less time than MMBS, DHM and JS. Although DHM is often more efficient than GCG early on, it is soon overtaken. This observation could be related to the diminishing return of adding bases, as shown in Figure 9(c), 10(c), and 11(c), since the addition of the first few bases yields a much steeper decrease in objective value than that achieved by later bases. Hence, after the initial stage, the computational efficiency of a solver turns to be dominated by better exploitation of the bases. That is, using added bases as an initialization for local improvement (as in GCG) appears to be more effective than simply finding an optimal conic combination of added bases (as in DHM).

Also notice from Figures 9(c), 10(c) and 11(c) that MMBS and GCG are almost identical in terms of the number of iterations required to achieve a given objective value. However, the manifold based local optimization employed by MMBS at each iteration makes its overall computational time worse than that of DHM and JS. This is a different outcome from the matrix completion experiments we saw earlier, and seems to suggest that the manifold approach in MMBS is more effective for the least square loss than for the logistic loss.

### 6.3 Multi-view Dictionary Learning

Our third set of experiments focuses on optimizing the objective for the multi-view dictionary learning model (67) presented in Section 4.4. Here we compared the GCG approach outlined in Section 4.4 with a straightforward *local* solver that is based on block coordinate descent (BCD), which alternates between: (a) fixing $H$ and optimizing $A$ and $B$ (with norm ball constraints), then (b) fixing $A$ and $B$ and optimizing $H$ (with nonsmooth regularizer). Both steps were solved by FISTA (Beck and Teboulle, 2009), initialized with the solution from the previous iteration. We tested several values of $t$, the dimensionality of the latent subspace of the solution, in these experiments. The GCG approach used Algorithm 3 with the local objective (62) instantiated by the norm expression in (68) and solved by L-BFGS with max iteration set to 20. We initialized the $U$ and $V$ matrices with 20 columns and rows respectively, each drawn i.i.d. from a uniform distribution on the unit sphere. Based on the noise model we consider in these experiments, we started with a base loss $\ell$ set to $l_1(Z, Z^*) = \sum_{ij} |Z_{ij} - Z_{ij}^*|$ where $Z^*$ is the true underlying matrix; however, following (Becker et al., 2009, Eq 3.1), we smoothed this $l_1$ loss by adding a quadratic prox-function with strong convexity modulus $10^{-4}$. Under the resulting loss, the optimization over $A$ and $B$ is decoupled given $H$. We terminated all algorithms after a maximum of running time (3,600 seconds) is reached.

---

21. The three data sets were all downloaded from http://lear.inrialpes.fr/src/jsgd/.

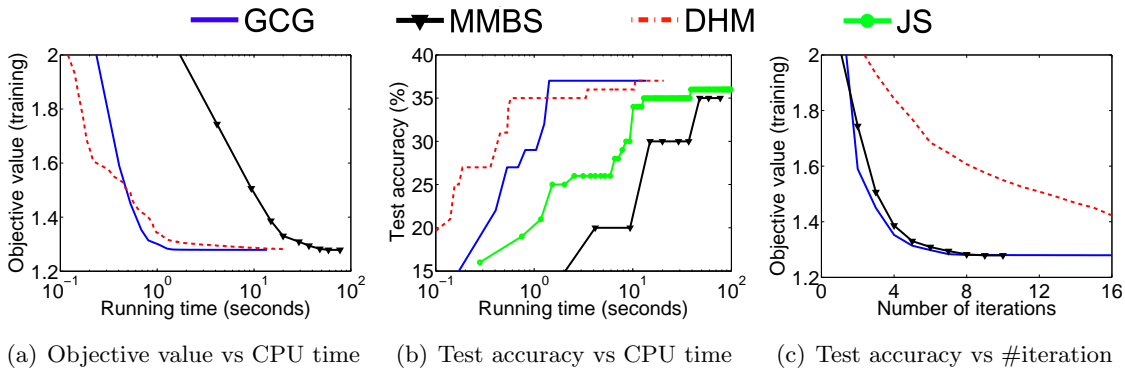

(a) Objective value vs CPU time   (b) Test accuracy vs CPU time   (c) Test accuracy vs #iteration

Figure 9: Multi-class classification on `Fungus10`

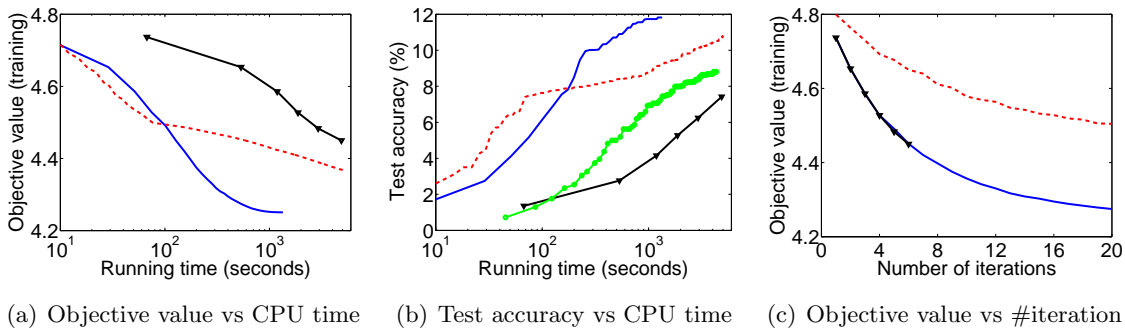

(a) Objective value vs CPU time   (b) Test accuracy vs CPU time   (c) Objective value vs #iteration

Figure 10: Multi-class classification on `Fungus134`

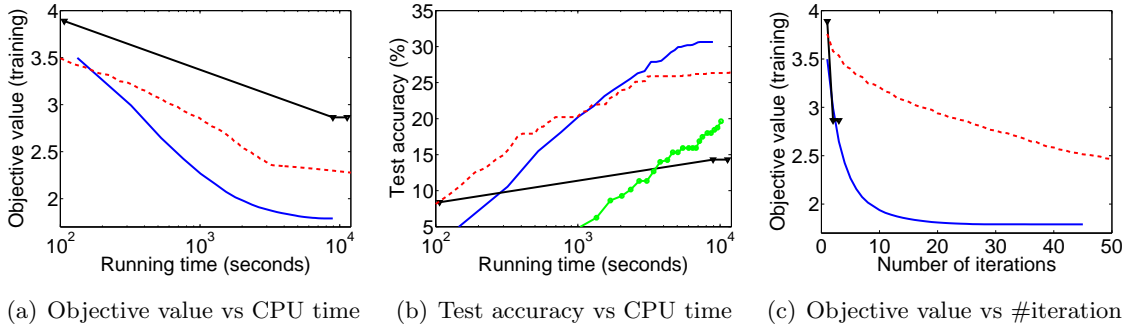

(a) Objective value vs CPU time   (b) Test accuracy vs CPU time   (c) Objective value vs #iteration

Figure 11: Multi-class classification on `k1024`

**Synthetic data.** We first constructed experiments on a synthetic data set that fulfills the conditional independence assumption (67) of the multi-view model (as discussed in Section 4.4). Here the $x$- and $y$-views have basis matrices $A^* \in \mathbb{R}^{n_1 \times t^*}$ and $B^* \in \mathbb{R}^{n_2 \times t^*}$ respectively, where all columns of $A^*$ were drawn i.i.d. from the $l_2$ sphere of radius $\beta = 1$, and all columns of $B^*$ were drawn i.i.d. from the $l_2$ sphere of radius $\gamma = 5$. We set $n_1 = 800$, $n_2 = 1000$, and $t^* = 10$ (so that the dictionary size is indeed small). The latent representation matrix $H^* \in \mathbb{R}^{t^* \times m}$ has all elements drawn i.i.d. from the zero-mean unit-variance Gaussian distribution. The number of training examples is set to $m = 200$, which is much lower than the number of features $n_1$ and $n_2$. Then the clean versions of $x$-view and $y$-view were generated by $X^* = A^* H^*$ and $Y^* = B^* H^*$ respectively, while the noisy versions were obtained by adding i.i.d. noise to 15% entries of $X^*$ and $Y^*$ that were selected uniformly at random, yielding $\tilde{X}$ and $\tilde{Y}$ respectively. The noise is uniformly distributed in $[0, 10]$. Given the noisy observations, the denoising task is to find reconstructions $X$ and $Y$ for the two views, such that the error $\|X - \tilde{X}\|_{\mathrm{F}}^2 + \|Y - \tilde{Y}\|_{\mathrm{F}}^2$ is small. The composite training problem is formulated in (67).

The results of comparing the modified GCG method with local optimization are presented in Figures 12 to 14, where the regularization parameter $\lambda$ has been varied among $\{80, 60, 40\}$. In each case, GCG optimizes the objective value significantly faster than BCD for all values of $t$ chosen for BCD. The differences between methods are even more evident when considering the reconstruction error. In particular, although the local optimization in GCG already achieves a low objective value in the first outer iteration (see Figure 12(a) and 13(a)), the reconstruction error remains high, thus requiring further effort to reduce it.

When $t = 10$ which equals the true rank $t^* = 10$, BCD is more likely to get stuck in a poor local minimum when $\lambda$ is large, for example as shown in Figure 12(a) for $\lambda = 80$. This occurs because when optimizing $H$, the strong shrinkage in the proximal update stifles major changes in $H$ (and in $A$ and $B$ consequently). This problem is exacerbated when the rank of the solution is overly restricted, leading to even worse local optima. On the other hand, under-regularization does indeed cause overfitting in this case; for example, as shown in Figure 14 for $\lambda = 40$. Interestingly, BCD gets trapped in local optima for all $t$, while GCG eventually escapes suboptimal solutions and eventually finds a way to considerably reduce the objective value. However, this simply creates a dramatic *increase* in the test reconstruction error, a typical phenomenon of overfitting.

**Image denoising.** We then conducted multi-view dictionary learning experiments on natural image data. The data set we considered is based on the Extended Yale Face Database B[22] (Georghiades et al., 2001). In particular, the data set we used consists of grey level face images of 28 human subjects, each with 9 poses and 64 lighting conditions. To construct the data set, we set the $x$-view to a fixed lighting (+000E+00) and the $y$-view to a different fixed lighting (+000E+20). We obtained paired views by randomly drawing a subject and a pose under the two lighting conditions. The underlying assumption is that each lighting has its own set of bases ($A$ and $B$) and each (person, pose) pair has the same latent representation for the two lighting conditions. All images were down-sampled to 50-by-50, meaning

---

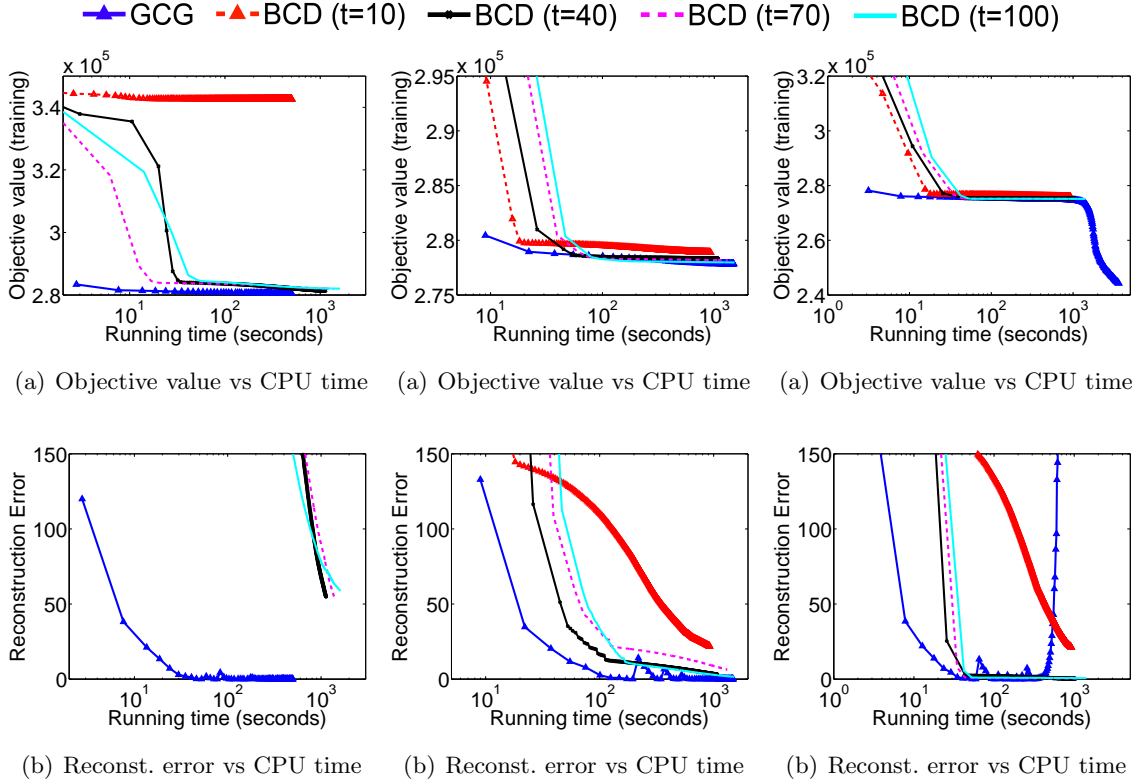

(a) Objective value vs CPU time    (a) Objective value vs CPU time    (a) Objective value vs CPU time

(b) Reconst. error vs CPU time    (b) Reconst. error vs CPU time    (b) Reconst. error vs CPU time

Figure 12: Denoising on    Figure 13: Denoising on    Figure 14: Denoising on
synthetic data    synthetic data    synthetic data
$\lambda = 80$    $\lambda = 60$    $\lambda = 40$

$n_1 = n_2 = 2500$, and we randomly selected $t = 50$ pairs of (person, pose) for training. The $x$-view was kept clean, while pixel errors were added to the $y$-view, randomly setting 10% of the pixel values to 1. This noise model mimics natural phenomena such as occlusion and loss of pixel information from image transfer. The goal is again to enable appropriate reconstruction of a noisy image using other views. Naturally, we set $\beta = \gamma = 1$ since the data is in the $[0, 1]$ interval.

Figures 15 to 17 show the results comparing GCG with the local BCD for different choices of $t$ for BCD. Clearly the value of $\lambda$ that yields the lowest reconstruction error is 40. However, here BCD simply gets stuck at the trivial local minimum $H = \mathbf{0}$, as shown in Figure 15(a). This occurs because under the random initialization of $A$ and $B$, the proximal update simply sets $H$ to $\mathbf{0}$ when the regularization is too strong. When $\lambda$ is decreased to 30, as shown in Figure 16(a), BCD can escape the trivial $H = \mathbf{0}$ solution for $t = 100$, but again it becomes trapped in another local minimum that remains much worse than the globally optimal solution found by GCG. When $\lambda$ is further reduced to 20, the BCD algorithms finally escape the trivial $H = \mathbf{0}$ point for all $t$, but still fail to find an acceptable solution. Although GCG again converges to a much better global solution in terms of the objective value, overfitting can still occur: in Figure 17(b) the reconstruction error of GCG increases significantly after an initial decline.

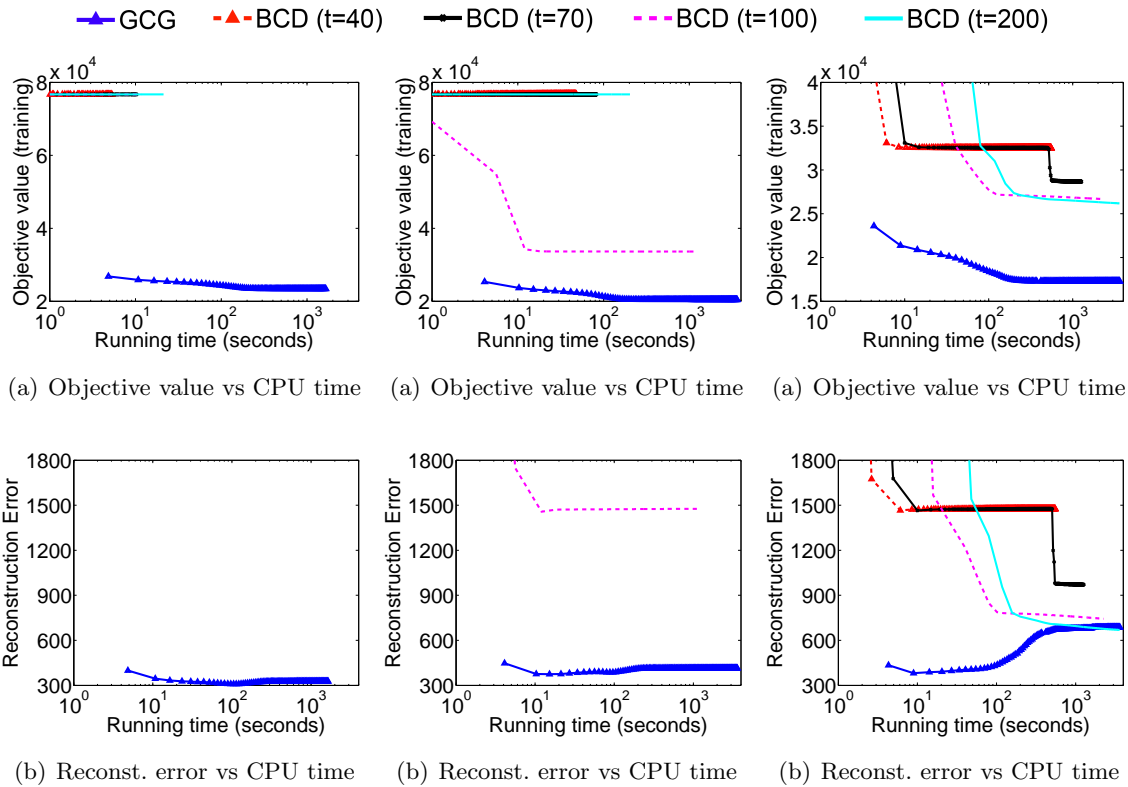

(a) Objective value vs CPU time     (a) Objective value vs CPU time     (a) Objective value vs CPU time

(b) Reconst. error vs CPU time     (b) Reconst. error vs CPU time     (b) Reconst. error vs CPU time

Figure 15: Denoising on     Figure 16: Denoising on     Figure 17: Denoising on
            face data                     face data                   face data
            $\lambda = 40$                       $\lambda = 30$                     $\lambda = 20$

## 6.4 Structured Sparse Estimation: CUR-like Matrix Factorization

Finally, we study the empirical performance of GCG on optimizing structured sparse models of the form discussed in Section 5. We considered an example of group lasso that is inspired by CUR matrix factorization (Mahoney and Drineas, 2009). Given a data matrix $X \in \mathbb{R}^{n \times d}$, the goal is to compute an approximate factorization $X \approx CUR$, such that $C$ contains a subset of $c$ columns from $X$ and $R$ contains a subset of $r$ rows from $X$. Mairal et al. (2011, §5.3) proposed a convex relaxation of this problem:

$$\min_W \frac{1}{2} \|X - XWX\|_F^2 + \lambda \left( \sum_i \|W_{i:}\|_\infty + \sum_j \|W_{:j}\|_\infty \right). \tag{89}$$

Conveniently, the regularizer fits the development of Section 5.3, with $p = 1$ and the groups defined to be the rows and columns of $W$. To evaluate different methods, we used four gene-expression data sets: SRBCT, Brain_Tumor_2, 9_Tumor, and Leukemia2, of sizes $83 \times 2308$, $50 \times 10367$, $60 \times 5762$ and $72 \times 11225$ respectively.[23] The data matrices were first centered columnwise and then rescaled to have unit Frobenius norm.

We compared three algorithms: GCG_TUM, consisting of the GCG Algorithm 4 with the polar operator developed in Section 5; GCG_Secant, consisting of GCG with the polar operator of (Mairal et al., 2011, Algorithm 2); and APG. The polar operator in GCG_TUM used L-BFGS to find an optimal $w^*_{cr}$ for the smoothed version of (87) given in Proposition 17, with smoothing parameter $\epsilon$ set to $10^{-3}$. The polar operator in GCG_Secant was implemented with a mex wrapper of a max-flow package.[24] Both GCG methods use a totally corrective update, which allows efficient optimization by L-BFGS-B via pre-computing $X\mathbb{P}_{J_{\bar{g}}}(\mathbf{g}_k)X$. The PU in APG uses the routine mexProximalGraph from the SPAMS package.[25]

**Results.** We tested with values for $\lambda$ chosen from the set $\{3 \cdot 10^{-3}, 10^{-3}, 10^{-4}\}$, which lead to increasingly dense solutions (for $\lambda = 10^{-2}$ it turns out that the optimal solution is $W = \mathbf{0}$, hence only smaller values of $\lambda$ are interesting). Figure 18 shows that $\lambda = 10^{-3}$ yields moderately sparse solutions; results for $\lambda = 3 \cdot 10^{-3}$ and $10^{-4}$ are given in Figures 19 and 20 respectively. In these experiments, GCG_TUM proves to be an order of magnitude faster than GCG_Secant in computing the polar. Although as (Mairal et al., 2011) observe, network flow algorithms often find solutions far more quickly in practice than what their theoretical bounds suggest, given the efficiency of the totally corrective update, nearly all of the computation of GCG_Secant is still devoted to the polar operator. Here, the more efficient polar operator used by GCG_TUM leads to a reduction of *overall* optimization time by at least 50%. Finally, APG is always slower than GCG_Secant by an order of magnitude, with the PU consuming the vast majority of the computation time.

## 7. Conclusion

We have presented a unified treatment of the generalized conditional gradient (GCG) algorithm that covers several useful extensions. After a thorough treatment of its convergence properties, we illustrated how an extended form of GCG can efficiently handle a gauge regularizer in particular. We illustrated the application of these GCG methods in two general application areas: low rank learning and structured sparse estimation. In each case, we considered a generic convex relaxation procedure and focused on the efficient computation of the polar operator—one of the key steps to obtaining an efficient implementation of GCG. To further accelerate its convergence, we interleaved the standard GCG update with a fixed rank subspace optimization, which greatly improves performance in practice without affecting the theoretical convergence properties. Extensive experiments on matrix completion, multi-class classification, multi-view learning and overlapping group lasso confirmed the superiority of the modified GCG approach in practice.

For future work, several interesting extensions are worth investigating. The current form of GCG is inherently framed as a batch method, where each iteration requires gradient of the loss to be computed over *all* training data. When the training data is large, it is preferable to only use a random subset, or distribute the training data across several processors. Another interesting direction is to interpolate between the polar operator and proximal updates; that is, one could incorporate a mini-batch of atoms at each iteration. For example, for trace norm regularization, one can consider adding $k$ instead of just the 1 largest singular vectors. It remains an open question whether such a "mini-batch" of atoms can provably accelerate the overall optimization, but some recent work (Hsieh and Olsen, 2014) has demonstrated promise in practice. Finally, many multi-structured estimation problems employ multiple nonsmooth regularizers (Richard et al., 2013), where algorithms have been developed for computing the associated PU e.g. (Bach et al., 2012a; Yu, 2013a). Developing efficient polar operators for these cases has not yet been investigated, but could prove to be very useful in practice.

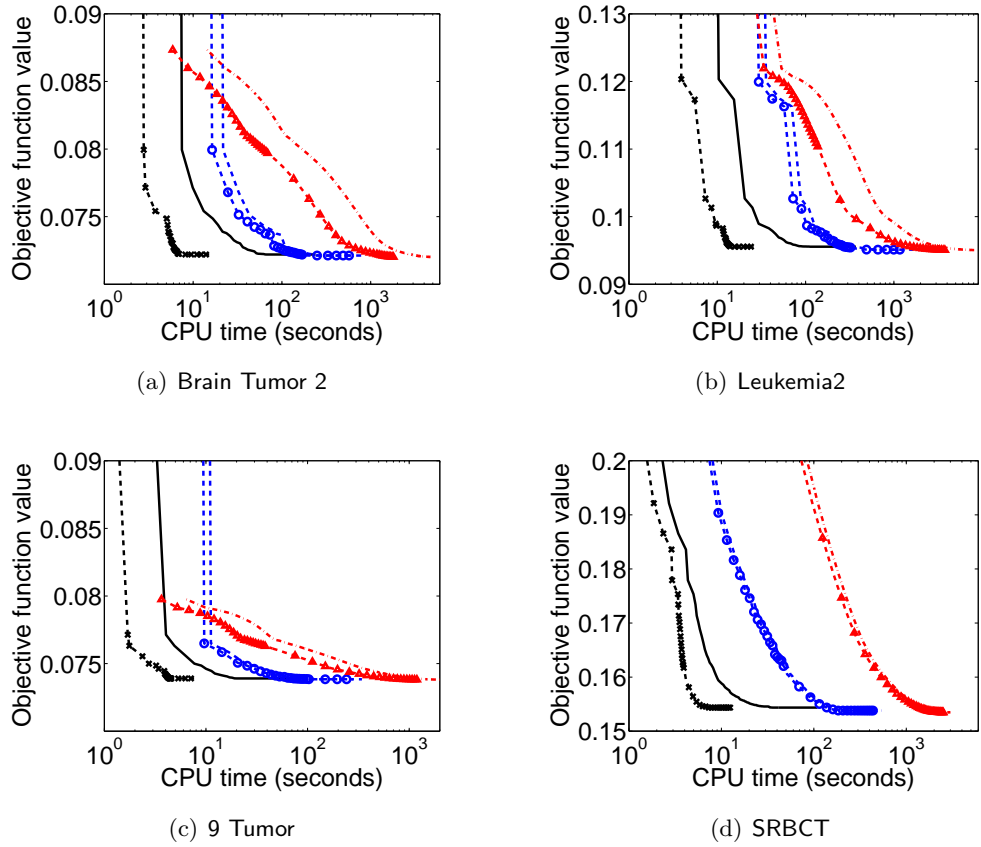

Figure 18: Results for convex CUR matrix factorization with $\lambda = 10^{-3}$.

# Acknowledgments

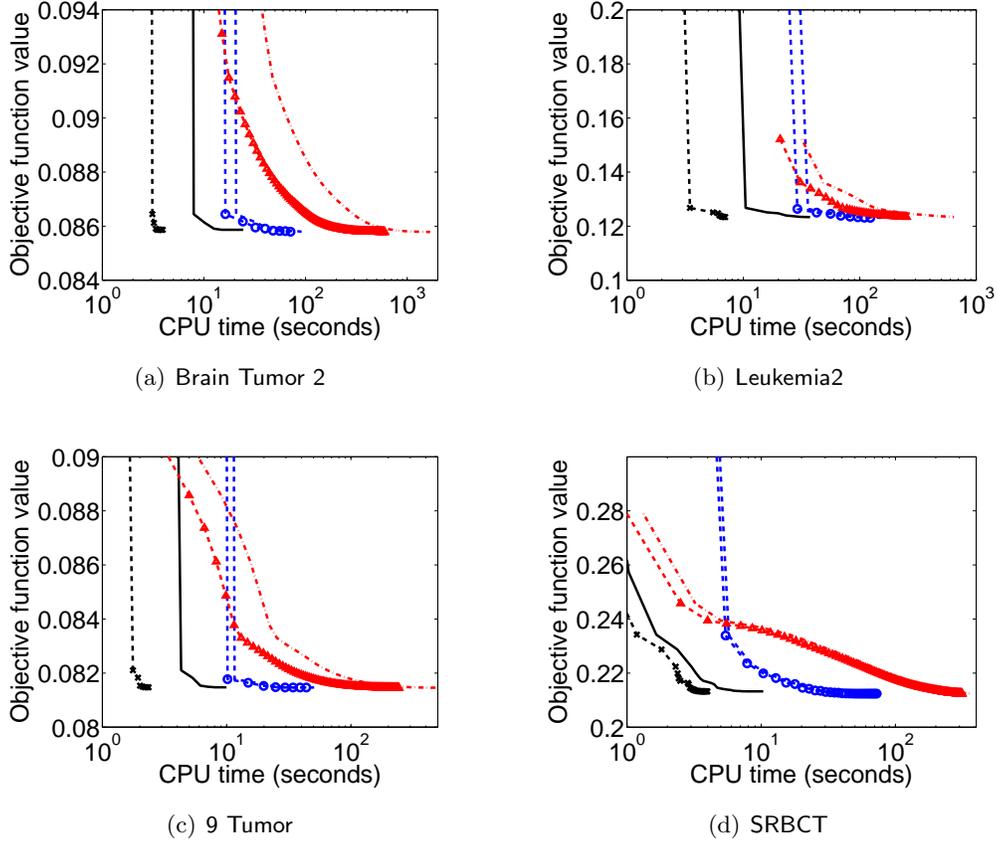

Figure 19: Results for convex CUR matrix factorization with $\lambda = 3 \cdot 10^{-3}$

This research was supported by the Alberta Innovates Center for Machine Learning (AICML), the Canada Research Chairs program, and NSERC. Most of this work was carried out when Yaoliang Yu and Xinhua Zhang were both with the University of Alberta and AICML. NICTA is funded by the Australian Government as represented by the Department of Broadband, Communications and the Digital Economy and the Australian Research Council through the ICT Centre of Excellence program. We thank Csaba Szepesvári for many constructive suggestions and also thank the reviewers who provided helpful comments on earlier versions of (parts of) this paper.

# Appendix A. Proofs for Section 3

## A.1 Proof of Proposition 3

(a): Let $C$ be the closure of the sequence $\{\mathbf{w}_t\}_t$. Due to the compactness assumption on the sublevel set and the monotonicity of $F(\mathbf{w}_t)$, $C$ is compact. Moreover $C \subseteq \operatorname{dom} f$ since for all cluster points, say $\mathbf{w}$, of $\mathbf{w}_t$ we have from the closedness of $F$ that $F(\mathbf{w}) \le$ $\liminf F(\mathbf{w}_{t_k}) \leq F(\mathbf{w}_0) < \infty$. Since $-\nabla \ell$ is continuous, $-\nabla \ell(C)$ is a compact subset of $\mathrm{int}(\mathrm{dom}\, f^*)$. Note that $f^*$ is continuous on the interior of its domain, therefore its subdifferential is locally bounded on $-\nabla \ell(C)$, see e.g. (Borwein and Vanderwerff, 2010, Proposition 4.1.26). A standard compactness argument then establishes the boundedness of $(\partial f^*)(-\nabla \ell(C))$. Thus $\{\mathbf{d}_t\}_t$ is bounded.

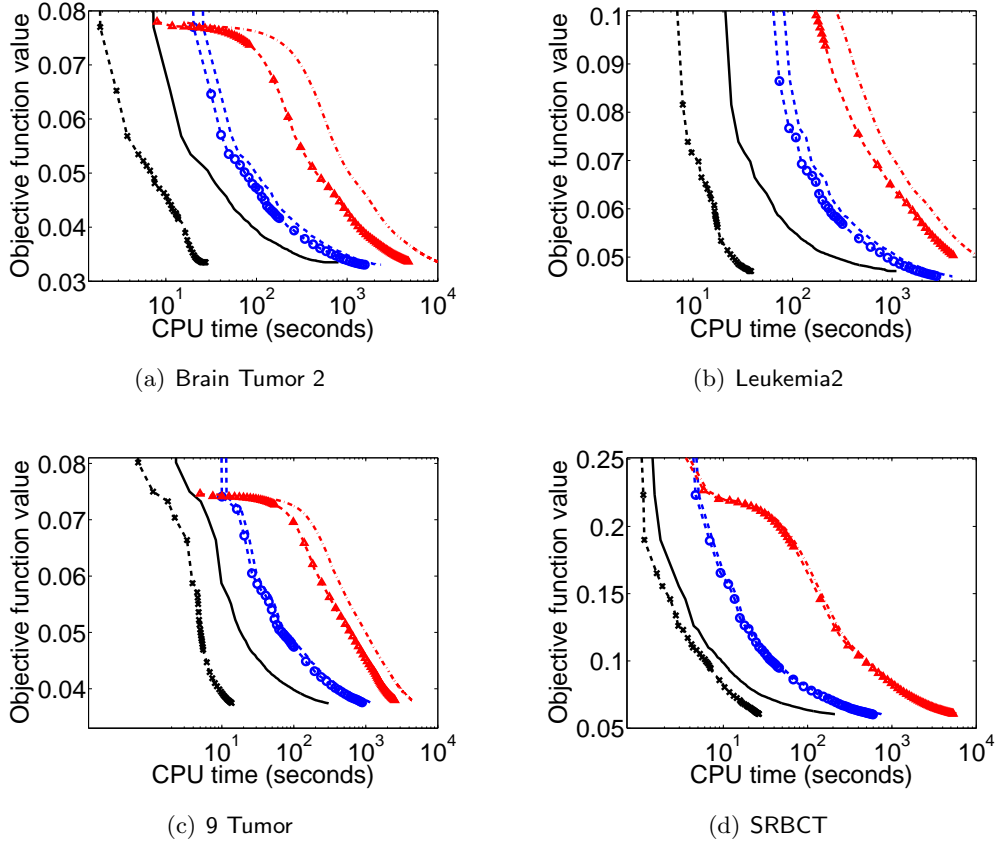

Figure 20: Results for convex CUR matrix factorization with $\lambda = 10^{-4}$

(b): Note first that the boundedness of $\{\mathbf{w}_t\}_t$ follows immediately from the boundedness assumption on the sublevel set, thanks to the monotonic property of $F(\mathbf{w}_t)$. Since $\nabla \ell$ is uniformly continuous, the set $\{-\nabla \ell(\mathbf{w}_t)\}_t$ is again bounded. On the other hand, we know from (Borwein and Vanderwerff, 2010, Theorem 4.4.13, Proposition 4.1.25) that $f$ is super-coercive iff $\partial f^*$ maps bounded sets into bounded sets. Therefore $\{\mathbf{d}_t\}_t$ is again bounded.

(c): This is clear.

## A.2  Proof of Proposition 4

By definition,

$$\mathsf{G}(\mathbf{w}) = F(\mathbf{w}) - \inf_{\mathbf{d}} \left\{ \ell(\mathbf{w}) + \langle \mathbf{d} - \mathbf{w}, \nabla \ell(\mathbf{w}) \rangle + f(\mathbf{d}) \right\}.$$

Letting $\mathbf{d} = \mathbf{w}$ we have $\mathsf{G}(\mathbf{w}) \geq 0$, whereas equality is attained iff $\mathbf{d} = \mathbf{w}$ is a minimizer of the convex problem

$$\inf_{\mathbf{d}} \left\{ \ell(\mathbf{w}) + \langle \mathbf{d} - \mathbf{w}, \nabla \ell(\mathbf{w}) \rangle + f(\mathbf{d}) \right\}.$$

Note that minimality with respect to $\mathbf{d}$ is achieved iff $\mathbf{0} \in \nabla \ell(\mathbf{w}) + \partial f(\mathbf{d})$; hence $\mathbf{d} = \mathbf{w}$ is optimal iff $\mathbf{0} \in \nabla \ell(\mathbf{w}) + \partial f(\mathbf{w})$; i.e. the necessary condition (7).

When $\ell$ is additionally convex, it follows that for all $\mathbf{w} \in \mathcal{B}$

$$\inf_{\mathbf{d}} \left\{ \ell(\mathbf{w}) + \langle \mathbf{d} - \mathbf{w}, \nabla \ell(\mathbf{w}) \rangle + f(\mathbf{d}) \right\} \leq \inf_{\mathbf{d}} \left\{ \ell(\mathbf{d}) + f(\mathbf{d}) \right\}, \tag{90}$$

hence

$$\begin{aligned}
G(\mathbf{w}) &= F(\mathbf{w}) - \inf_{\mathbf{d}} \left\{ \ell(\mathbf{w}) + \langle \mathbf{d} - \mathbf{w}, \nabla \ell(\mathbf{w}) \rangle + f(\mathbf{d}) \right\} \\
&\geq F(\mathbf{w}) - \inf_{\mathbf{d}} \left\{ \ell(\mathbf{d}) + f(\mathbf{d}) \right\} \\
&= F(\mathbf{w}) - \inf_{\mathbf{d}} F(\mathbf{d}).
\end{aligned}$$

### A.3 Proof of Theorem 5

Due to the step size rule, the subroutine is also `Monotone`, hence Algorithm 1 always makes monotonic progress. We now strengthen this observation by looking more carefully at the step size rule. Since the subroutine is `Descent`, we have

$$\begin{aligned}
F(\mathbf{w}_{t+1}) &\leq F(\tilde{\mathbf{w}}_{t+1}) = F((1 - \eta_t)\mathbf{w}_t + \eta_t \mathbf{d}_t) \\
&\leq \min_{0 \leq \eta \leq 1} \ell((1 - \eta)\mathbf{w}_t + \eta \mathbf{d}_t) + (1 - \eta) f(\mathbf{w}_t) + \eta f(\mathbf{d}_t) \\
&= \min_{0 \leq \eta \leq 1} \ell(\mathbf{w}_t) + \eta \langle \mathbf{d}_t - \mathbf{w}_t, \nabla \ell(\mathbf{u}_t) \rangle + (1 - \eta) f(\mathbf{w}_t) + \eta f(\mathbf{d}_t) \\
&= \min_{0 \leq \eta \leq 1} F(\mathbf{w}_t) + \eta [\langle \mathbf{d}_t - \mathbf{w}_t, \nabla \ell(\mathbf{u}_t) \rangle - f(\mathbf{w}_t) + f(\mathbf{d}_t)], \tag{91}
\end{aligned}$$

where, using the mean value theorem, $\mathbf{u}_t$ is some vector lying between $\mathbf{w}_t$ and $(1-\eta)\mathbf{w}_t + \eta \mathbf{d}_t$. As $\eta \to 0$, $\mathbf{u}_t \to \mathbf{w}_t$, hence by continuity, $\langle \mathbf{d}_t - \mathbf{w}_t, \nabla \ell(\mathbf{u}_t) \rangle - f(\mathbf{w}_t) + f(\mathbf{d}_t) \to -\mathsf{G}(\mathbf{w}_t) \leq 0$ (recall from Proposition 4 that we always have $\mathsf{G}(\mathbf{w}) \geq 0$ for all $\mathbf{w}$). If $\mathsf{G}(\mathbf{w}_t) = 0$ we have nothing to prove, otherwise we have $F(\mathbf{w}_{t+1}) < F(\mathbf{w}_t)$ for all sufficiently small $\eta > 0$. This completes the proof for the first claim.

To prove the second claim, we clearly can assume that $\mathsf{G}(\mathbf{w}_t) \neq 0$ for all $t$ sufficiently large and that $F(\mathbf{w}_t)$ converges to a finite limit. Rearrange (91):

$$F(\mathbf{w}_t) - F(\mathbf{w}_{t+1}) \geq \max_{0 \leq \eta \leq 1} \eta [\mathsf{G}(\mathbf{w}_t) - \langle \mathbf{d}_t - \mathbf{w}_t, \nabla \ell(\mathbf{u}_t) - \nabla \ell(\mathbf{w}_t) \rangle]. \tag{92}$$

By assumption $\{\mathbf{d}_t\}_t, \{\mathbf{w}_t\}_t$ are bounded, thus $\{\mathbf{u}_t\}_t$ is bounded too (as $\mathbf{u}_t$ is some convex combination of $\mathbf{d}_t$ and $\mathbf{w}_t$). On the other hand, since $\nabla \ell$ is assumed to be uniformly continuous on bounded sets, when $\eta$ is sufficiently small, say $\eta = \tilde{\eta} > 0$, we have $\mathbf{u}_t$ sufficiently close to $\mathbf{w}_t$ such that

$$\langle \mathbf{d}_t - \mathbf{w}_t, \nabla \ell(\mathbf{u}_t) - \nabla \ell(\mathbf{w}_t) \rangle \leq \|\mathbf{d}_t - \mathbf{w}_t\| \ \|\nabla \ell(\mathbf{u}_t) - \nabla \ell(\mathbf{w}_t)\|^{\circ} \leq \epsilon,$$

for some (arbitrary) $\epsilon > 0$. Crucially, $\tilde{\eta}$ does not depend on $t$ (although it may depend on $\epsilon$), thanks to the boundedness and the uniform continuity assumption. Therefore for $t$ sufficiently large we have from (92) that

$$\tilde{\eta}[\mathsf{G}(\mathbf{w}_t) - \epsilon] \leq F(\mathbf{w}_t) - F(\mathbf{w}_{t+1}) \leq \epsilon\tilde{\eta},$$

where the last inequality is because $F(\mathbf{w}_t)$ converges to a finite limit. Thus $\mathsf{G}(\mathbf{w}_t) \leq 2\epsilon$. Since $\epsilon$ is arbitrary, $\mathsf{G}(\mathbf{w}_t) \to 0$ and the proof is complete.

### A.4 Proof of Theorem 6

Not surprisingly the step size is chosen to minimize the quadratic upper bound:

$$F(\mathbf{w}_{t+1}) \leq \min_{\eta \in [0,1]} \ell(\mathbf{w}_t) + \eta \langle \mathbf{d}_t - \mathbf{w}_t, \nabla\ell(\mathbf{w}_t) \rangle + \frac{L\eta^2}{2} \|\mathbf{d}_t - \mathbf{w}_t\|^2 + (1-\eta)f(\mathbf{w}_t) + \eta f(\mathbf{d}_t)$$

$$= \min_{\eta \in [0,1]} F(\mathbf{w}_t) - \eta\mathsf{G}(\mathbf{w}_t) + \frac{L\eta^2}{2} \|\mathbf{w}_t - \mathbf{d}_t\|^2.$$

Therefore if $G(\mathbf{w}_t) > 0$, through minimizing $\eta$ in the above we have $F(\mathbf{w}_{t+1}) < F(\mathbf{w}_t)$. On the other hand, if $\mathsf{G}(\mathbf{w}_t) = 0$ for some $t$, then $\eta_t = 0$, resulting in $\mathbf{w}_{t+1} = \mathbf{w}_t$. Thus the algorithm will not change its iterate afterwards.

Assume that $\mathsf{G}(\mathbf{w}_t) \neq 0$ for any $t$ and that $F(\mathbf{w}_t)$ converges to a finite limit (otherwise there is nothing to prove). Analyzing the step size in each case separately, we have

$$\min\left\{\frac{\mathsf{G}^2(\mathbf{w}_t)}{2L\|\mathbf{w}_t - \mathbf{d}_t\|^2}, \frac{\mathsf{G}(\mathbf{w}_t)}{2}\right\} \leq F(\mathbf{w}_t) - F(\mathbf{w}_{t+1}) \to 0.$$

Due to the boundedness assumption in Assumption 3, we know $\mathsf{G}(\mathbf{w}_t) \to 0$.

### A.5 Proof of Theorem 7

Since the subroutine is `Relaxed`,

$$F(\mathbf{w}_{t+1}) \leq \ell(\mathbf{w}_t) + \eta_t \langle \mathbf{d}_t - \mathbf{w}_t, \nabla\ell(\mathbf{w}_t) \rangle + \frac{L\eta_t^2}{2} \|\mathbf{d}_t - \mathbf{w}_t\|^2 + (1-\eta_t)f(\mathbf{w}_t) + \eta_t f(\mathbf{d}_t)$$

$$= F(\mathbf{w}_t) - \eta_t\mathsf{G}(\mathbf{w}_t) + \eta_t\varepsilon_t + \frac{L\eta_t^2}{2} \|\mathbf{w}_t - \mathbf{d}_t\|^2.$$

Summing the indices from $s = 0$ to $s = t$ leads to (19).

The second claim follows simply from the observation that $\sum_t \eta_t\mathsf{G}(\mathbf{w}_t)$ is bounded from above under the given assumptions. Specifically, note that the divergence of the partial sum $H_t := \sum_{s=0}^t \eta_s \to \infty$ implies that $\sum_t \eta_t/H_t \to \infty$ and $\sum_t \eta_t/H_t^{1+\delta} < \infty$ for any $\delta > 0$, see e.g. (Hardy et al., 1952, Result 162).

### A.6 Proof of Theorem 8

Define $D = \text{dom } f$, which is closed and convex by assumption, and define $T(\mathbf{w}) := \partial f^*(\nabla\ell(\mathbf{w}))$. According to Proposition 2, the assumption of $f$ being $L$-strongly convex actually implies that $f^*$ is differentiable with a $1/L$-Lipschitz continuous gradient (note that $f \in \Gamma_0$ implies that $(f^*)^* = f$). Thus $T : D \to D$, being the composition of an $L$-Lipschitz continuous function $\nabla\ell$ and a $1/L$-Lipschitz continuous function $\nabla f^*$, is nonexpansive. In the motivation for GCG (given at the beginning of Section 3.1), it was pointed out that it merely consists of the fixed-point iteration $\mathbf{w}_{t+1} = (1 - \eta_t)\mathbf{w}_t + \eta_t T(\mathbf{w}_t)$. Therefore the claim follows immediately from the well-known Krasnosel'skiĭ-Mann theorem, see e.g. (Bauschke and Combettes, 2011, Theorem 5.14).

### A.7 Proof of Theorem 11

Note first that with the step size $\eta_t = 2/(t+2)$,

$$\mathbf{w}_{t+1} = \left(1 - \frac{2}{t+2}\right)\mathbf{w}_t + \frac{2}{t+2}\mathbf{d}_t = \frac{2}{(t+1)(t+2)}\sum_{s=0}^{t}(s+1)\mathbf{d}_s. \tag{93}$$

Using Jensen's inequality and the fact $\mathbf{d}_s \in \partial f^*(-\mathbf{g}_s)$ we have:

$$\begin{aligned}
\frac{(t+1)(t+2)}{2}[\ell^*(\bar{\mathbf{g}}_{t+1}) + f^*(-\bar{\mathbf{g}}_{t+1}) &- \ell^*(\mathbf{g}) - f^*(-\mathbf{g})] \\
&\leq \sum_{s=0}^{t}(s+1)[\ell^*(\mathbf{g}_s) + f^*(-\mathbf{g}_s) - \ell^*(\mathbf{g}) - f^*(-\mathbf{g})] \\
&\leq \sum_{s=0}^{t}(s+1)[\langle\mathbf{g}_s - \mathbf{g}, -\mathbf{d}_s\rangle + \ell^*(\mathbf{g}_s) - \ell^*(\mathbf{g})] \\
&= (t+1)[\langle\mathbf{g}_t - \mathbf{g}, -\mathbf{d}_t\rangle + \ell^*(\mathbf{g}_t) - \ell^*(\mathbf{g})] \\
&\quad + \sum_{s=0}^{t-1}(s+1)[\langle\mathbf{g}_s - \mathbf{g}, -\mathbf{d}_s\rangle + \ell^*(\mathbf{g}_s) - \ell^*(\mathbf{g})].
\end{aligned} \tag{94}$$

Conveniently, the term in (94), after some tedious algebra, simplifies to

$$(t+1)\langle\mathbf{w}_t - \mathbf{d}_t, \mathbf{g}_t - \mathbf{g}_{t+1}\rangle - \frac{(t+1)(t+2)}{2}\mathsf{D}_{\ell^*}(\mathbf{g}_{t+1}, \mathbf{g}_t) - \frac{(t+1)(t+2)}{2}\mathsf{D}_{\ell^*}(\mathbf{g}, \mathbf{g}_{t+1}) + \frac{t(t+1)}{2}\mathsf{D}_{\ell^*}(\mathbf{g}, \mathbf{g}_t),$$

where $\mathsf{D}_f(\mathbf{x}, \mathbf{y}) = f(\mathbf{x}) - f(\mathbf{y}) - \langle\mathbf{x} - \mathbf{y}, \nabla\partial f(\mathbf{y})\rangle$ is the Bregman divergence induced by $f$. Since $\nabla\ell$ is $L$-Lipschitz continuous, $\ell^*$ is $1/L$-strongly convex (cf. Proposition 2), hence $\mathsf{D}_{\ell^*}(\mathbf{g}_{t+1}, \mathbf{g}_t) \geq \frac{1}{2L}\|\mathbf{g}_{t+1} - \mathbf{g}_t\|_\circ^2$. Applying the Cauchy-Schwarz inequality we have

$$\begin{aligned}
\langle\mathbf{w}_t - \mathbf{d}_t, \mathbf{g}_t - \mathbf{g}_{t+1}\rangle - \frac{(t+2)}{2}\mathsf{D}_{\ell^*}(\mathbf{g}_{t+1}, \mathbf{g}_t) &\leq \|\mathbf{w}_t - \mathbf{d}_t\|\,\|\mathbf{g}_t - \mathbf{g}_{t+1}\|_\circ - \frac{t+2}{4L}\|\mathbf{g}_{t+1} - \mathbf{g}_t\|_\circ^2 \\
&\leq \frac{L\|\mathbf{w}_t - \mathbf{d}_t\|^2}{t+2}.
\end{aligned}$$

Repeating the argument we have the desired result:

$$\frac{(t+1)(t+2)}{2}[\ell^*(\bar{\mathbf{g}}_{t+1}) + f^*(-\bar{\mathbf{g}}_{t+1}) - \ell^*(\mathbf{g}) - f^*(-\mathbf{g})] \leq \sum_{s=0}^{t}\frac{s+1}{s+2}L\|\mathbf{w}_s - \mathbf{d}_s\|^2.$$

### A.8 Proof of Theorem 9

Since the subroutine is `Relaxed` and the subproblem (10) is solved up to some $\varepsilon_t \geq 0$,

$$F(\mathbf{w}_{t+1}) \leq \ell(\mathbf{w}_t) + \eta_t \langle \mathbf{d}_t - \mathbf{w}_t, \nabla\ell(\mathbf{w}_t)\rangle + \frac{L\eta_t^2}{2}\|\mathbf{d}_t - \mathbf{w}_t\|^2 + (1-\eta_t)f(\mathbf{w}_t) + \eta_t f(\mathbf{d}_t)$$

$$\leq F(\mathbf{w}_t) - \eta_t \mathsf{G}(\mathbf{w}_t) + \frac{L\eta_t^2}{2}\|\mathbf{d}_t - \mathbf{w}_t\|^2 + \eta_t \varepsilon_t$$

$$= F(\mathbf{w}_t) - \eta_t \mathsf{G}(\mathbf{w}_t) + \eta_t^2(\varepsilon_t/\eta_t + \tfrac{L}{2}\|\mathbf{d}_t - \mathbf{w}_t\|^2).$$

Define $\Delta_t := F(\mathbf{w}_t) - F(\mathbf{w})$ and $\mathsf{G}_t := \mathsf{G}(\mathbf{w}_t)$. Thus

$$\Delta_{t+1} \leq \Delta_t - \eta_t \mathsf{G}_t + \eta_t^2(\varepsilon_t/\eta_t + \tfrac{L}{2}\|\mathbf{d}_t - \mathbf{w}_t\|^2), \tag{95}$$

$$\Delta_t \leq \mathsf{G}_t. \tag{96}$$

Plug (96) into (95) and expand:

$$\Delta_{t+1} \leq \pi_t(1-\eta_0)\Delta_0 + \sum_{s=0}^{t} \frac{\pi_t}{\pi_s}\eta_s^2(\varepsilon_s/\eta_s + \tfrac{L}{2}\|\mathbf{a}_s - \mathbf{w}_s\|^2). \tag{97}$$

To prove the second claim, we have from (95)

$$\eta_t \mathsf{G}_t \leq \Delta_t - \Delta_{t+1} + \eta_t^2(\varepsilon_t/\eta_t + \tfrac{L}{2}\|dec_t - \mathbf{w}_t\|^2).$$

Summing from $k$ to $t$:

$$\left(\min_{k \leq s \leq t} \mathsf{G}_s\right)\sum_{s=k}^{t}\eta_s \leq \sum_{s=k}^{t}\eta_s \mathsf{G}_s \leq F(\mathbf{w}_k) - F(\mathbf{w}_{t+1}) + \sum_{s=k}^{t}\eta_s^2(\varepsilon_s/\eta_s + \tfrac{L}{2}\|\mathbf{d}_s - \mathbf{w}_s\|^2).$$

Rearranging completes the proof.

### A.9 Proof of Corollary 10

Since $\eta_t = 2/(t+2)$, we have $\eta_0 = 1$ and $\pi_t = \frac{2}{(t+1)(t+2)}$. For the first claim all we need to verify is that, by induction, $\frac{1}{(t+1)(t+2)}\sum_{s=0}^{t}\frac{s+1}{s+2} \leq \frac{1}{t+4}$.

We prove the second claim by a sequence of calculations similar to that in (Freund and Grigas, 2013). First, using (21) and (22) with $t = 1$ and $\mathbf{w} = \mathbf{w}_2$:

$$\tilde{\mathsf{G}}_1^1 \leq \frac{1}{\eta_1}\left(F(\mathbf{w}_1) - F(\mathbf{w}_2) + \frac{\eta_1^2}{2}(\delta + L_F)\right)$$

$$\leq \frac{3}{2}\left(\frac{1}{2}(\delta + L_F) + \frac{2}{9}(\delta + L_F)\right) = \frac{13}{12}(\delta + L_F),$$

proving (23) for $t = 1, 2$ (note that $\tilde{\mathsf{G}}_1^2 \leq \tilde{\mathsf{G}}_1^1$ by definition). For $t \geq 3$, we consider $k = t/2 - 1$ if $t$ is even and $k = (t+1)/2 - 1$ otherwise. Clearly $k \geq 1$, hence

$$\sum_{s=k}^{t}\eta_s = 2\sum_{s=k}^{t}\frac{1}{s+2} \geq 2\int_{k-1}^{t}\frac{1}{s+2}\mathrm{d}s = 2\ln\frac{t+2}{k+1} \geq 2\ln 2,$$

$$\sum_{s=k}^{t}\eta_s^2/2 = 2\sum_{s=k}^{t}\frac{1}{(s+2)^2} \leq 2\int_{k}^{t+1}\frac{1}{(s+2)^2}\mathrm{d}s = 2\left(\frac{1}{k+2} - \frac{1}{t+3}\right).$$

Using again (21), and (22) with $t = k$ and $\mathbf{w} = \mathbf{w}_{t+1}$:

$$\tilde{\mathsf{G}}_1^t \leq \tilde{\mathsf{G}}_k^t \leq \frac{\delta + L_F}{2\ln 2}\left(\frac{2}{k+3} + \frac{2}{k+2} - \frac{2}{t+3}\right) \leq \frac{\delta + L_F}{\ln 2}\left(\frac{2}{t+2} + \frac{1}{t+3}\right) \leq \frac{3(\delta + L_F)}{t\ln 2}.$$

### A.10 Proof of Theorem 12

A proof can be based on the simple observation that the objective value $F^\star$ satisfies

$$F^\star := \inf_{\mathbf{w}} \{\ell(\mathbf{w}) + f(\mathbf{w})\} = \inf_{\mathbf{w}, \rho : \kappa(\mathbf{w}) \leq \rho} \ell(\mathbf{w}) + h(\rho). \tag{98}$$

Note that if $\rho$ were known, the theorem could have been proved as before. The key idea behind the step size choice (41) is to ensure that Algorithm 2 performs almost the same as if it knew the *unknown* but fixed constant $\rho = \kappa(\mathbf{w})$ beforehand.

Consider an arbitrary $\mathbf{w}$ and let $\rho = \kappa(\mathbf{w})$. We also use the shorthand $\hat{F}_t := \ell(\mathbf{w}_t) + h(\rho_t) \geq \ell(\mathbf{w}_t) + f(\mathbf{w}_t) = F(\mathbf{w}_t)$. Then, the following chain of inequalities can be verified:

$$
\begin{aligned}
\hat{F}_{t+1} &:= \ell(\mathbf{w}_{t+1}) + h(\rho_{t+1}) \\
&\leq \hat{F}_t + \langle \theta_t \mathbf{a}_t - \eta_t \mathbf{w}_t, \nabla\ell(\mathbf{w}_t) \rangle + \tfrac{L}{2} \|\theta_t \mathbf{a}_t - \eta_t \mathbf{w}_t\|^2 - \eta_t h(\rho_t) + \eta_t h(\theta_t/\eta_t) \\
&\leq \hat{F}_t + \eta_t \left\langle \tfrac{\rho}{\alpha_t} \mathbf{a}_t - \mathbf{w}_t, \nabla\ell(\mathbf{w}_t) \right\rangle + \tfrac{L\eta_t^2}{2} \left\| \tfrac{\rho}{\alpha_t} \mathbf{a}_t - \mathbf{w}_t \right\|^2 - \eta_t h(\rho_t) + \eta_t h(\rho/\alpha_t) \\
&\leq \hat{F}_t + \min_{\mathbf{z} : \kappa(\mathbf{z}) \leq \rho} \eta_t \langle \mathbf{z} - \mathbf{w}_t, \nabla\ell(\mathbf{w}_t) \rangle + \eta_t \rho \varepsilon_t + \tfrac{L\eta_t^2}{2} \left\| \tfrac{\rho}{\alpha_t} \mathbf{a}_t - \mathbf{w}_t \right\|^2 - \eta_t h(\rho_t) + \eta_t h(\rho/\alpha_t) \\
&= \hat{F}_t + \min_{\mathbf{z} : \kappa(\mathbf{z}) \leq \rho} \eta_t \langle \mathbf{z} - \mathbf{w}_t, \nabla\ell(\mathbf{w}_t) \rangle - \eta_t h(\rho_t) + \eta_t h(\rho) \\
&\qquad\quad + \eta_t^2 \underbrace{\left( \tfrac{L}{2} \left\| \tfrac{\rho}{\alpha_t} \mathbf{a}_t - \mathbf{w}_t \right\|^2 + (\rho \varepsilon_t + h(\rho/\alpha_t) - h(\rho))/\eta_t \right)}_{:= \delta_t} \\
&= \hat{F}_t + \eta_t^2 \delta_t - \eta_t \underbrace{\left[ \langle \mathbf{w}_t, \nabla\ell(\mathbf{w}_t) \rangle + h(\rho_t) - \min_{\mathbf{z} : \kappa(\mathbf{z}) \leq \rho} \langle \mathbf{z}, \nabla\ell(\mathbf{w}_t) \rangle + h(\rho) \right]}_{:= \hat{\mathsf{G}}(\mathbf{w}_t)},
\end{aligned}
$$

where the first inequality is because the subroutine is `Relaxed`, the second inequality follows from the minimality of $\theta_t$ in (41), and the third inequality is due to the choice of $\mathbf{a}_t$ in Line 3 of Algorithm 2.

Recall that $\rho = \kappa(\mathbf{w})$. Moreover, due to the convexity of $\ell$,

$$
\begin{aligned}
\hat{\mathsf{G}}(\mathbf{w}_t) &= \langle \mathbf{w}_t, \nabla\ell(\mathbf{w}_t) \rangle + h(\rho_t) - \min_{\mathbf{z} : \kappa(\mathbf{z}) \leq \rho} \left( \langle \mathbf{z}, \nabla\ell(\mathbf{w}_t) \rangle + h(\rho) \right) \\
&= h(\rho_t) - h(\kappa(\mathbf{w})) + \max_{\mathbf{z} : \kappa(\mathbf{z}) \leq \kappa(\mathbf{w})} \langle \mathbf{w}_t - \mathbf{z}, \nabla\ell(\mathbf{w}_t) \rangle \\
&\geq h(\rho_t) - h(\kappa(\mathbf{w})) + \max_{\mathbf{z} : \kappa(\mathbf{z}) \leq \kappa(\mathbf{w})} \ell(\mathbf{w}_t) - \ell(\mathbf{z}) \\
&\geq \hat{F}_t - F(\mathbf{w}).
\end{aligned}
$$

Thus we have retrieved the recursion:

$$
\begin{aligned}
\hat{F}_{t+1} - F(\mathbf{w}) &\leq \hat{F}_t - F(\mathbf{w}) - \eta_t \hat{\mathsf{G}}(\mathbf{w}_t) + \eta_t^2 \delta_t, \\
\hat{F}_t - F(\mathbf{w}) &\leq \hat{\mathsf{G}}(\mathbf{w}_t).
\end{aligned}
$$

Summing the indices as in the proof of Theorem 9 and noting that $F(\mathbf{w}_t) \leq \hat{F}_t$ for all $t$ completes the proof.

## Appendix B. Proofs for Section 4

### B.1 Proof of Proposition 14

We first note that the atomic set

$$\mathcal{A} := \{\mathbf{u}\mathbf{v}^\top : \mathbf{u} \in \mathbb{R}^m, \mathbf{v} \in \mathbb{R}^n, \|\mathbf{u}\|_c \leq 1, \|\mathbf{v}\|_r \leq 1\}$$

is compact, so is its convex hull conv $\mathcal{A}$. It is also easy to see that $\mathcal{A}$ is a connected set (by considering the continuous map $(\mathbf{u}, \mathbf{v}) \mapsto \mathbf{u}\mathbf{v}^\top$).

Recall from (30) that

$$\kappa(W) = \inf\left\{\rho \geq 0 : W = \rho \sum_i \sigma_i \mathbf{a}_i, \sigma_i \geq 0, \sum_i \sigma_i = 1, \mathbf{a}_i \in \mathcal{A}\right\} \tag{99}$$

$$= \inf\left\{\rho \geq 0 : W \in \rho \text{ conv } \mathcal{A}\right\}.$$

Since conv $\mathcal{A}$ is compact with $\mathbf{0} \in \text{int}(\text{conv } \mathcal{A})$ we know the infimum above is attained. Thus there exist $\rho \geq 0$ and $C \in \text{conv } \mathcal{A}$ so that $W = \rho C$ and $\kappa(W) = \rho$. Applying Caratheodory's theorem (for connected sets) we know $C = \sum_{i=1}^t \sigma_i \mathbf{u}_i \mathbf{v}_i^\top$ for some $\sigma_i \geq 0$, $\sum_i \sigma_i = 1$, $\mathbf{u}_i \mathbf{v}_i^\top \in \mathcal{A}$ and $t \leq mn$. Let $U = \sqrt{\rho}\left[\sqrt{\sigma_1}\mathbf{u}_1, \ldots, \sqrt{\sigma_t}\mathbf{u}_t\right]$, and $V = \sqrt{\rho}\left[\sqrt{\sigma_1}\mathbf{v}_1, \ldots, \sqrt{\sigma_t}\mathbf{v}_t\right]^\top$. We then have $W = UV$ and

$$\kappa(W) = \rho \geq \frac{1}{2}\sum_{i=1}^t \left(\|U_{:i}\|_c^2 + \|V_{i:}\|_r^2\right),$$

which proves one side of Proposition 14.

On the other hand, consider any $U$ and $V$ that satisfy

$$W = UV = \sum_{j=1}^t \|U_{:j}\|_c \|V_{j:}\|_r \sum_{i=1}^t \frac{\|U_{:i}\|_c \|V_{i:}\|_r}{\sum_{j=1}^t \|U_{:j}\|_c \|V_{j:}\|_r} \frac{U_{:i}}{\|U_{:i}\|_c} \frac{V_{i:}}{\|V_{i:}\|_r},$$

assuming w.l.o.g. that $\|U_{:i}\|_c \|V_{i:}\|_r \neq 0$ for all $1 \leq i \leq t$. This, together with the definition (99), gives us the other half of the equality:

$$\kappa(W) \leq \sum_{i=1}^t \|U_{:i}\|_c \|V_{i:}\|_r.$$

The proof is completed by making the elementary observation

$$\sum_{i=1}^t \|U_{:i}\|_c \|V_{i:}\|_r \leq \frac{1}{2}\sum_{i=1}^t \left(\|U_{:i}\|_c^2 + \|V_{i:}\|_r^2\right).$$

### B.2 Proof of Proposition 16

By (59), the square of the polar is

$$(\kappa^\circ(G))^2 = \max\left\{\mathbf{c}^\top Z Z^\top \mathbf{c} : \|\mathbf{a}\|_2 \leq 1, \|\mathbf{b}\|_2 \leq 1\right\}, \tag{100}$$

where recall that $\mathbf{c} = \begin{bmatrix} \mathbf{a} \\ \mathbf{b} \end{bmatrix}$ is the concatenation of $\mathbf{a}$ and $\mathbf{b}$. Since we are maximizing a convex quadratic function over a convex set in (100), there must be a maximizer attained at an extreme point of the feasible region. Therefore the problem is equivalent to

$$(\kappa^\diamond(G))^2 = \max\left\{\mathbf{c}^\top Z Z^\top \mathbf{c} : \|\mathbf{a}\|_2 = 1, \|\mathbf{b}\|_2 = 1\right\} \tag{101}$$

$$= \max\left\{\operatorname{tr}(SZZ^\top) : \operatorname{tr}(SI_1) = 1, \operatorname{tr}(SI_2) = 1, S \succeq 0, \operatorname{rank}(S) = 1\right\}, \tag{102}$$

where $S = \mathbf{c}\mathbf{c}^\top$. Key to this proof is the fact that the rank one constraint can be dropped without breaking the equivalence. To see why, notice that without this rank constraint the feasible region is the intersection of the positive semi-definite cone with two hyperplanes, hence the rank of all its extreme points must be upper bounded by one (Pataki, 1998). Furthermore the linearity of the objective implies that there must be a maximizing solution attained at an extreme point of the feasible region. Therefore we can continue by

$$(\kappa^\diamond(G))^2 = \max\left\{\operatorname{tr}(SZZ^\top) : \operatorname{tr}(SI_1) = 1, \operatorname{tr}(SI_2) = 1, S \succeq 0\right\}$$

$$= \max_{S \succeq 0} \min_{\mu_1, \mu_2 \in \mathbb{R}} \operatorname{tr}(SGG^\top) - \mu_1(\operatorname{tr}(SI_1) - 1) - \mu_2(\operatorname{tr}(SI_2) - 1)$$

$$= \min_{\mu_1, \mu_2 \in \mathbb{R}} \max_{S \succeq 0} \operatorname{tr}(SGG^\top) - \mu_1(\operatorname{tr}(SI_1) - 1) - \mu_2(\operatorname{tr}(SI_2) - 1)$$

$$= \min\left\{\mu_1 + \mu_2 : \mu_1, \mu_2 \in \mathbb{R}, GG^\top \preceq \mu_1 I_1 + \mu_2 I_2\right\} \tag{103}$$

$$= \min\left\{\mu_1 + \mu_2 : \mu_1 \geq 0, \mu_2 \geq 0, \left\|D_{\mu_2/\mu_1}G\right\|_{\mathrm{sp}}^2 \leq \mu_1 + \mu_2\right\} \tag{104}$$

$$= \min\left\{\left\|D_\mu G\right\|_{\mathrm{sp}}^2 : \mu \geq 0\right\},$$

where the third equality is due to the Lagrangian duality, the fifth equality follows from the equivalence $\left\|D_{\mu_2/\mu_1}G\right\|_{\mathrm{sp}}^2 \leq \mu_1 + \mu_2 \iff D_{\mu_2/\mu_1}GG^\top D_{\mu_2/\mu_1} \preceq (\mu_1 + \mu_2)I \iff GG^\top \preceq (\mu_1 + \mu_2)D_{\mu_2/\mu_1}^{-1}D_{\mu_2/\mu_1}^{-1} = \mu_1 I_1 + \mu_2 I_2$, and the last equality is obtained through the re-parameterization $\mu = \mu_2/\mu_1, \nu = \mu_1 + \mu_2$ and the elimination of $\nu$.

### B.3 Solving the multi-view polar

This section shows how to efficiently compute the polar (dual norm) in the multi-view setting.

Let us first backtrack in (104), which produces the optimal $\mu_1$ and $\mu_2$ via

$$\mu_1 + \mu_2 = \|D_\mu G\|_{\mathrm{sp}}^2, \qquad \text{and} \qquad \frac{\mu_2}{\mu_1} = \mu. \tag{105}$$

Then the KKT condition for the step from (103) to (104) can be written as

$$M \succeq \mathbf{0}, \qquad M\mathbf{c} = \mathbf{0}, \qquad \|\mathbf{a}\|_2 = \|\mathbf{b}\|_2 = 1, \quad \text{where} \quad M := \mu_1 I_1 + \mu_2 I_2 - GG^\top. \tag{106}$$

By duality, this is a sufficient and necessary condition for $\mathbf{a}$ and $\mathbf{b}$ to be the solution of the polar operator. Given the optimal $\mu_1$ and $\mu_2$, the first condition $M \succeq \mathbf{0}$ must have been satisfied already. Let the null space of $M$ be spanned by an orthonormal basis $\{\mathbf{q}_1, \ldots, \mathbf{q}_k\}$. Then we can parameterize $\mathbf{c}$ as $\mathbf{c} = Q\boldsymbol{\alpha}$ where $Q = [\mathbf{q}_1, \ldots, \mathbf{q}_k] =: \begin{bmatrix} Q_1 \\ Q_2 \end{bmatrix}$. By (106),

$$0 = \|\mathbf{a}\|_2^2 - \|\mathbf{b}\|_2^2 = \boldsymbol{\alpha}' \left( Q_1^\top Q_1 - Q_2^\top Q_2 \right) \boldsymbol{\alpha}. \tag{107}$$

Let $P\Sigma P^\top$ be the eigen-decomposition of $Q_1^\top Q_1 - Q_2^\top Q_2$ with $\Sigma$ being diagonal and $P$ being unitary. Let $\mathbf{v} = P^\top \boldsymbol{\alpha}$. Then $\mathbf{c} = QP\mathbf{v}$ and (107) simply becomes $\mathbf{v}^\top \Sigma \mathbf{v} = 0$. In addition, (106) also implies $2 = \|\mathbf{c}\|_2^2 = \mathbf{v}^\top P^\top Q^\top QP\mathbf{v} = \mathbf{v}^\top \mathbf{v}$. So finally the optimal solution to the polar problem (101) can be completely characterized by

$$\left\{ \begin{bmatrix} \mathbf{a} \\ \mathbf{b} \end{bmatrix} = QP\mathbf{v} : \mathbf{v}^\top \Sigma \mathbf{v} = 0, \|\mathbf{v}\|_2^2 = 2 \right\}, \tag{108}$$

which is simply a linear system after a straightforward change of variable. The major computational cost for recovery is to find the null space of $M$, which can be achieved by QR decomposition. The complexity of eigen-decomposition for $Q_1^\top Q_1 - Q_2^\top Q_2$ depends on the dimension of the null space of $M$, which is usually low in practice.

## Appendix C. Proofs for Section 5

### C.1 Vertices of $Q$ must be scalar multiples of those of $P$

First note that if $\mathbf{0} \notin P$, we have nothing to prove since $P = Q$. So we assume $\mathbf{0} \in P$ below.

Consider an arbitrary vertex $\mathbf{q} \in Q$. Clearly $\mathbf{q} \neq \mathbf{0}$ and $\mathbf{q} \in P$, hence $\mathbf{q} = \sum_{i=1}^n \alpha_i \cdot \mathbf{p}^{(i)}$, where $n \geq 1$, $\alpha_i > 0$, $\langle \mathbf{1}, \boldsymbol{\alpha} \rangle \leq 1$, and $\mathbf{p}^{(i)}$ are nonzero vertices of $P$. Clearly $\mathbf{p}^{(i)} \in Q$ as $\mathbf{p}^{(i)} \in P$ and $l_i := \langle \mathbf{1}, \mathbf{p}^{(i)} \rangle \geq 1$. It suffices to show $n = 1$. To prove by contradiction, suppose $n \geq 2$.

**(i)** If $\langle \mathbf{1}, \boldsymbol{\alpha} \rangle = 1$, then $\mathbf{q}$ is a convex combination of at least two points in $Q$, hence it cannot be a vertex.

**(ii)** If $\langle \mathbf{1}, \mathbf{q} \rangle = \sum_i \alpha_i l_i = 1$, then $\mathbf{q} = \sum_{i=1}^n (\alpha_i l_i) \frac{\mathbf{p}^{(i)}}{l_i}$. But $\frac{\mathbf{p}^{(i)}}{l_i} \in Q$ as $\frac{\mathbf{p}^{(i)}}{l_i} = \frac{1}{l_i}\mathbf{p}^{(i)} + (1 - \frac{1}{l_i})\mathbf{0} \in P$ and $\langle \mathbf{1}, \mathbf{p}^{(i)} \rangle = l_i \geq 1$. Again contradiction.

**(iii)** If $\langle \mathbf{1}, \mathbf{q} \rangle > 1$ and $\langle \mathbf{1}, \boldsymbol{\alpha} \rangle < 1$, then $\beta := \frac{1}{\langle \mathbf{1}, \mathbf{q} \rangle} < 1 < \frac{1}{\langle \mathbf{1}, \boldsymbol{\alpha} \rangle} =: \gamma$. Clearly $\beta\mathbf{q} \in Q$ because $\beta\mathbf{q} = \beta\mathbf{q} + (1 - \beta)\mathbf{0} \in P$ and $\langle \mathbf{1}, \beta\mathbf{q} \rangle = 1$. Also $\gamma\mathbf{q} \in Q$ as $\gamma\mathbf{q} = \frac{\sum_{i=1}^n \alpha_i \mathbf{p}^{(i)}}{\sum_{i=1}^n \alpha_i} \in P$ and $\langle \mathbf{1}, \gamma\mathbf{q} \rangle = \frac{\sum_{i=1}^n \alpha_i \langle \mathbf{1}, \mathbf{p}^{(i)} \rangle}{\langle \mathbf{1}, \boldsymbol{\alpha} \rangle} \geq \frac{\sum_{i=1}^n \alpha_i}{\langle \mathbf{1}, \boldsymbol{\alpha} \rangle} = 1$. So $\mathbf{q}$ lies between two points in $Q$: $\beta\mathbf{q}$ and $\gamma\mathbf{q}$. Contradiction.

Therefore $n = 1$, which completes the proof.

To summarize, we have proved that if $\mathbf{q}$ is a vertex of $Q$ but not a vertex of $P$, then it must sum to 1 and be a scalar multiple of some vertex of $P$.

### C.2 Proof of Proposition 17

The proof is based on Proposition 2. We assume that $r$, the upper bound on the number of groups each variable can belong to, is greater than 1 since otherwise the problem is trivial.

*Proof.* Note that there are $n$ variables which we index by $i$ and there are $\ell$ groups (subsets of variables) which we index by $G$. The input vector $\tilde{\mathbf{w}} \in \mathbb{R}^n \times \mathbb{R}^\ell$.

Let $l_i$ be the number of groups that contain variable $i$, and $\mathcal{S}_i := \{\mathbf{s} \in \mathbb{R}_+^{l_i} : \langle \mathbf{1}, \mathbf{s} \rangle = 1\}$ be the $(l_i - 1)$-dimensional simplex. Using the well-known variational representation of the max function, we rewrite the (negated) objective $h$ in (87) as

$$h(\tilde{\mathbf{w}}) = \sum_{i \in [n]} \tilde{g}_i \max_{\boldsymbol{\alpha}^{(i)} \in \mathcal{S}_i} \left\{ -\sum_{G : i \in G} \alpha_G^{(i)} \tilde{w}_G \right\} = \max_{\boldsymbol{\alpha}^{(i)} \in \mathcal{S}_i} \sum_{i \in [n]} \sum_{G : i \in G} -\tilde{g}_i \tilde{w}_G \alpha_G^{(i)}, \tag{109}$$

which is to be minimized. Here the second equality follows from the separability of the variables $\boldsymbol{\alpha}^{(i)}$. Fix $\epsilon > 0$ and denote $c := \frac{\epsilon}{n \log r}$. Consider

$$h_\epsilon(\tilde{\mathbf{w}}) = \max_{\boldsymbol{\alpha}^{(i)} \in \mathcal{S}_i} \sum_{i \in [n]} \sum_{G : i \in G} \left( -\tilde{g}_i \tilde{w}_G \alpha_G^{(i)} - c \cdot \alpha_G^{(i)} \log \alpha_G^{(i)} \right),$$

i.e., we add to $h$ the scaled entropy function $-c \sum_{i \in [n], G : i \in G} \alpha_G^{(i)} \log \alpha_G^{(i)}$ whose negation is known to be strongly convex on the simplex (w.r.t. the $\ell_1$-norm) (Nesterov, 2005). Since the entropy is nonnegative, we have for any $\tilde{\mathbf{w}}$, $h(\tilde{\mathbf{w}}) \leq h_\epsilon(\tilde{\mathbf{w}})$ and moreover

$$h_\epsilon(\tilde{\mathbf{w}}) - h(\tilde{\mathbf{w}}) \leq c \max_{\boldsymbol{\alpha}^{(i)} \in \mathcal{S}_i} \sum_{i \in [n]} \sum_{G : i \in G} -\alpha_G^{(i)} \log \alpha_G^{(i)} \leq c \cdot n \log r = \epsilon,$$

where the last inequality is due to the well-known upper bound of the entropy over the probability simplex, i.e. entropy attains its maximum when all odds are equally likely. Therefore $h(\tilde{\mathbf{w}}) - h_\epsilon(\tilde{\mathbf{w}}) \in (-\epsilon, 0]$, and we have proved part (ii) of Proposition 17.

By straightforward calculation

$$h_\epsilon(\tilde{\mathbf{w}}) = \sum_{i \in [n]} \max_{\boldsymbol{\alpha}^{(i)} \in \mathcal{S}_i} \sum_{G : i \in G} \left( -\tilde{g}_i \tilde{w}_G \alpha_G^{(i)} - c \cdot \alpha_G^{(i)} \log \alpha_G^{(i)} \right)$$

$$= c \sum_{i \in [n]} \log \sum_{G : i \in G} \exp \left( -\frac{\tilde{g}_i \tilde{w}_G}{c} \right), \tag{110}$$

$$\frac{\partial}{\partial \tilde{w}_G} h_\epsilon(\tilde{\mathbf{w}}) = -\sum_{i : i \in G} \tilde{g}_i p_i(G), \quad \text{where} \quad p_i(G) := \frac{\exp\left( -\frac{\tilde{g}_i \tilde{w}_G}{c} \right)}{\sum_{\tilde{G} : i \in \tilde{G}} \exp\left( -\frac{\tilde{g}_i \tilde{w}_{\tilde{G}}}{c} \right)}. \tag{111}$$

Hence $h_\epsilon(\tilde{\mathbf{w}})$ can be computed in $O(nr)$ time (since the second summation in (110) contains at most $r$ terms). Similarly all $\{p_i(G) : i \in [n], i \in G\}$ can be computed in $O(nr)$ time. Therefore part (iii) of Proposition 17 is established.

Finally, to bound the Lipschitz constant of the gradient of $h_\epsilon$, we observe that $h_\epsilon(\tilde{\mathbf{w}}) = \eta^*(A\tilde{\mathbf{w}})$, where $\eta^*$ is the Fenchel conjugate of the scaled negative entropy

$$\eta(\boldsymbol{\alpha}) = c \sum_{i \in [n]} \sum_{G : i \in G} \alpha_G^{(i)} \log \alpha_G^{(i)},$$

and $A$ is defined as the matrix satisfying

$$\langle \boldsymbol{\alpha}, A\tilde{\mathbf{w}} \rangle = \sum_{i \in [n]} \sum_{G : i \in G} -\alpha_G^{(i)} \tilde{g}_i \tilde{w}_G.$$

It is known that the scaled negative entropy $\eta$ is strongly convex with modulus $c$ (w.r.t. the $\ell_1$-norm). Furthermore, employing $\ell_1$ norm on $\boldsymbol{\alpha}$ and $\ell_2$ norm on $\tilde{\mathbf{w}}$, the operator norm of the matrix $A$ can be bounded as

$$
\begin{aligned}
\|A\|_{2,1} &:= \max_{\boldsymbol{\alpha}: \|\boldsymbol{\alpha}\|_1 = 1} \max_{\tilde{\mathbf{w}}: \|\tilde{\mathbf{w}}\|_2 = 1} \langle \boldsymbol{\alpha}, A\tilde{\mathbf{w}} \rangle = \max_{\tilde{\mathbf{w}}: \|\tilde{\mathbf{w}}\|_2 = 1} \max_{\boldsymbol{\alpha}: \|\boldsymbol{\alpha}\|_1 = 1} \sum_{i \in [n]} \sum_{G: i \in G} -\alpha_G^{(i)} \tilde{g}_i \tilde{w}_G \\
&\leq \left( \max_{i \in [n]} \tilde{g}_i \right) \cdot \max_{\tilde{\mathbf{w}} \geq 0: \|\tilde{\mathbf{w}}\|_2 = 1} \max_{\boldsymbol{\alpha} \geq 0: \|\boldsymbol{\alpha}\|_1 = 1} \sum_{i \in [n]} \sum_{G: i \in G} \alpha_G^{(i)} \tilde{w}_G \\
&\leq \left( \max_{i \in [n]} \tilde{g}_i \right) \cdot \max_{\boldsymbol{\alpha} \geq 0: \|\boldsymbol{\alpha}\|_1 = 1} \sum_{i \in [n]} \sum_{G: i \in G} \alpha_G^{(i)} = \max_{i \in [n]} \tilde{g}_i = \|\tilde{\mathbf{g}}\|_\infty .
\end{aligned}
$$

The equality is obviously attainable. Therefore by Theorem 1 of (Nesterov, 2005), $h_\epsilon(\tilde{\mathbf{w}}) = \eta^*(A\tilde{\mathbf{w}})$ has Lipschitz continuous gradient w.r.t. the $\ell_2$ norm, and the Lipschitz constant is $\frac{1}{c} \|A\|_{2,1}^2 = \frac{1}{\epsilon} \|\tilde{\mathbf{g}}\|_\infty^2 \, n \log r$. This completes our proof of part (i) of Proposition 17. $\qquad\square$

### C.3 Recovery of An Integral Solution

This section discusses the subtle issue of recovering an (approximate) integral solution of the polar operator (80) from that of the smoothed LP reformulation.

Recall that our ultimate goal is to find an integral solution of the polar operator (80):

$$
\lambda^* := \max_{\mathbf{0} \neq \mathbf{w} \in P} \frac{\langle \tilde{\mathbf{g}}, \mathbf{w} \rangle}{\langle \mathbf{b}, \mathbf{w} \rangle}.
$$

As we showed in Section 5.2, the optimal objective value is exactly equal to that of (82), for which the smoothing technique is able to find an $\epsilon$ accurate solution efficiently. That means we can obtain some $\lambda_\epsilon$ (smooth objective function value) with the guarantee that $\lambda_\epsilon \in [\lambda^* - \epsilon, \lambda^*]$. With such $\lambda_\epsilon$ in hand, we now show how to find an $\epsilon$ accurate solution for (80), i.e. a $\mathbf{w}_\epsilon \in P \backslash \{\mathbf{0}\}$ such that

$$
\frac{\langle \tilde{\mathbf{g}}, \mathbf{w}_\epsilon \rangle}{\langle \mathbf{b}, \mathbf{w}_\epsilon \rangle} \geq \lambda^* - \epsilon. \tag{112}
$$

Indeed, this is simple according to the following proposition.

**Proposition 18.** *Given $\lambda_\epsilon \in [\lambda^* - \epsilon, \lambda^*]$, find*

$$
\mathbf{w}_\epsilon := \arg \max_{\mathbf{w} \in P \backslash \{\mathbf{0}\}} \left\{ \langle \tilde{\mathbf{g}}, \mathbf{w} \rangle - \lambda_\epsilon \langle \mathbf{b}, \mathbf{w} \rangle \right\}. \tag{113}
$$

*Then $\mathbf{w}_\epsilon$ satisfies (112).*

*Proof.* By the definition of $\lambda^*$, $\max_{\mathbf{w} \in P \backslash \{\mathbf{0}\}} \{ \langle \tilde{\mathbf{g}}, \mathbf{w} \rangle - \lambda^* \langle \mathbf{b}, \mathbf{w} \rangle \} = 0$. As $\lambda_\epsilon \leq \lambda^*$, so $\max_{\mathbf{w} \in P \backslash \{\mathbf{0}\}} \{ \langle \tilde{\mathbf{g}}, \mathbf{w} \rangle - \lambda_\epsilon \langle \mathbf{b}, \mathbf{w} \rangle \} \geq 0$. This implies $\frac{\langle \tilde{\mathbf{g}}, \mathbf{w}_\epsilon \rangle}{\langle \mathbf{b}, \mathbf{w}_\epsilon \rangle} \geq \lambda_\epsilon \geq \lambda^* - \epsilon$. $\qquad\square$

Note that (113) is exactly the submodular minimization problem that the secant method of (Obozinski and Bach, 2012) is based on. This step is computationally expensive and has to be solved for multiple values of $\lambda_t$ in the secant method. By contrast, our strategy needs to solve this problem only once.

For group sparsity, (113) boils down to a max-flow problem which is again expensive. Fortunately, by exploiting the structure of the problem it is possible to design a heuristic rounding procedure. For convenience, let us copy (87), the LP reformulation, here:

$$\max_{\tilde{\mathbf{w}}} \ \sum_{i \in [n]} \tilde{g}_i \min_{G : i \in G \in \mathcal{G}} \tilde{w}_G, \quad \text{subject to} \quad \tilde{\mathbf{w}} \geq 0, \ \sum_{G \in \mathcal{G}} b_G \tilde{w}_G = 1.$$

A solution $\tilde{\mathbf{w}}$ corresponds to an integral solution to the polar operator if and only if $\tilde{w}_G \in \{0, c\}$ where $c$ ensures the normalization $\sum_{G \in \mathcal{G}} b_G \tilde{w}_G = 1$. By solving the smoothed objective, we obtain a solution $\tilde{\mathbf{w}}^*$ which does not necessarily satisfy this condition. However, a smaller value of the component $\tilde{w}_G^*$ does suggest a higher likelihood for $\tilde{w}_G$ to be 0. Therefore, we can sort $\{w_G^*\}$ and set the $w_G$ of the smallest $k$ groups to 0 ($k$ ranging from 0 to $|\mathcal{G}| - 1$), and the $w_G$ for the remaining groups were set to a common value that fulfills the normalization constraint. The best $k$ can be selected by brute-force enumeration. Moreover, by exploiting the structure of the objective, we can design an algorithm that accomplishes the enumeration in $O(nr)$ time. Finally, we can check the suboptimality of resulting solution: If we do not decrease the objective value by say $\epsilon$, then this rounded integral solution is $2\epsilon$ sub-optimal for the polar operator (80); otherwise we fall back to (113), which we observe rarely happens in practice.

## Appendix D. Efficient Hessian-vector Multiplication for Logistic Loss

Recall in multi-class classification, each training example $\mathbf{x}_i$ has a label $y_i \in [C]$. Denote $X = (\mathbf{x}_1, \ldots, \mathbf{x}_m) \in \mathbb{R}^{n \times m}$, and $W \in \mathbb{R}^{n \times C}$ is the weight matrix where each column corresponds to a class. Then the average logistic loss can be written as

$$-\ell(W) = \frac{1}{m} \sum_{i=1}^{m} \log \frac{\exp\left(\mathbf{x}_i^\top W_{:,y_i}\right)}{\sum_{c=1}^{C} \exp\left(\mathbf{x}_i^\top W_{:,c}\right)} = \frac{1}{m} \sum_{i=1}^{m} \log \frac{\exp\left(\langle \mathbf{x}_i \otimes \mathbf{1}_{y_i}^\top, W \rangle\right)}{\sum_c \exp\left(\langle \mathbf{x}_i \otimes \mathbf{1}_c^\top, W \rangle\right)},$$

where $\langle A, B \rangle = \mathrm{tr}(A^\top B)$, and $\otimes$ is the Kronecker product. The gradient is

$$-\nabla \ell(W) = \frac{1}{m} \sum_{i=1}^{m} \left[ \mathbf{x}_i \otimes \mathbf{1}_{y_i}^\top - \sum_c P_{ci} \left( \mathbf{x}_i \otimes \mathbf{1}_c^\top \right) \right], \quad \text{where } P_{ci} := \frac{\exp\left(\langle \mathbf{x}_i \otimes \mathbf{1}_c^\top, W \rangle\right)}{\sum_d \exp\left(\langle \mathbf{x}_i \otimes \mathbf{1}_d^\top, W \rangle\right)}.$$

Denote the probability matrix $P = (P_{ci})_{C \times n}$. The Hessian in $W$ is a four mode tensor.

To simplify notation, we adopt the vector representation of $W$ and denote $\mathbf{w} = \text{vec}(W)$. Then the loss, gradient, and Hessian can be written as

$$-\ell(\mathbf{w}) = \frac{1}{m} \sum_{i=1}^{m} \log \frac{\exp\left(\langle \mathbf{1}_{y_i} \otimes \mathbf{x}_i, \mathbf{w} \rangle\right)}{\sum_c \exp\left(\langle \mathbf{1}_c \otimes \mathbf{x}_i, \mathbf{w} \rangle\right)},$$

$$-\nabla\ell(\mathbf{w}) = \frac{1}{m} \sum_{i=1}^{m} \left[ \mathbf{1}_{y_i} \otimes \mathbf{x}_i - \sum_c P_{ci} \left(\mathbf{1}_c \otimes \mathbf{x}_i\right) \right], \text{ where } P_{ci} := \frac{\exp\left(\langle \mathbf{1}_c \otimes \mathbf{x}_i, \mathbf{w} \rangle\right)}{\sum_d \exp\left(\langle \mathbf{1}_d \otimes \mathbf{x}_i, \mathbf{w} \rangle\right)},$$

$$-\nabla^2\ell(\mathbf{w}) = \frac{1}{m} \sum_{i=1}^{m} \left[ \sum_c P_{ci}(\mathbf{1}_c \otimes \mathbf{x}_i)(\mathbf{1}_c \otimes \mathbf{x}_i)^\top - \sum_{c,d} P_{ci}P_{di}(\mathbf{1}_c \otimes \mathbf{x}_i)(\mathbf{1}_d \otimes \mathbf{x}_i)^\top \right]$$

$$= \frac{1}{m} \sum_{i=1}^{m} \left[ \sum_c P_{ci}\left((\mathbf{1}_c \otimes \mathbf{1}_c^\top) \otimes \mathbf{x}_i \otimes \mathbf{x}_i^\top\right) - \sum_{c,d} P_{ci}P_{di}\left((\mathbf{1}_c \otimes \mathbf{1}_d^\top) \otimes \mathbf{x}_i \otimes \mathbf{x}_i^\top\right) \right]$$

$$= \frac{1}{m} \sum_{i=1}^{m} Q^{(i)} \otimes \mathbf{x}_i \otimes \mathbf{x}_i^\top, \text{ where } Q^{(i)} := \sum_c P_{ci}(\mathbf{1}_c \otimes \mathbf{1}_c^\top) - \sum_{c,d} P_{ci}P_{di}(\mathbf{1}_c \otimes \mathbf{1}_d^\top).$$

Clearly $Q^{(i)}$ is symmetric. Given a direction $D \in \mathbb{R}^{n \times C}$, the directional Hessian is given by

$$-H_D := \nabla^2\ell(\mathbf{w}) \cdot \text{vec}(D) = \frac{1}{m} \sum_{i=1}^{m} \left[ Q^{(i)} \otimes \mathbf{x}_i \otimes \mathbf{x}_i^\top \right] \text{vec}(D) = \frac{1}{m} \sum_{i=1}^{m} \text{vec}\left(\mathbf{x}_i^\top D(Q^{(i)} \otimes \mathbf{x}_i^\top)\right).$$

Let $R := D^\top X \in \mathbb{R}^{C \times m}$, we can continue by

$$-m\,H_D = \sum_{i=1}^{m} \text{vec}\left(R_{:i}^\top(Q^{(i)} \otimes \mathbf{x}_i^\top)\right) = \sum_{i=1}^{m} \text{vec}\left((Q^{(i)} \otimes \mathbf{x}_i)R_{:i}\right) = \sum_{i=1}^{m} \text{vec}\left((Q^{(i)}R_{:i}) \otimes \mathbf{x}_i\right)$$

$$= \sum_{i=1}^{m} \text{vec}(J_{:i} \otimes \mathbf{x}_i) = \text{vec}(XJ^\top), \quad \text{where} \quad J := [Q^{(1)}R_{:1}, \ldots, Q^{(m)}R_{:m}]. \tag{114}$$

Finally we work out $J$. Note that the $c$-th element of $Q^{(i)}R_{:i}$ is $P_{ci}R_{ci} - P_{ci}\sum_d P_{di}R_{di}$. So

$$J = S - P \circ (\mathbf{1} \otimes (\mathbf{1}^\top S)), \quad \text{where} \quad S := P \circ R \text{ and } \circ \text{ denotes Hadamard product.} \tag{115}$$

In summary, we compute (a) $R = D^\top X$, (b) $S = P \circ R$, (c) $J$ as in (115), and (d) use (114).